\documentclass[10pt]{article}
\usepackage[latin1]{inputenc}
\usepackage{amssymb, amsmath, amsthm}
\usepackage{makeidx}
\usepackage[totalwidth=450pt, totalheight=650pt]{geometry}

\newtheorem{ex}{Example}[subsection]
\newtheorem{theo}[ex]{Theorem}

\newtheorem{prop}[ex]{Proposition}
\newtheorem{coro}[ex]{Corollary}
\newtheorem{lemma}[ex]{Lemma}
\newtheorem{defi}[ex]{Definition}

\newcommand{\F}{\mathbb{F}}

\newcommand{\N}{\mathbf{N}}
\newcommand{\M}{\mathbf{M}}
\newcommand{\E}{\mathbf{E}}

\newcommand{\CC}{\mathbb{C}}

\newcommand{\ie}{\emph{i.e.}\;}

\newcommand{\FF}{\mathfrak{F}}
\newcommand{\OO}{\mathcal{O}}
\newcommand{\Ig}{\mathfrak{I}}

\newcommand{\GG}{\mathcal{G}}

\newcommand{\HH}{\mathcal{H}}
\newcommand{\EE}{\mathcal{E}}

\newcommand{\XX}{\mathfrak{X}}

\newcommand{\pp}{\mathfrak{p}}
\newcommand{\qq}{\mathfrak{q}}

\newcommand{\aaa}{\mathfrak{a}}
\newcommand{\bbb}{\mathfrak{b}}
\newcommand{\ccc}{\mathfrak{c}}
\newcommand{\ddd}{\mathfrak{d}}

\newcommand{\fff}{\mathfrak{f}}

\newcommand{\NN}{\mathbb{N}}
\newcommand{\ZZ}{\mathbb{Z}}
\newcommand{\QQ}{\mathbb{Q}}
\newcommand{\UU}{\mathbb{U}}

\newcommand{\alg}{\mathrm{alg}}
\newcommand{\ab}{\mathrm{ab}}
\newcommand{\abi}{{\mathrm{ab},\infty}}

\newcommand{\Ci}{\mathbf{C}_\infty}

\newcommand{\tor}{\mathrm{tor}}
\newcommand{\adm}{\mathrm{adm}}

\newcommand{\RR}{\mathbb{R}}

\newcommand{\abs}[1]{\left\vert#1\right\vert}

\newcommand{\norm}[1]{\left\Vert#1\right\Vert}

\newcommand{\cen}[1]{\begin{center} #1 \end{center}}

\newcommand{\Hay}{H(\sgn)}

\newcommand{\I}{\mathrm{I}}
\newcommand{\II}{\mathrm{II}}
\newcommand{\III}{\mathrm{III}}

\newcommand{\n}{\newline}

\newcounter{c_liste}
\newenvironment{mylist}{\begin{list}{{\textup{(\arabic{c_liste})}}}{\usecounter{c_liste}\leftmargin 5mm \itemsep 0mm}}{\end{list}}

\DeclareMathOperator{\re}{Re}
\DeclareMathOperator{\im}{Im}

\DeclareMathOperator{\Ker}{Ker}

\DeclareMathOperator{\id}{id}

\DeclareMathOperator{\Tr}{Tr}
\DeclareMathOperator{\tr}{tr}
\DeclareMathOperator{\Gal}{Gal}
\DeclareMathOperator{\Card}{Card}

\DeclareMathOperator{\Aut}{Aut}
\DeclareMathOperator{\End}{End}

\DeclareMathOperator{\sgn}{sgn}
\DeclareMathOperator{\ann}{ann}

\renewcommand{\indexname}{Index of notations}

\makeindex

\begin{document}

\title{Bost-Connes type systems for function fields}
\author{Benoît Jacob}
\date{}
\maketitle

\vskip 1 cm

\begin{abstract}
We describe a construction which associates to any function field $k$ and any place $\infty$ of $k$ a $C^*$-dynamical system $(C_{k,\infty},\sigma_t)$ that is analogous to the Bost-Connes system associated to $\QQ$ and its archimedian place. Our construction relies on Hayes' explicit class field theory in terms of sign-normalized rank one Drinfel'd modules. We show that $C_{k,\infty}$ has a faithful continuous action of $\Gal(K/k)$, where $K$ is a certain field constructed by Hayes, such that $k^\abi\subset K\subset k^\ab$, where $k^\abi$ is the maximal abelian extension of $k$ that is totally split at $\infty$. We classify the extremal KMS$_\beta$ states of $(C_{k,\infty},\sigma_t)$ at any temperature $0<1/\beta<\infty$ and show that a phase transition with spontaneous symmetry breaking occurs at temperature $1/\beta=1$. At high temperature $1/\beta\geqslant 1$, there is a unique KMS$_\beta$ state. At low temperature $1/\beta<1$, the space of extremal KMS$_\beta$ states is principal homogeneous under $\Gal(K/k)$. Each such state is of type $\I_\infty$ and the partition function is the Dedekind zeta function $\zeta_{k,\infty}$. Moreover, we construct a ``rational'' $*$-subalgebra $\HH$, we give a presentation of $\HH$ and of $C_{k,\infty}$, and we show that the values of the low-temperature extremal KMS$_\beta$ states at certain elements of $\HH$ are related to special values of partial zeta functions.\n

{\bf Erratum:} This article wrongly claims that at high temperature $1/\beta\geqslant 1$, the unique KMS$_\beta$ state is of type $\III_{q^{-\beta}}$, where $q$ is the cardinal of the constant subfield of $k$. It has been shown by Neshveyev and Rustad \cite{NesRus12} that the correct type is $\III_{q^{-\beta d_\infty}}$ where $d_\infty$ is the degree of the place $\infty$. The original statements have been kept for reference, but errata have been inserted next to them.
\end{abstract}

\pagebreak

\tableofcontents{}

\pagebreak

\addcontentsline{toc}{section}{\indexname}

\printindex

\addcontentsline{toc}{section}{Introduction}
\section*{Introduction}

\paragraph*{Statement of the main results.} Let \index{k@$k$} $k$ be any global function field. Let \index{0@$\infty$} $\infty$ be any place of $k$. In this paper, we will 
associate to the pair $(k,\infty)$ a $C^*$-dynamical system $(C_{k,\infty}, (\sigma_t))$.\n

Our system aims to be an analog of the Bost-Connes (BC for short) system associated to $\QQ$, cf. Bost and Connes \cite{BosCon95}. The partition function of the BC system is the Riemann zeta function without the $\Gamma$-factor at infinity. Similarly, we will check (Lemma \ref{xr36}) that the partition function of our system is the zeta function of the field $k$ without the factor corresponding to the place $\infty$ of $k$.\\

The BC system admits $\Gal(\QQ^\ab/\QQ)$ as symmetry group. Similarly, we will check (Proposition \ref{galsym}) that our system has $\Gal(K/k)$ as symmetry group (meaning that $\Gal(K/k)$ acts continuously and faithfully on $C_{k,\infty}$, commuting with the flow $\sigma_t$), where $K$ is a field having the following property:
$$k^\abi\subset K \subset k^\ab,$$
where \index{k_ab_i@$k^\abi$} $k^\abi$ is the maximal abelian extension of $k$ that is totally split at $\infty$. The field $K$ is generated over $k$ by coefficents and torsion points of certain rank one Drinfel'd modules; this is part of David R. Hayes' explicit class field theory for function fields, cf. Hayes \cite{Hay74}, \cite{Hay79} and \cite{Hay92}, which we will quickly review. If $\infty'$ is any place of $k$ other than $\infty$, we have (cf. \cite{Hay74}, Theorem 7.2)
$$k^{\ab,\infty'}\cdot k^\abi=k^\ab.$$

We will construct our $C^*$-algebra $C_{k,\infty}$ as the maximal $C^*$-algebra of a certain groupoid $\GG$. We will also (Proposition \ref{presentC}) give a presentation of $C_{k,\infty}$ as a $C^*$-algebra.\\

For any temperature $1/\beta\in\RR^*_+$, let \index{K_beta@$K_\beta$ ($\beta\in\RR^*_+$)} $K_\beta$ be space of KMS$_\beta$ states of $\left(C_{k,\infty},(\sigma_t)\right)$, endowed with the weak$^*$ topology. By Bratteli and Robinson \cite{BraRob}, II, Theorem 5.3.30, the space $K_\beta$ is a compact simplex (in particular, it is convex). Let \index{E_K_beta@$\EE(K_\beta)$ ($\beta\in\RR^*_+$)}$\EE(K_\beta)$ denote the subspace of extreme points of $K_\beta$. The elements of $\EE(K_\beta)$ are called the \emph{extremal} KMS$_\beta$ states. By loc. cit., a KMS$_\beta$ state is extremal if, and only if it is a factor state. Thus, $\EE(K_\beta)$ is equal to the space of KMS$_\beta$ factor states.\\

We will classify the KMS$_\beta$ states of our system for any temperature $1/\beta\in\RR^*_+$:\\

At low temperature $1/\beta<1$, we will prove (Theorem \ref{propkms}) that $\EE(K_\beta)$ is principal homogeneous\footnote{Let $G$ be topological group acting on a topological space $X$. One says that $X$ is \emph{principal homogeneous} under $G$ if, for any $x\in X$, the map $g\mapsto gx$ is a homeomorphism $G\rightarrow X$.} under $\Gal(K/k)$. The states in $\EE(K_\beta)$ are of type $\I_\infty$ (Proposition \ref{factorIi}).\\

At high temperature $1/\beta\geqslant 1$, we will prove (Theorem \ref{xr61}) that there exists a unique KMS$_\beta$ state.\\

{\bf Erratum:} This article wrongly claims that at high temperature $1/\beta\geqslant 1$, the unique KMS$_\beta$ state is of type $\III_{q^{-\beta}}$, where $q$ is the cardinal of the constant subfield of $k$. It has been shown by Neshveyev and Rustad \cite{NesRus12} that the correct type is $\III_{q^{-\beta d_\infty}}$ where $d_\infty$ is the degree of the place $\infty$. The original statements have been kept for reference, but errata have been inserted next to them.\\

We will construct a dense $*$-subalgebra $\HH$ which gives an arithmetic structure to our dynamical system, as in \cite{BosCon95}. For example, we will show (Theorem \ref{xr40}) that evaluating low-temperature extremal KMS$_\beta$ states on certain elements of the subalgebra $\HH$ gives rise to formulas involving special values of partial zeta functions.\\

Many of our proofs are adapted from \cite{BosCon95}, and we have also borrowed several ideas from Harari and Leichtnam \cite{HarLei97}.

\paragraph*{Outline.} This paper is divided into four sections. In Section 1, we first review definitions and results in the arithmetic of function fields and in the analytic theory of Drinfel'd modules. We review Hayes' explicit class field theory for function fields, in terms of sign-normalized rank one Drinfel'd modules. We choose once and for all a sign-function $\sgn$, and Hayes' theory provides us with a finite set $H(\sgn)$ of Drinfel'd modules with special arithmetic properties. In particular, their coefficients and torsion points generate the extension $K/k$ which we mentionned above. In the rest of this paper, the only Drinfel'd modules which we consider are the elements of $H(\sgn)$.\\

In Section 2, we do the actual construction of the $C^*$-dynamical system $(C_{k,\infty}, (\sigma_t))$. From the finite set $H(\sgn)$ provided by Hayes' theory, we construct a compact topological space $X$ in the following way: for any $\phi\in H(\sgn)$, let $X_\phi$ denote the dual group of the discrete group of torsion points of the Drinfel'd module $\phi$. Let $X$ be the disjoint union of the $X_\phi$, where $\phi$ runs over $H(\sgn)$. The compact space $X$ is endowed with a natural action of the semigroup $\Ig_\OO$ of ideals. This gives rise to a groupoid $\GG$, and the $C^*$-algebra $C_{k,\infty}$ is obtained as the maximal groupoid $C^*$-algebra of $\GG$. The flow $(\sigma_t)$ is then easy to define.\\

In Section 3, we prove a number of results about the algebraic structure of $(C_{k,\infty}, (\sigma_t))$. We introduce a $*$-subalgebra $\HH$ which plays the rôle of the rational Hecke algebra in the paper \cite{BosCon95}. We prove that $\HH$ is dense in $C_{k,\infty}$, and we give a presentation of $\HH$ as a $*$-algebra and of $C_{k,\infty}$ as a $C^*$-algebra. We then study an action of $\Gal(K/k)$ on $C_{k,\infty}$, and compute the fixed-point subalgebra $C_1$. The rest of this section is devoted to miscellaneous arithmetical results which we use in the last section.\\

In Section 4, for any temperature $1/\beta\in\RR^*_+$, we describe the space $\EE(K_\beta)$ of extremal KMS$_\beta$ states (endowed with the weak$^*$ topology), and we compute the type of all such states. We first construct a KMS$_\beta$ state $\varphi_\beta$, and show that it is the unique $\Gal(K/k)$-invariant KMS$_\beta$ state. We then show that the action of $\Gal(K/k)$ on $\EE(K_\beta)$ is transitive and continuous. Thus, in order to describe $\EE(K_\beta)$, it is enough to find an element of $\EE(K_\beta)$ and to describe its orbit under $\Gal(K/k)$. At low temperature $1/\beta<1$, we associate to any admissible character $\chi$ a Gibbs state $\varphi_{\beta,\chi}$ in the regular representation at $\chi$. We prove that the map $\chi\mapsto\varphi_{\beta,\chi}$ is a homeomorphism from the space $X^\adm$ of admissible characters to $\EE(K_\beta)$. We also prove that both spaces are principal homogeneous under $\Gal(K/k)$. We check that the states in $\EE(K_\beta)$ are of type $\I_\infty$, that the 
partition function is the Dedekind zeta function $\zeta_{k,\infty}$, and we compute the values of the $\varphi_{\beta,\chi}$ at some points of $\HH$, in terms of special values at $\beta$ of partial zeta functions of $k$. At high temperature $1/\beta\geqslant 1$, we prove that $\EE(K_\beta)=\{\varphi_\beta\}$.

\paragraph*{Literature on Bost-Connes type constructions.} The 1995 paper \cite{BosCon95} has inspired many mathematicians. Unfortunately, it would be impossible to mention all of them here; we refer to Section 1.4 of Connes and Marcolli \cite{ConMar04} for a more complete summary. M. Laca, N. Larsen, I. Raeburn and others have investigated in a number of papers (see for instance \cite{ArlLacRae97}, \cite{Lac98}, \cite{LacRae00}, \cite{LacFra04}, \cite{LarRae02}) the semigroup crossed product and Hecke algebra aspects of the BC construction and generalizations of it. In 1997, D. Harari and E. Leichtnam have obtained in \cite{HarLei97} a system with spontaneous symmetry breaking for any global field. In 1999, P. Cohen has obtained in \cite{Coh99} a system for number fields whose partition function is the Dedekind zeta function. In 2002, S. Neshveyev has given in \cite{Nes02} a new proof of the unicity of the KMS$_\beta$ state at high temperature. In 2004, A. Connes and M. Marcolli have introduced in \cite{
ConMar04} the noncommutative space of $\QQ$-lattices up to scaling and commensurability, allowing for a comprehensive reformulation of the BC construction, and have studied the case of rank $2$. In 2005, A. Connes, M. Marcolli and N. Ramachandran  have obtained in \cite{ConMarRam05}, \cite{ConMarRam06} the ``good'' system for quadratic imaginary number fields, and have studied its relation to complex multiplication of elliptic curves. The same year, E. Ha and F. Paugam have extended in \cite{HaPau05} the Connes-Marcolli setting to arbitrary Shimura varieties. Lastly, in the paper \cite{ConConMar05}, A. Connes, C. Consani and M. Marcolli have introduced the notion of an \emph{endomotive}, putting the BC construction into a much wider perspective which also includes A. Connes' spectral realization \cite{Con99} of the zeroes of the Riemann zeta function.

\paragraph*{Acknowledgements. }Sergey Neshveyev and Simen Rustad \cite{NesRus12} found and corrected the important mistake I made in Theorem \ref{xr80} and in the lemmas used to prove it. I have incorporated errata in the relevant places of this article while keeping the original wrong statements for reference.\\

I thank \'Eric Leichtnam for giving me this research subject, for many helpful comments on early versions of this paper, and for his explanations on operator algebras. I thank Matilde Marcolli for encouraging me to give a talk on this material at MPI, Bonn, in October 2005. During this research, I enjoyed the excellent working environment of the ``Projet Algèbres d'Opérateurs'' at Jussieu, and I thank \'Etienne Blanchard for letting me talk in its seminar. I am very grateful to Alain Connes, David R. Hayes, Georges Skandalis and Stefaan Vaes who kindly answered mathematical questions. I also thank Cécile Armana, Pierre Fima, Eugene Ha, Cyril Houdayer and Frédéric Paugam for helpful discussions.

\paragraph*{Remarks. } Here are two interesting remarks that people made at the end of my MPI talk.\\

\noindent {\bf 1.} As Alain Connes pointed out, our system lacks one feature of the BC system: fabulous states. The reason for that is obvious: values of states are elements of $\CC$, so the symmetry group $\Gal(K/k)$ doesn't act naturally on them. Obtaining fabulous states would require to have a theory of dynamical systems of positive characteristic, where states would take values in some field of positive characteristic. Note that even though the low temperature extremal KMS$_\beta$ states of our system do not have the fabulous property, they have interesting special values (Theorem \ref{xr40}).\\

\noindent {\bf 2.} Arkady Kholodenko mentionned that it might be possible to adapt his work on $2+1$ gravity \cite{Kho01} in order to obtain zeta functions of function fields as partition functions, and that Drinfel'd modules should play a rôle.

\paragraph*{Notations. } In this paper, \index{N@$\NN$ (nonnegative integers)} $\NN$ denotes the set of nonnegative integers, \index{N_star@$\NN^*$ (positive integers)} $\NN^*$ denotes the set of positive integers, and \index{Rplusstar@$\RR^*_+$ (positive real numbers)} $\RR^*_+$ denotes the set of positive real numbers. Thus, $0\in\NN$, $0\not\in\NN^*$, and $0\not\in\RR^*_+$. For any Hilbert space $H$, we let \index{B_H@$B(H)$ ($H$ a Hilbert space)} $B(H)$ denote the algebra of all bounded linear operators on $H$. For any set $X$, we write \index{B_l_2@$B\ell^2$} $B\ell^2(X)$ for $B(\ell^2(X))$. For any $x\in\RR$, we set
\index{1@$\lfloor x\rfloor$ ($x\in\RR$)} $$\lfloor x\rfloor\:=\:\max\,\{n\in\ZZ\:;\;\;n\leqslant x\}.$$
For any predicate $P$, we define \index{1P@$1_P$ ($P$ a predicate)} $1_P$ to be equal to $1$ if $P$ is true, and $0$ if $P$ is false. Thus, we have, for any two 
predicates $P$ and $Q$:$$1_{P\:\text{and}\:Q}=1_P 1_Q.$$

\section{Function fields, Drinfel'd modules, and Hayes' explicit class field theory}

\subsection{Function fields}

\noindent Here are three equivalent definitions of a \emph{function field}:
\begin{itemize}
\item A field which is a finite extension of $\F_p(T)$, for some prime number $p$.
\item A global field of positive characteristic.
\item The field $K(C)$ of rational functions on a projective curve $C$ over a finite field. The curve $C$ can always be chosen to be smooth.
\end{itemize}

Thus, global fields fall into two categories: those of characteristic $0$ are the number fields, and those of positive characteristic are the function fields.\\

\noindent Recall that at the beginning of this paper, we chose a function field $k$ and a place $\infty$ of $k$.\\

\noindent Function fields have many similarities with number fields. An important part of algebraic number theory works in the same way for all global fields.\\

The analog of the Dedekind ring of integers is defined as follows. According to the third definition of a function field, view $k$ as the field $K(C)$ of rational functions on a smooth projective curve $C$ over a finite field. View $\infty$ as a closed point of $C$. Let $\index{O@$\OO$}\OO$ be the subring of $k$ of all functions having no pole away from $\infty$. In other words, $\OO$ is the ring of regular functions on the affine curve $C-\{\infty\}$. Note that $k=K(C)$ is the field of fractions of $\OO$.

\paragraph*{Example.} $k=\F_p(T)$ and $\infty$ is the place corresponding to an absolute value $\abs{\cdot}$ such that $\abs{T}>1$. The subring $\OO$ is then the polynomial ring $\F_p[T]$.\\

Call \emph{finite} the places of $k$ other than $\infty$. We have a natural bijection
\cen{finite places of $k$ $\;\;\longleftrightarrow\;\;$ maximal ideals of $\OO$.}

Let \index{p@$p$} $p$ denote the characteristic of $k$. The range of the unique unital ring morphism $\ZZ\rightarrow k$ is a finite field with $p$ elements; we note it \index{F_p@$\F_p$}$\F_p$. The algebraic closure of $\F_p$ in $k$ is called the \emph{constant subfield} of $k$. Let \index{q@$q$} $q$ denote its cardinal. Of course, $q$ is a power of $p$. We let \index{F_q@$\F_q$} $\F_q$ denote the constant subfield of $k$. An element of $k$ is said to be \emph{constant} if it belongs to $\F_q$.\\

For any place $\pp$ of $k$, we let \index{Np@$\N\pp$ ($\pp$ a place of $k$)} $\N\pp$ denote the cardinal of the residue field of $\pp$. Thus, $\N\pp=q^{n_\pp}$ for some positive integer $n_\pp$ called the \emph{degree} of $\pp$. Note that, if $\pp$ is finite, then the residue field is the quotient $\OO/\pp$. \\

The rest of this subsection is a review of a few well-known theorems about function fields, which will be used in the proofs of our classification of KMS$_\beta$ states. These theorems are: the strong approximation theorem, Weil's ``Riemann Hypothesis for curves'', and the abelian case of the \v{C}ebotarev density theorem for the natural density. The first one will be used in subsection \ref{admsubsec}, which in turn will be used in the classification of KMS$_\beta$ states at low temperature. The two other ones will be used  in the classification of KMS$_\beta$ states at high temperature.\\

Let \index{A_f@$A_f$} $A_f$ denote the ring of finite adèles of $k$. This is the restricted product of the $k_\pp$ 
with respect to the $\OO_\pp$, where $\pp$ runs over all finite places of $k$. Let $\iota_f\::\:k\hookrightarrow A_f$ be the diagonal embedding.

\begin{theo} {\bf Strong approximation theorem.}
\label{strongapprox}
The field $\iota_f(k)$ is dense in $A_f$.
\end{theo}
\begin{proof}
See Cassels and Fröhlich \cite{CasFro67}, Chapter II, \S 15, page 67.
\end{proof}
This is contrasted with the fact that if $\iota\::\:k\hookrightarrow A$ is the diagonal embedding into the full ring of adèles, then $\iota(k)$ is discrete in $A$. Note that
$$A = A_f \times k_\infty,$$
where \index{k_i@$k_\infty$} $k_\infty$ is the completion of $k$ at $\infty$.\\

Let us now recall Weil's ``Riemann hypothesis for curves'' theorem. The \emph{genus} of a function field is the genus of any projective smooth curve of which it is the function field. For the statement of the following theorem, we temporarily forget that we already chose a function field $k$ and defined $q$ as the cardinal of its constant subfield.

\begin{theo} {\bf (A. Weil, the Riemann Hypothesis for curves).}
\label{riemann}
Let $k$ be a function field of genus $g$. Let $q$ be the cardinal of its constant subfield. Let $N$ be the number of places of $k$ with norm $q$ (\ie with degree $1$). Then
$$ q - 2g\sqrt{q} + 1\leqslant N\leqslant q + 2g\sqrt{q} + 1.$$
\end{theo}
\begin{proof}
Weil's original proof is published in \cite{Wei48}.
\end{proof}

Let us now come back to the function field $k$ that we fixed at the beginning of this paper. Let \index{g@$g$ (genus of $k$)} $g$ denote the genus of $k$.\\

Given an integer $n\geqslant 1$, one may ask how to obtain a result similar to Theorem \ref{riemann} for places of $k$ with norm $q^n$ (\ie with degree $n$). Note that one cannot replace $q$ by $q^n$ in Theorem \ref{riemann}. Here one has to be wary of the distinction between closed points, which correspond to places of $k$, and geometric points, which correspond to places of suitable extensions of $k$. The following corollary will be used in subsections \ref{uniqkms} and \ref{typeiii}.

\begin{coro}
\label{riemann2}
For any $n\geqslant 1$, let \index{Q_k_q_n@$Q(k,q^n)$ ($n\in\NN^*$)} $Q(k,q^n)$ denote the number of places of $k$ with norm $q^n$, and let \index{P_k_q_n@$P(k,q^n)$ ($n\in\NN^*$)} $P(k,q^n)$ denote the number of places $k$ with norm $\leqslant q^n$. The following estimates hold when $n\rightarrow\infty$:
\begin{eqnarray}
Q(k,q^n) & = & \frac{q^n}{n}\,+\,O\left(q^{n/2}\right)\label{riemann3},\\
P(k,q^n) & \sim & \frac{q}{q-1}\,\cdot\,\frac{q^n}{n}\cdot\label{riemann4}
\end{eqnarray}
\end{coro}
\begin{proof}
For any $n\geqslant 1$, let $k_n=k\otimes_{\F_q} \F_{q^n}$. Note that the constant subfield of $k_n$ is $\F_{q^n}$. Let $N_n$ denote the number of places of $k_n$ with norm $q^n$. By Theorem \ref{riemann} applied to the function field $k_n$, we have
\begin{equation}
\label{riemann5}
q^n - 2gq^{n/2} + 1\leqslant N_n\leqslant q^n + 2gq^{n/2} + 1.
\end{equation}
Let $n\geqslant 1$. One easily checks that for any $m\mid n$ there is a bijection
\cen{places of $k$ with norm $q^m$ $\;\;\longleftrightarrow\;\;$ $\Gal(k_n/k)$-orbits with cardinal $m$ of places of $k_n$ with norm $q^n$.}
Thus, we have
\begin{equation}
\label{orbits}
N_n=\sum_{m\mid n} m\,Q(k,q^m).
\end{equation}
This gives
$$nQ(k,q^n)=N_n-\sum_{m\mid n,\:m\leqslant n/2} m\,Q(k,q^m).$$
By equation (\ref{orbits}), we have $mQ(k,q^m)\leqslant N_m$, so we find
\begin{eqnarray*}
N_n\geqslant n\,Q(k,q^n)& \geqslant & N_n-\sum_{m\mid n,\:m\leqslant n/2} N_m.\\
& \geqslant & N_n-(n/2)N_{\lfloor n/2 \rfloor}.
\end{eqnarray*}
Applying the inequality (\ref{riemann5}), we get
$$q^n + 2gq^{n/2} + 1 \geqslant n\,Q(k,q^n) \geqslant q^n - 2gq^{n/2} + 1 -(n/2)(q^{n/2} + 2gq^{n/4} + 1),$$
and the estimate (\ref{riemann3}) follows. From the estimate (\ref{riemann3}), using the equality
$$P(k,q^n)=\sum_{m=1}^n Q(k,q^m),$$
one can obtain the estimate (\ref{riemann4}) by an elementary computation.
\end{proof}

For any $s\in\CC$ with $\re s>1$, put
$$\zeta_k(s)\:=\:\prod_\pp \frac{1}{1-\N\pp^{-s}}$$
where the product is taken over all places of $k$. One shows that $\zeta_k$ can be continued to a meromorphic function on $\CC$. Note that $\zeta_k$ is periodic, with period $2\pi i/\log q$. The inequation (\ref{riemann5}) for all $n\geqslant 1$ is then equivalent to the statement that all zeroes of $\zeta_k$ have real part $1/2$. One defines the \emph{zeta function without the factor at} $\infty$, noted \index{zeta_k_inf@$\zeta_{k,\infty}$} $\zeta_{k,\infty}$, to be the meromorphic continuation of the function defined when $\re s>1$ by
$$\zeta_{k,\infty}(s)\:=\:\prod_{\pp\neq\infty} \frac{1}{1-\N\pp^{-s}}\:=\:(1-\N\infty^{-s})\zeta_k(s)\cdot$$
Note that when $\re s > 1$, we have
$$\zeta_{k,\infty}(s)\:=\:\sum_{\aaa\in\Ig_\OO} \frac{1}{\N\aaa^s}\cdot$$

Let us now recall a version of the \v{C}ebotarev density theorem.\\

Let $S$ denote the set of all places of $k$. A set $P$ of places of $k$ is said to \emph{have a Dirichlet density} if the following limit exists in $\RR$:
\index{d_p@$d(P)$ ($P$ a set of places)} $$d(P)=\lim_{s\rightarrow 1_+} \frac{\sum_{\pp\in P} \N\pp^{-s}}{\sum_{\pp\in S} \N\pp^{-s}}\:.$$
Moreover, $P$ is said to \emph{have a natural density} if the following limit exists in $\RR$:
\index{delta_p@$\delta(P)$ ($P$ a set of places)} $$\delta(P)=\lim_{N\rightarrow +\infty} \frac{\Card \{\pp\in P,\;\;\N\pp \leqslant N\}}{\Card \{\pp\in S,\;\;\N\pp \leqslant N\}}\cdot$$
If a set $P$ has a natural density, then it also has a Dirichlet density, and $d(P)=\delta(P)$.

\begin{theo} {\bf \v{C}ebotarev density theorem, abelian case, for the natural density.}
\label{cebotarev}
Let $L$ be a finite abelian extension of $k$. Let $\sigma\in\Gal(L/k)$. Let $P$ denote the set of all places $\pp$ of $k$ unramified in $L$ and such that $\sigma_\pp=\sigma$, where $\sigma_\pp=(\pp,L/k)\in\Gal(L/k)$ is the Artin automorphism of $L$ associated to $\pp$. Then $P$ has natural density $\delta(P)=1/[L:k]$. Therefore, it also has Dirichlet density $d(P)=1/[L:k]$.
\end{theo}

{\bf Erratum:} Sergey Neshveyev kindly points out that the above formulation of the Cebotarev density theorem is wrong, as for function fields it only holds for the Dirichlet density, not for the natural density. This only affects the type III computation mentioned in the Erratum in the Absract. Consequently, the corollary below is wrong, but that also doesn't affect anything else than the type III computation already mentioned in that Erratum.

\begin{coro}
\label{xr2}
Let $L$ be a finite abelian extension of $k$. Let $\sigma\in\Gal(L/k)$. For any $n\geqslant 1$, let $P(L/k,q^n,\sigma)$ denote the number of places of $\pp$ of $k$ unramified in $L$ such that $\N\pp\leqslant q^n$ and $\sigma_\pp=\sigma$, where $\sigma_\pp=(\pp,L/k)\in\Gal(L/k)$ is the Artin automorphism of $L$ associated to $\pp$. Let $Q(L/k,q^n,\sigma)$ denote the number of places of $\pp$ of $k$ unramified in $L$ such that $\N\pp=q^n$ and $\sigma_\pp=\sigma$. The following estimates hold when $n\rightarrow\infty$:
\begin{eqnarray}
P(L/k,q^n,\sigma) & \sim & \frac{q}{(q-1)\,[L:k]}\,\cdot\,\frac{q^n}{n}\label{riemann6}\;,\\
Q(L/k,q^n,\sigma) & \sim & \frac{1}{[L:k]}\,\cdot\,\frac{q^n}{n}\:\cdot\label{riemann7}
\end{eqnarray}
\end{coro}
\begin{proof}
The estimate (\ref{riemann6}) follows from Theorem \ref{cebotarev} and the estimate (\ref{riemann4}). We have
$$Q(L/k,q^n,\sigma)=P(L/k,q^n,\sigma)-P(L/k,q^{n-1},\sigma),$$
so
$$\frac{Q(L/k,q^n,\sigma)}{P(k,q^n)}\:=\:\frac{P(L/k,q^n,\sigma)}{P(k,q^n)}\,-\,\frac{P(L/k,q^{n-1},\sigma)}{P(k,q^{n-1})}\cdot\frac{P(k,q^{n-1})}{P(k,q^n)}\cdot$$
Hence
$$\frac{Q(L/k,q^n,\sigma)}{P(k,q^n)}\:\xrightarrow{n\rightarrow\infty}\: \frac{1}{[L:k]} \,-\, \frac{1}{[L:k]} \cdot\frac{1}{q}\:=\:\frac{q-1}{q[L:k]}\cdot$$
Applying the estimate (\ref{riemann4}) to that, we get the estimate (\ref{riemann7}).
\end{proof}

\subsection{Drinfel'd modules over $\Ci$}

Our references in this subsection are \cite{Hay92} and Chapter IV of Goss \cite{Gos91}.\\

Recall that the maximal abelian extension of a quadratic imaginary number field is generated by the $j$-invariant and the torsion points of a suitable elliptic curve over $\CC$. One wishes to develop a similar theory for function fields. Thus, one looks for good analogs of $\CC$ and of the notion of an elliptic curve over $\CC$. The analog of the field $\CC$ has been well known for a long time, and is what we will note $\Ci$. The analog of the notion of an elliptic curve over $\CC$ is going to be the notion of a Drinfel'd module over $\Ci$.\\

We begin with describing the analog of $\CC$. Let \index{k_i@$k_\infty$} $k_\infty$ be the completion of $k$ at $\infty$. The problem is that $k_\infty$ is not algebraically closed. Take an algebraic closure $k_\infty^\alg/k_\infty$. One shows that $\infty$ extends uniquely to a place of $k_\infty^\alg$. Then, the problem is that $k_\infty^\alg$ is not complete. So let \index{C_i@$\Ci$} $\Ci$ \label{index_Ci} denote the completion of $k_\infty^\alg$ at $\infty$. The field $\Ci$ is both complete and algebraically closed.\\

Let us choose once and for all an imbedding $\iota\::\:k\hookrightarrow \Ci$, and use it to view $k$ as a subfield of $\Ci$.

\paragraph*{Lattices.} We are now ready to introduce Drinfel'd modules. The most concrete way to introduce elliptic curves over $\CC$ is to first define \emph{lattices} in $\CC$. Similarly, we are going to first define lattices in $\Ci$.\\

Recall that $\OO$ is the subring of integers of $k$, defined in the previous subsection. A subgroup $L\subset\Ci$ is said to be \emph{discrete} if there exists a neighborhood $U$ of $0$ in $\Ci$ such that $U\cap L=\{0\}$.

\begin{defi}
An $\OO$-\emph{lattice} in $\Ci$ is a discrete, finitely generated $\OO$-submodule of $\Ci$.
\end{defi}
\noindent We will say ``lattice'' instead of ``$\OO$-lattice in $\Ci$''.\\

That is an abstract definition, but in this paper, we will only have to deal with a special case of lattices, rank one lattices, for which there is a very concrete definition. Let us first define the \emph{rank} of a lattice.\\

Let $L$ be a lattice. As $\Ci$ is a field containing $\OO$, it is obviously a torsion-free $\OO$-module. Hence $L$ is also torsion-free. As $\OO$ is a Dedekind ring, the $\OO$-module $L$, being finitely generated and torsion-free, is automatically projective, so there exist an integer $r\geqslant 1$ and ideals $\aaa_1,\ldots,\aaa_r\in\Ig_\OO$ such that $L$ is isomorphic as an $\OO$-module to $\aaa_1\oplus\cdots\oplus\aaa_r$. 

\begin{defi}
The integer $r$ above is called the \emph{rank} of $L$.
\end{defi}
\noindent Let \index{IO@$\Ig_\OO$} $\Ig_\OO$ be the semigroup of all nonzero ideals of $\OO$, under the usual multiplication law of ideals.\\

\noindent For rank one lattices, we have the following result:

\cen{\emph{A subset of $\Ci$ is a rank one lattice if, and only if it is of the form $\xi\aaa$ with $\xi\in\Ci^*$ and $\aaa\in \Ig_\OO$.}}

\paragraph*{The Drinfel'd module associated to a lattice.} Let $L$ be a lattice (of any rank). Remember the following product formula:
$$\forall z\in \CC,\;\;\;\sin z=z\prod_{t\in \pi\ZZ-\{0\}}(1-z/t).$$
Similarly, let us define a function $e_L\::\:\Ci\rightarrow\Ci$ by the following formula:
\index{eL@$e_L$ ($L$ a lattice)} $$\forall x\in \Ci,\;\;\;e_L(x)=x\prod_{\ell\in L-\{0\}}(1-x/\ell).$$
One shows that this product converges for all $x$. The function $e_L$ should be called the ``sinus function associated to $L$'', but authors have decided to call it the ``exponential function associated to $L$''. We have
\begin{equation}
\label{e_L_add}
\forall x,y\in\Ci,\;\;\;e_L(x+y)=e_L(x) + e_L(y),
\end{equation}
and
\begin{equation}
\label{e_L_mod}
\forall a\in\OO,\;\forall x\in\Ci,\;\;\;e_L(ax) = \phi^L_a (e_L(x)),
\end{equation}
where $\phi^L_a\in\Ci[X]$ is the polynomial given by the following formula if $a\neq 0$:
$$\phi^L_a=aX\prod_{0\neq \ell\in a^{-1}L/L} (1-X/e_L(\ell)),$$
and $\phi^L_0=0$. Note that if $a$ is a nonzero constant (that is, $a\in\F_q^*$), then it is invertible in $\OO$, and hence $a^{-1}L=L$. Thus, one has
\begin{equation}
\label{phifq}
\forall a\in\F_q,\;\;\;\phi_a = aX.
\end{equation}
As we will shortly see, this allows to check that for any $a\in\OO$, the polynomial $\phi^L_a$ is $\F_q$-linear, which means that in can have nonzero coefficients only in degrees that are powers of $q$.\\

Equation (\ref{e_L_add}) is an analog of the classical formula for $\sin(x+y)$, not of the formula for $\exp(x+y)$. The fact that $e_L$ is additive, while $\sin$ isn't, is a phenomenon typical of characteristic $p$ algebra, just like the additivity of the Frobenius map $x\mapsto x^p$. The polynomials $\phi^L_a$ can be viewed as analogs of the classical Chebycheff polynomials of trigonometry.\\

One shows, by analytic means, that $e_L$ induces a bijection
$$e_L\;:\;\Ci/L\rightarrow \Ci.$$
So this is a group isomorphism. Use it to transport the $\OO$-module structure
of $\Ci/L$ to a new $\OO$-module structure on $\Ci$, which we note
\index{phi_Ci@$\phi(\Ci)$ ($\phi$ a Drinfel'd module)} $\phi^L(\Ci)$. Thus, $\phi^L(\Ci)$ is the $\OO$-module that is equal to $\Ci$ as an additive group, and whose $\OO$-module structure is given by
$$(a,x)\mapsto \phi^L_a(x).$$
Thus, by definition, the map $e_L$ is an isomorphism of $\OO$-modules
\begin{equation}
\label{e_L_isom}
e_L\;:\;\Ci/L\;\rightarrow\;\phi^L(\Ci).
\end{equation}
The \emph{Drinfel'd module associated to $L$} is the map
\begin{eqnarray}
\phi^L\;:\;\OO & \rightarrow & \Ci[X],\nonumber\\
a & \mapsto & \phi^L_a.\nonumber
\end{eqnarray}

\paragraph*{Definition of a Drinfel'd module over $\Ci$.} The map $\phi^L$ that we have just defined satisfies
\begin{eqnarray}
\forall a,b\in\OO,\;\;\;\phi^L_{a+b}& =& \phi^L_a + \phi^L_b. \label{phiadd}\\
\forall a,b\in\OO,\;\;\;\phi^L_{ab}& =& \phi^L_a \circ \phi^L_b\;=\;\phi^L_b\circ\phi^L_a, \label{phimul}
\end{eqnarray}
Let $\tau=X^q$ and, for $n\geqslant 0$, $\tau^n=X^{q^n}$. In particular, $\tau^0=X$. Let $\Ci\{\tau\}$ denote the (noncommutative) $ \Ci$-algebra whose underlying vector space is the $\Ci$-linear span of the $\tau^n$, for $n\geqslant 0$, and where the ``multiplication'' law is the composition law $\circ$. Note that $\Ci\{\tau\}$ consists exactly of those polynomials that are $\F_q$-linear. Combining equations (\ref{phifq}) and (\ref{phimul}), one gets that for all $a\in\OO$, the polynomial $\phi^L_a$ is $\F_q$-linear:
$$\forall a\in\OO,\;\;\;\phi^L_a\in\Ci\{\tau\},$$
and one also gets that the map $\OO\rightarrow \Ci\{\tau\}$, $a\mapsto\phi^L_a$ is $\F_q$-linear. Thus, it is a morphism of $\F_q$-algebras
\begin{eqnarray}
\phi^L\;:\;\OO & \rightarrow & \Ci\{\tau\} \nonumber\\
a & \mapsto & \phi^L_a. \nonumber
\end{eqnarray}
Let $$D\;:\;\Ci\{\tau\}\rightarrow \Ci$$
be the derivative-at-$0$ map. In other words, $D$ is the $\Ci$-linear map defined by $D(\tau^0)=1$ and $D(\tau^n)=0$ for any $n\geqslant 1$. We have
$$\forall a\in\OO,\;\;\;D\left(\phi^L_a\right)=a.$$

This leads to the general definition of a Drinfel'd module over $\Ci$:

\begin{defi}
Let $\phi\::\:\OO\rightarrow\Ci\{\tau\},\; a\mapsto\phi_a$ \index{phi_a@$\phi_a$ ($\phi$ a Drinfel'd module , $a\in\OO$)}be a morphism of $\F_q$-algebras. Then $\phi$ is a \emph{Drinfel'd module over $\Ci$} if and only if
\begin{mylist}
\item For all $a\in\OO$, $D(\phi_a)=a$.
\item $\phi$ is non-trivial, \ie $\phi$ is not the map $a\mapsto a\tau^0$.
\end{mylist}
\end{defi}

To any lattice $L$ of any rank, we have associated a Drinfel'd module over $\Ci$, which we noted $\phi^L$. The uniformization theorem states that any Drinfel'd module over $\Ci$ comes from a unique lattice. Thus, the map $L\mapsto\phi^L$ is a bijection between lattices and Drinfel'd modules over $\Ci$.\\

The rank of a Drinfel'd module over $\Ci$ is the rank of the associated lattice.

\paragraph*{Action of the ideals.} For any Drinfel'd module $\phi$ over $\Ci$ and any $\aaa\in\Ig_\OO$, we define the polynomial $\phi_\aaa\in\Ci\{\tau\}$ as follows. Let $I_{\aaa,\phi}$ be the left ideal of $\Ci\{\tau\}$ generated by the $\phi_a$, for $a\in\aaa$. One can show that every left ideal of $\Ci\{\tau\}$ is principal, so there exists a unique monic \index{phi_aaa@$\phi_\aaa$ ($\phi$ a Drinfel'd module, $\aaa\in\Ig_\OO$)} $\phi_\aaa\in\Ci\{\tau\}$ such that $I_{\aaa,\phi}=\Ci\{\tau\} \phi_\aaa$.\\

For any Drinfel'd module $\phi$ over $\Ci$ and any $a\in\OO$, we define an element $\mu_\phi(a)\in\Ci^*$ by
\index{muphia@$\mu_\phi(a)$ ($\phi$ a Drinfel'd module, $a\in\OO$)} $$\mu_\phi(a)= \text{leading (highest-degree) coefficient of the polynomial}\;\phi_a.$$
Note that if $\aaa$ is a principal ideal of $\OO$, for any $a\in\OO$ such that $\aaa=a\OO$, we have $$\phi_\aaa=\mu_\phi(a)^{-1}\phi_a.$$
It is easy to see that for any $b\in\OO$, we have $I_{\aaa,\phi}\phi_b\subset I_{\aaa,\phi}$. Thus, for any $b\in\OO$, we have $\phi_\aaa\phi_b\in I_{\aaa,\phi}$, so there is a unique $\phi_b'\in\Ci\{\tau\}$ such that $$\phi_\aaa\phi_b=\phi_b'\phi_\aaa.$$
One shows that the map $b\mapsto\phi_b'$ is a Drinfel'd module over $\Ci$. We will note it\index{aaastarphi@$\aaa*\phi$ ($\phi$ a Drinfel'd module, $\aaa\in\Ig_\OO$)}\index{0@$*$ (as in ``$\aaa*\phi$'')} $\aaa*\phi$. \label{index_actiddrinf} For any two $\aaa,\bbb\in\Ig_\OO$, we have
$$\aaa*(\bbb*\phi) = (\aaa\bbb)*\phi.$$
Thus, $(\aaa,\phi)\mapsto \aaa*\phi$ is an action of $\Ig_\OO$ on the set of all Drinfel'd modules over $\Ci$.\\

Let \index{F_O@$\FF_\OO$}$\FF_\OO$ be the enveloping (``Grothendieck'') group of the abelian semigroup $\Ig_\OO$. The abelian group $\FF_\OO$ may be realized concretely as the group of fractional ideals of $k$ with respect to the Dedekind ring $\OO$. One shows that the action of $\Ig_\OO$ on the set of Drinfel'd modules over $\Ci$ extends to an action of $\FF_\OO$. One also has the equality
\begin{equation}
\label{phi_ab}
\phi_{\aaa\bbb}=(\bbb*\phi)_\aaa\phi_\bbb.
\end{equation}

\paragraph*{Torsion points.} Let $\phi\::\:\OO\rightarrow\Ci\{\tau\}$, $a\mapsto\phi_a$ be a Drinfel'd module over $\Ci$. Remember that $\phi(\Ci)$ is the $\OO$-module that is equal to $\Ci$ as an abelian group, and whose $\OO$-module structure is given by
$$(a,x)\mapsto \phi_a(x).$$
Let \index{phi_Citor@$\phi(\Ci)^\tor$ ($\phi$ a Drinfel'd module)} $\phi(\Ci)^\tor$ denote the $\OO$-torsion submodule of $\phi(\Ci)$. In other words, an element $x\in\phi(\Ci)$ is in $\phi(\Ci)^\tor$ if and only if $\phi_a(x)=0$ for some nonzero $a\in\OO$.\\

 For any $a\in\OO$, let \index{phi_a_croc@$\phi[a]$ ($\phi$ a Drinfel'd module, $a\in\OO$)} $\phi[a]=\ker\phi_a$.  For any $\aaa\in\Ig_\OO$, let \index{phi_aaa_croc@$\phi[\aaa]$ ($\phi$ a Drinfel'd module, $\aaa\in\Ig_\OO$)} $\phi[\aaa]=\ker\phi_\aaa$.
Under the bijection given by equation (\ref{e_L_isom}), the sets $\phi(\Ci)^\tor$, $\phi[a]$ and $\phi[\aaa]$ are identified to the following subsets of $\Ci/L$:
\begin{eqnarray}
e_L^{-1}(\phi(\Ci)^\tor) & = & kL/L \nonumber\\
\forall a \in\OO-\{0\},\;\;e_L^{-1}(\phi[a]) & = & a^{-1}L/L \nonumber\\
\forall \aaa \in\Ig_\OO,\;\;e_L^{-1}(\phi[\aaa]) & = & \aaa^{-1}L/L, \nonumber
\end{eqnarray}
where $\aaa^{-1}$ is the inverse of $\aaa$ as a fractional ideal with respect to $\OO$, \ie
$$\aaa^{-1}\;=\;\{x\in k,\;\;\;x\aaa\subset\OO\}.$$
The following equalities follow from the definitions:
\begin{eqnarray}
\forall a\in\OO,\;\;\;\phi[a] &=&\phi[a\OO],\nonumber\\
\forall \aaa\in\Ig_\OO,\;\;\;\phi[\aaa]&=&\bigcap_{a\in\aaa} \phi[a],\nonumber\\
\phi(\Ci)^\tor&=&\bigcup_{a\in\OO} \phi[a],\nonumber\\
\phi(\Ci)^\tor&=&\bigcup_{\aaa\in\Ig_\OO} \phi[\aaa],\nonumber
\end{eqnarray}
and
\begin{equation}
\label{xr5}
\forall \aaa,\bbb\in\Ig_\OO,\;\;\;\aaa\mid\bbb \Leftrightarrow \phi[\aaa] \subset \phi[\bbb].
\end{equation}
One also checks that, for all $\aaa,\bbb \in \Ig_\OO$,
\begin{equation}
\label{xr6}\phi[\aaa]\bigcap\phi[\bbb]=\phi[\aaa+\bbb],
\end{equation}
\begin{equation}
\label{xr7}\phi[\aaa]+\phi[\bbb]=\phi\left[\aaa\bigcap\bbb\right].
\end{equation}
We have 
\begin{equation}
\label{cardphi}
\forall\aaa\in\Ig_\OO,\;\;\;\Card\phi[\aaa]=(\N\aaa)^r,
\end{equation}
where $r$ is the rank of $\phi$ and \index{Naaa@$\N\aaa$ ($\aaa\in\Ig_\OO$)} $\N\aaa$ is the absolute norm of $\aaa$, \ie $\N\aaa$ is the cardinal of $\OO/\aaa$.\\

Let $\aaa\in\Ig_\OO$. By construction, $\phi_\aaa$ is an $\OO$-module morphism
$$\phi_\aaa\;:\;\phi(\Ci)\rightarrow (\aaa*\phi)(\Ci).$$
For any $\bbb\in\Ig_\OO$, Let $\phi_\aaa\vert_{\phi[\bbb]}$ denote the restriction of $\phi_\aaa$ to $\phi[\bbb]$.

\begin{lemma}
\label{xr8}
Let $\phi$ be a Drinfel'd module over $\Ci$. Let $\aaa,\bbb\in\Ig_\OO$. Let $\ddd=\aaa+\bbb$ be the gcd of $\aaa$ and $\bbb$. We have:
\begin{eqnarray}
\Ker\left(\phi_\aaa\vert_{\phi[\bbb]}\right) & = & \phi[\ddd], \nonumber\\
\im\left(\phi_\aaa\vert_{\phi[\bbb]}\right) & = & (\aaa*\phi)[\ddd^{-1}\bbb]. \nonumber
\end{eqnarray} 
\end{lemma}
\begin{proof}
First equality: we have $\Ker\left(\phi_\aaa\vert_{\phi[\bbb]}\right)=\phi[\aaa]\bigcap\phi[\bbb],$ so the result 
follows from equation (\ref{xr6}).\n

\noindent Second equality: let $r$ denote the rank of $\phi$. We have
$$\Card\left(\im(\phi_\aaa\vert_{\phi[\bbb]})\right)=\Card(\phi[\bbb])/\Card\left(\Ker\left(\phi_\aaa\vert_{\phi[\bbb]}
\right)\right) = \Card(\phi[\bbb])/\Card(\phi[\ddd]) = (\N\bbb)^r/(\N\ddd)^r$$ and
$$\Card((\aaa*\phi)[\ddd^{-1}\bbb])=\N(\ddd^{-1}\bbb)^r,$$
so the two cardinals are equal, so it is enough to show one inclusion. Let 
$x\in\im\left(\phi_\aaa\vert_{\phi[\bbb]}\right)$. It is enough to show 
that $(\aaa*\phi)_{\ddd^{-1}\bbb}(x)=0$.  Let $y\in\phi[\bbb]$ such that $\phi_\aaa(y)=x$. Let $\ccc=\aaa\bigcap\bbb$ be the lcm. We have $\ddd^{-1}\bbb=\aaa^{-1}\ccc$. 
But$$(\aaa*\phi)_{\aaa^{-1}\ccc}(x)=(\aaa*\phi)_{\aaa^{-1}\ccc}(\phi_\aaa(y)),$$so, by equation (\ref{phi_ab}), we 
get$$(\aaa*\phi)_{\aaa^{-1}\ccc}(x)=\phi_\ccc(y).$$But $\bbb\mid\ccc$ and $y\in\phi[\bbb]$, so $y\in\phi[\ccc]$, so \[(\aaa*\phi)_{\aaa^{-1}\ccc}(x)=0.\qedhere \]
\end{proof}

\begin{coro}
\label{xr9}
Let $\phi$ be a Drinfel'd module over $\Ci$. Let $\aaa,\bbb\in\Ig_\OO$. For all $\lambda\in\phi[\bbb]$, there exists $\mu\in(\aaa^{-1}*\phi)[\aaa\bbb]$ such that 
$$(\aaa^{-1}*\phi)_\aaa(\mu)=\lambda.$$\end{coro}
\begin{proof}
Let $ \psi=\aaa^{-1}*\phi$. Let $\bbb_2=\aaa\bbb$. Let $\ddd_2=\aaa$, so that $\ddd_2$ is the gcd of $\aaa$ and $\bbb_2$. By Lemma \ref{xr8}, we 
have$$\im\left(\psi_\aaa \vert_{\psi[\bbb_2]}\right)  =  \phi[\ddd_2^{-1}\bbb_2],$$
so
\[\im\left(\psi_\aaa \vert_{\psi[\aaa\bbb]}\right)  =  \phi[\bbb].\qedhere \]
\end{proof}

\begin{coro}
\label{idsurj}
Let $\phi$ be a Drinfel'd module over $\Ci$. For all $\aaa\in\Ig_\OO$, the map $$(\aaa^{-1}*\phi)_\aaa\;:\;(\aaa^{-1}*\phi)(\Ci)^\tor\rightarrow \phi(\Ci)^\tor$$ is 
surjective.\end{coro}

\subsection{Hayes' explicit class field theory}
\label{xr10}

In this subsection, we review D. R. Hayes' explicit class field theory for function fields, in terms of sign-normalized rank one Drinfel'd modules. We follow \cite{Hay92}, Part II, and \cite{Gos91}, Chapter VII. Recall that $k_\infty$ is the completion of $k$ at $\infty$. Let \index{F_inf@$\F_\infty$} $\F_\infty$ denote the constant subfield of $k_\infty$. The field  $\F_\infty$ is a finite extension of $\F_q$, and its degree is equal to the degree of the place $\infty$.

\begin{defi}
A \emph{sign function} on $k_\infty^*$ is a group morphism $\sgn\::\:k_\infty^*\rightarrow \F_\infty^*$ which induces the identity map on $\F_\infty^*$.
\end{defi}

Let us choose once and for all a sign-function \index{sgn@$\sgn$} $\sgn$ (by \cite{Hay92}, Corollary 12.2, the number of possible choices is equal to the cardinal of  $\F_\infty^*$). We let $\sgn(0)=0$, so that $\sgn$ becomes a function $k_\infty\rightarrow \F_\infty$.

\begin{defi}
A Drinfel'd module $\phi$ over $\Ci$ is said to be $\sgn$-\emph{normalized} if there exists an element $\sigma\in\Gal(\F_\infty/\F_q)$ such that
$$\forall a\in\OO,\;\;\;\mu_\phi(a) = \sigma (\sgn (a)).$$
\end{defi}

Let us now focus on the case of Drinfel'd modules of rank one:

\begin{defi}
Let \index{H_sgn@$\Hay$}$\Hay$ denote the set of $\sgn$-normalized rank one Drinfel'd modules over $\Ci$. The elements of $\Hay$ are also called \emph{Hayes modules} (for the triple $(k,\infty,\sgn)$).
\end{defi}

\begin{prop}
$\Hay$ is a finite set, and its cardinal \index{h_sgn@$h(\sgn)$}$h(\sgn)$ is given by:
$$h(\sgn) = \frac{\Card\F_\infty^*}{\Card\F_q^*}\cdot h(\OO),$$
where $h(\OO)$ is the class number of the Dedekind ring $\OO$.
\end{prop}
\begin{proof}
See \cite{Hay92}, Corollary 13.4.
\end{proof}

\begin{prop}
For any $\phi\in\Hay$ and any $\aaa\in\FF_\OO$, we have $\aaa*\phi\in\Hay$. Thus, $\FF_\OO$ acts on $\Hay$.
\end{prop}
\begin{proof}
See \cite{Hay92}, page 22.
\end{proof}

\begin{defi}
\label{Hplus}
Let $\phi\in\Hay$, and let $y\in\OO-\F_q$ (recall that $\F_q$ denotes the constant subfield of $k$). Let \index{H_plus@$H^+$} $H^+$ be the field generated over $k$ by the coefficients of $\phi_y$.
\end{defi}
One shows (see \cite{Hay92}, page 23) that $H^+$ does not depend on the choice of $\phi$ and $y$.

\begin{prop}
The extension $H^+/k$ is finite, abelian, and unramified away from $\infty$.
\end{prop}
\begin{proof}
See \cite{Hay92}, Propositions 14.1 and 14.4.
\end{proof}
One shows (see \cite{Hay92}, \S 15) that $H^+$ contains a subfield $H$ which plays the rôle of the Hilbert class field for the pair $(k,\infty)$.\\

Here is a concrete picture of the Galois group $\Gal(H^+/k)$. First, let $\mathcal{P}^+_\OO$ be the following subgroup of $\FF_\OO$ :
$$\mathcal{P}^+_\OO\;=\;\{x\OO\;:\;x\in k,\;\sgn(x)=1\}.$$
We then have the following proposition:
\begin{prop}
The Artin map $(\cdot,H^+/k)$ induces an isomorphism from $\FF_\OO/\mathcal{P}^+_\OO$ to $\Gal(H^+/k)$.
\end{prop}
\begin{proof}
See \cite{Hay92}, Theorem 14.7.
\end{proof}

The Galois group $\Gal(H^+/k)$ acts on $\Hay$ by $(\sigma,\phi)\mapsto\sigma\phi$, where $\sigma\phi$ is defined by $(\sigma\phi)_a = \sigma(\phi_a)$ for all $a\in\OO$ (one checks that $\sigma\psi\in\Hay$).

\begin{theo}
\label{hayes1}
For any $\aaa\in\Ig_\OO$, if $\sigma_\aaa=(\aaa,H^+/k)\in\Gal(H^+/k)$ denotes the Artin automorphism of $H^+$ associated to $\aaa$, then we have
$$\forall\phi\in\Hay,\;\;\;\sigma_\aaa\phi=\aaa*\phi.$$
The set $\Hay$ is principal homogeneous under the action of $\Gal(H^+/k)$.
\end{theo}
\begin{proof}
See \cite{Hay92}, Theorems 13.8 and 14.7.
\end{proof}

\begin{defi}
\label{Kc}
For any $\phi\in\Hay$, let \index{K@$K$} $K$ denote the field generated over $H^+$ by the elements of $\phi(\Ci)^\tor$.  For any $\ccc\in\Ig_\OO$, let \index{K_fff@$K_\ccc$ ($\ccc\in\Ig_\OO$)}$K_\ccc$ denote the field generated over $H^+$ by the elements of $\phi[\ccc]$.
\end{defi}
One shows (see \cite{Hay92}, page 28) that $K$ and $K_\ccc$ are independent of the choice of $\phi$. The extension $K_\ccc/k$ is called the \emph{narrow ray class extension modulo} $\ccc$. By construction, we have
$$K\:=\:\bigcup_{\ccc\in\Ig_\OO} K_\ccc.$$
\begin{theo}
\label{hayes2}
For any $\ccc\in\Ig_\OO$, the extension $K_\ccc/k$ is finite, abelian, and unramified away from $\infty$ and the prime divisors of $\ccc$. Moreover, $K_\ccc$ contains the ray class field of $k$ of conductor $\ccc$ totally split at $\infty$.
For any $\aaa\in\Ig_\OO$ prime to $\ccc$, if $\sigma_\aaa=(\aaa,K_\ccc/k)\in\Gal(K_\ccc/k)$ denotes the Artin automorphism of $K_\ccc$ associated to $\aaa$, then we have
$$\forall\phi\in\Hay,\;\;\forall\lambda\in\phi[\ccc],\;\;\;\sigma_\aaa\lambda=\phi_\aaa(\lambda).$$
\end{theo}
\begin{proof}
See \cite{Hay92}, page 28, or \cite{Hay79}, section 8.
\end{proof}

In particular, this shows that
$$k^\abi\subset K\subset k^\ab,$$ where $k^\abi$ is the maximal abelian extension of $k$ that is totally split at $\infty$.\\

Let us give a concrete picture of the Galois group $\Gal(K_\ccc/k)$, for $\ccc\in\Ig_\OO$. Let $\FF_\OO(\ccc)$ denote the subgroup of $\FF_\OO$ of all fractional ideals that are prime to $\ccc$, and let
$$\mathcal{P}^+_\OO(\ccc)\;=\;\{x\OO\;:\;x\in k,\;\sgn(x)=1,\; x\equiv 1\;\mathrm{mod}\;\ccc\}.$$
We then have the following proposition:
\begin{prop}
The Artin map $(\cdot,K_\ccc/k)$ induces an isomorphism from $\FF_\OO(\ccc)/\mathcal{P}^+_\OO(\ccc)$ to $\Gal(K_\ccc/k)$.
\end{prop}
\begin{proof}
See \cite{Hay92}, page 28.
\end{proof}
Moreover, the Galois group $\Gal(K_\ccc/H^+)$ has an even simpler description: one can check (loc. cit.) that it is isomorphic to the group of invertible elements in $\OO/\ccc$.

\section{Construction of the $C^*$-dynamical system $(C_{k,\infty}, (\sigma_t))$}

\subsection {The space $X$ of characters}\label{Xsubsec}

For any $\phi\in\Hay$, let \index{X_phi@$X_\phi$}$X_\phi$ be the dual group of the discrete abelian torsion group $\phi(\Ci)^\tor$. Thus, an element of $X_\phi$ is a character of $\phi(\Ci)^\tor$. The group $X_\phi$ is profinite:
$$X_\phi=\lim_{\leftarrow\aaa} \widehat{\phi[\aaa]},$$
where $\aaa$ runs over $\Ig_\OO$ ordered by divisibility. Let $X$ be the (disjoint) union of the $X_\phi$:\index{X@$X$}
$$X\;=\; \bigcup_{\phi\in\Hay} X_\phi.$$

Note that the elements of $X$ are reminiscent of characters in \cite{HarLei97} and of $\QQ$-lattices (or $k$-lattices) in 
\cite{ConMar04} and \cite{ConMarRam05}, \cite{ConMarRam06}.

\begin{lemma}
\label{UUp}
For any character $\chi\in X$, we have
$$\im \chi \subset \UU_p,$$
where $\UU_p$ is the group of $p$-th roots of unity in $\CC$.
\end{lemma}
\begin{proof}
Recall that for any $\phi\in\Hay$, as a group, $\phi(\Ci)$ is equal to $\Ci$, which is a field of characteristic $p$. Thus, for all $\lambda\in\phi(\Ci)^\tor$, we have $\chi(\lambda)^p =\chi(p\lambda) =\chi(0) =1$.
\end{proof}

\begin{lemma}
\label{Xcompact}
$X$ is compact (and Hausdorff).
\end{lemma}
\begin{proof}
For any $\phi\in\Hay$, the group $X_\phi$ is profinite, hence compact. As $\Hay$ is finite, $X$ is compact.
\end{proof}

We define an action of $\Ig_\OO$ on $X$  by
\index{chi_a@$\chi^\aaa$ ($\chi\in X$, $\aaa\in\Ig_\OO$)} 
\begin{equation}
\label{xr11}
\forall \aaa\in\Ig_\OO,\;\;\forall \phi\in\Hay,\;\; \forall \chi\in X_\phi,\;\;
\chi^\aaa=\chi\circ(\aaa^{-1}*\phi)_\aaa.
\end{equation}
Recall that $(\aaa^{-1}*\phi)_\aaa$ is a map from $(\aaa^{-1}*\phi)(\Ci)$ to $\phi(\Ci)$. Thus, if $\chi\in X_\phi$, then $\chi^\aaa\in X_{\aaa^{-1}*\phi}$. Note that equation (\ref{phi_ab}) guarantees that this is a semigroup action of $\Ig_\OO$.\\

The exponent notation ($\chi^\aaa$) is inspired by what happens with characters of $\QQ/\ZZ$. These characters may be composed with the map $\phi_n\::\:x\mapsto nx$, for any $n\in\NN^*$. By definition of a character, we have $\chi\circ\phi_n=\chi^n$. In our case, $\NN^*$ is replaced by $\Ig_\OO$ and the maps $\phi_n$ are replaced by the $\phi_\aaa$.\\

We define an 
action of $\Gal(K/k)$ on $X$  by
\index{sigma_chi@$\sigma \chi$ ($\sigma\in\Gal(K/k)$, $\chi\in X$)}
\begin{equation}
\label{galX}
\forall \sigma\in\Gal(K/k),\;\;\forall\chi\in X,\;\;\;\sigma\chi=\chi\circ\sigma,
\end{equation}

One checks that the actions of $\Gal(K/k)$ and of $\Ig_\OO$ on $X$ commute with one another.

\begin{lemma}
\label{idinj}
For all $\aaa\in\Ig_\OO$, the map $X\rightarrow X$, $\chi\mapsto \chi^\aaa$ is injective.
\end{lemma}
\begin{proof}
Let $\chi_1,\chi_2\in X$ such that $\chi_1^\aaa=\chi_2^\aaa$. For $i=1,2$, let $\phi^i$ be such that $\chi_i\in X_{\phi^i}$. By definition, we have $\chi_i^\aaa\in X_{\aaa^{-1}*\phi^i}$, so
$\aaa^{-1}*\phi^1=\aaa^{-1}*\phi^2$, so $\phi^1=\phi^2$. Let $\phi=\phi^1=\phi^2$. We have $$\chi_1\circ(\aaa^{-1}*\phi)_\aaa=\chi_2\circ(\aaa^{-1}*\phi)_\aaa.$$
Corollary \ref{idsurj} then shows that $\chi_1=\chi_2$.
\end{proof}

\begin{coro}
\label{xr12}
Let $\aaa_1,\aaa_2,\bbb_1,\bbb_2\in\Ig_\OO$ be such that $\aaa_1^{-1}\aaa_2=\bbb_1^{-1}\bbb_2$.
\begin{mylist}
\item Let $\chi_1,\chi_2\in X$. We have:
$$\chi_1^{\aaa_1}=\chi_2^{\aaa_2}\;\;\Leftrightarrow\;\;\chi_1^{\bbb_1}=\chi_2^{\bbb_2}.$$
\item Let $\chi_1,\chi_2,\chi_3\in X$. We have
$$\chi_1^{\aaa_1}=\chi_2^{\aaa_2}\;\;\textup{and}\;\;\chi_3^{\bbb_1}=\chi_2^{\bbb_2}\;\;\Rightarrow\;\; \chi_1=\chi_3.$$
\end{mylist}
\end{coro}
\begin{proof}
Let us first prove (1). Suppose that $\chi_1^{\aaa_1}=\chi_2^{\aaa_2}$. We have $\chi_1^{\aaa_1\bbb_2}=\chi_2^{\aaa_2\bbb_2}$. But $\aaa_2\bbb_1=\aaa_1\bbb_2$, so $\chi_1^{\aaa_2\bbb_1}=\chi_2^{\aaa_2\bbb_2}$, so
$$(\chi_1^{\bbb_1})^{\aaa_2}=(\chi_2^{\bbb_2})^{\aaa_2},$$
so Lemma \ref{idinj} gives $\chi_1^{\bbb_1}=\chi_2^{\bbb_2}$, which proves one implication, and the other implication follows by swapping $\aaa_i$ with $\bbb_i$ for $i=1,2$.\\

Let us now prove (2). We have
$$\chi_1^{\aaa_1\bbb_1}= \chi_2^{\aaa_2\bbb_1} =\chi_2^{\aaa_1\bbb_2}=\chi_3^{\aaa_1\bbb_1},$$
so Lemma \ref{idinj} gives $\chi_1=\chi_3$.
\end{proof}

Corollary \ref{xr12} allows to extend the action of $\Ig_\OO$ on $X$ to a partially defined action of $\FF_\OO$, as follows:

\begin{defi}
For any $\chi\in X$, let \index{F_chi@$\FF_\chi$ ($\chi\in X$)} $\FF_\chi$ denote the set of all $\ccc\in\FF_\OO$ such that there exists $\chi_1\in X$ such that
\begin{equation}
 \label{xr13}
\chi_1^{\aaa_1}=\chi^{\aaa_2}.
\end{equation}
for some $\aaa_1,\aaa_2\in\Ig_\OO$ such that $\ccc=\aaa_1^{-1}\aaa_2$. By Corollary \ref{xr12} (1),
the existence of $\chi_1$ only depends on $\chi$ and $\ccc$, and does not depend on the choice of $\aaa_1,\aaa_2\in\Ig_\OO$ such that $\ccc=\aaa_1^{-1}\aaa_2$. By Corollary \ref{xr12} (2), the character $\chi_1$, when it exists, is uniquely determined by $\chi$ and $\ccc$. When $\ccc\in\FF_\chi$, we define a character $\chi^\ccc$ by \index{chi_c@$\chi^\ccc$ ($\chi\in X$, $\ccc\in\FF_\OO$)}
$$\chi^\ccc=\chi_1.$$
\end{defi}

The partially defined map $\FF_\OO\times X\rightarrow X$, $(\ccc,\chi)\mapsto \chi^\ccc$, should be regarded as a \emph{partially defined} group action of $\FF_\OO$ on $X$. For any $(\ccc,\chi)$, the character $\chi^\ccc$ is defined if, and only if $\ccc\in\FF_\chi$. For any $\ccc_1,\ccc_2$ in $\FF_\chi$, if $\ccc_1\ccc_2\in\FF_\chi$, one checks that $\chi^{\ccc_1\ccc_2}=(\chi^{\ccc_1})^{\ccc_2}.$ Of course, when $\ccc\in\Ig_\OO$, the character $\chi^\ccc$ is just the one that was defined in equation (\ref{xr11}).\\

For any $\chi$, we have $\Ig_\OO\subset\FF_\chi$. Characters $\chi\in X$ for which this inclusion is an equality ($\FF_\chi=\Ig_\OO$) will be called \emph{admissible}, and will play an important rôle later (see subsection \ref{admsubsec}).\n

Note that we obviously have
\begin{equation}
\label{xr14}
\forall \chi\in X,\;\;\;\forall \aaa\in\Ig_\OO,\;\;\;\FF_{\chi^\aaa}=\aaa^{-1}\FF_\chi.
\end{equation}

\begin{lemma}
\label{xr15}
Let $\chi\in X$. Let $\phi\in\Hay$ such that $\chi\in X_\phi$. For any $\aaa\in\Ig_\OO$, we have
$$\aaa^{-1}\in\FF_\chi\:\:\:\Longleftrightarrow\:\:\:\forall\lambda\in\phi[\aaa],\:\:\chi(\lambda)=1.$$
When this is the case, the character $\chi^{\aaa^{-1}}$ is given by
$$\forall\lambda\in(\aaa*\phi)(\Ci)^\tor,\:\:\:\chi^{\aaa^{-1}}(\lambda)=(\N\aaa)^{-1}\sum_{\phi_\aaa(\mu)=\lambda} 
\chi(\mu).$$
\end{lemma}
\begin{proof}
If $\aaa^{-1}\in\FF_\chi$, then there exists $\chi_1\in X$ such that $\chi=\chi_1^\aaa$. Thus, $\forall\lambda\in\phi[\aaa]$, 
we have $\chi(\lambda)=\chi_1(\phi_\aaa(\lambda))$, but $\phi_\aaa(\lambda)=0$, so $\chi(\lambda)=1$.\n

Now suppose that $\forall\lambda\in\phi[\aaa]$, $\chi(\lambda)=1$. For all $\lambda\in (\aaa*\phi)(\Ci)^\tor$, set
$$\chi_1(\lambda)=(\N\aaa)^{-1}\sum_{\phi_\aaa(\mu)=\lambda} \chi(\mu).$$
Let us show that this defines a character $\chi_1$ of $(\aaa*\phi)(\Ci)^\tor$. Let $\lambda\in(\aaa*\phi)(\Ci)^\tor$. By Lemma \ref{idsurj}, 
there exists $\mu_1\in\phi(\Ci)^\tor$ such that $\phi_\aaa(\mu_1)=\lambda$. We have 
$$\chi_1(\lambda)=(\N\aaa)^{-1}\sum_{\mu_0\in\phi[\aaa]} \chi(\mu_0+\mu_1)=(\N\aaa)^{-1}\left(\sum_{\mu_0\in\phi[\aaa]} 
\chi(\mu_0)\right)\chi(\mu_1).$$But we have $\chi(\mu_0)=1$ for all $\mu_0\in\phi[\aaa]$, and, by Equation (\ref{cardphi}), 
we have $\Card(\phi[\aaa])=\N\aaa$. Thus, we get
$$\forall\lambda\in (\aaa*\phi)(\Ci)^\tor,\:\:\forall 
\mu_1\:\:\text{such that}\:\:\phi_\aaa(\mu_1)=\lambda,\:\:\:\:\chi_1(\lambda)=\chi(\mu_1).$$
Now let $\lambda'\in 
(\aaa*\phi)(\Ci)^\tor$ and $\mu_1'$ such that $\phi_\aaa(\mu'_1)=\lambda'$. We have$$\lambda+\lambda'=\phi_\aaa(\mu_1)+ 
\phi_\aaa(\mu_1')=\phi_\aaa(\mu_1+\mu_1'),$$so $$\chi_1(\lambda+\lambda')=\chi_1(\phi_\aaa(\mu_1+\mu_1')) = 
\chi(\mu_1+\mu_1')=\chi(\mu_1)\chi(\mu_1'),$$
so $$\chi_1(\lambda+\lambda')=\chi_1(\lambda)\chi_1(\lambda'),$$
so $\chi_1(\lambda)^p=\chi_1(p\lambda)=\chi_1(0)=1,$ so
$$\forall\lambda\in (\aaa*\phi)(\Ci)^\tor,\;\;\;\chi_1(\lambda)\in\UU_p.$$
Hence, $\chi_1$ is a group morphism $(\aaa*\phi)(\Ci)^\tor\rightarrow\UU_p$, so $\chi_1\in X$, and we have by construction $\chi_1^\aaa=\chi$. Thus, we have $\aaa^{-1}\in\FF_\chi$ and $\chi^{\aaa^{-1}}=\chi_1$.\end{proof}

\begin{lemma}
\label{xr16}
For all $\chi\in X$, for all $\aaa,\bbb\in\Ig_\OO$ relatively prime, we have
$$\aaa^{-1}\bbb\in\FF_\chi\:\:\:\Longleftrightarrow\:\:\:\aaa^{-1}\in\FF_\chi.$$
\end{lemma}
\begin{proof}
Let $\phi$ be such that $\chi\in X_\phi$. We have $\aaa^{-1}\bbb\in\FF_\chi\:\Leftrightarrow\:\aaa^{-1}\in\FF_{\chi^\bbb}$. Lemma 
\ref{xr15} applied to $\chi^\bbb$ thus 
gives$$\aaa^{-1}\bbb\in\FF_\chi\:\:\:\Longleftrightarrow\:\:\:\forall\lambda\in(\bbb^{-1}*\phi)[\aaa],\:\:\chi((\bbb^{-1}*\phi)_\bbb(\lambda))
=1.$$But, as $\aaa$ and $\bbb$ are relatively prime, by Lemma \ref{xr8}, the map $\lambda\mapsto\bbb\lambda$ is a 
bijection from $(\bbb^{-1}*\phi)[\aaa]$ onto $\phi[\aaa]$. Thus we 
get$$\aaa^{-1}\bbb\in\FF_\chi\:\:\:\Longleftrightarrow\:\:\:\forall\lambda\in\phi[\aaa],\:\:\chi(\lambda)=1,$$and, by Lemma 
\ref{xr15}, this is equivalent to $\aaa^{-1}\in\FF_\chi$.\end{proof}
\begin{lemma}
\label{xr17}
For all $\chi\in X$, for all $\aaa,\bbb\in\Ig_\OO$ relatively prime, we have
$$(\aaa\bbb)^{-1}\in\FF_\chi\:\:\:\Longleftrightarrow\:\:\: \aaa^{-1}\in\FF_\chi\:\:\textup{and}\:\:\bbb^{-1}\in\FF_\chi.$$
\end{lemma}
\begin{proof}
Let $\phi$ be such that $\chi\in X_\phi$. By Lemma \ref{xr15}, the statement that we want to prove is equivalent to the following:
\begin{equation}
\label{xr18}
\forall\lambda\in\phi[\aaa\bbb],\:\:\chi(\lambda)=1\:\:\:\Longleftrightarrow\:\:\:\forall\lambda\in\phi[\aaa],\:\:
\chi(\lambda)=1\:\:\textup{and}\:\:\forall\lambda\in\phi[\bbb],\:\:\chi(\lambda)=1.
\end{equation}
By Equations (\ref{xr6}) and (\ref{xr7}), as $\aaa$ 
and $\bbb$ are relatively prime, we have  $$\phi[\aaa\bbb]=\phi[\aaa]\bigoplus\phi[\bbb],$$
so, for any $\lambda\in\phi[\aaa\bbb]$, there exists a unique pair $(\lambda_1,\lambda_2)\in\phi[\aaa]\times\phi[\bbb]$ such that $\lambda=\lambda_1+\lambda_2$. We have $\chi(\lambda)=\chi(\lambda_1)\chi(\lambda_2)$, so equation (\ref{xr18}) follows.
\end{proof}

\subsection{Construction of the groupoid $\GG$ and of the dynamical system $(C_{k,\infty},\,(\sigma_t))$}

Let $\GG$ be the following subset of $X\times\FF_\OO$:
\index{GG@$\GG$}
$$\GG\:=\:\left\{ (\chi,\ccc)\in X\times \FF_\OO,\:\: \ccc\in\FF_\chi\right\}.$$
We turn $\GG$ into a groupoid by endowing it with the groupoid law
$$(\chi_1,\ccc_1)\circ (\chi_2,\ccc_2) = (\chi_2, \ccc_1 \ccc_2)\:\:\:\:\:\text{if}\:\:\:\chi_1=\chi_2^{\ccc_2}$$
and the inverse map
$$(\chi,\ccc)^{-1}=(\chi^\ccc, \ccc^{-1}).$$
One checks that, under the identification $\GG^{(0)}= X\times\{1\}\simeq X$, the range and source maps $r$ and $s$ are respectively given by $r(\chi, \ccc) = \chi^\ccc$ and $s(\chi,\ccc)=\chi$.\\

The abelian group $\FF_\OO$ is endowed with the discrete topology. The groupoid $\GG$ is endowed with its topology as a subset of 
$X\times\FF_\OO$.

\begin{lemma}
$\GG$ is a locally compact groupoid.
\end{lemma}
\begin{proof}
$X\times\FF_\OO$ is locally compact by Lemma \ref{Xcompact}, and $\GG$ is a closed subset of it, so it is also locally 
compact. It is clear that the composition and inverse maps are continuous, so this is a locally compact 
groupoid.\end{proof}

The $C^*$-algebra $C_{k,\infty}$ that was advertised in the introduction of this paper is the maximal\footnote{Actually, it coincides with the reduced $C^*$-algebra because $\FF_\OO$ is an abelian group, but we will not need that in this paper.} $C^*$-algebra of the groupoid $\GG$. Let us quickly explain what that means.\\

\noindent For $\chi\in X$, let $\GG_\chi$ denote the fiber of $s$ above $\chi$, that is, \index{GG_chi@$\GG_\chi$ ($\chi\in X$)}
$$\GG_\chi=\{\chi\}\times\FF_\chi,$$
so $\GG_\chi$ is discrete and is in bijection with $\FF_\chi$.\n

Let \index{C_cGG@$C_c(\GG)$} $C_c(\GG)$ denote the convolution algebra of continuous maps $\GG\rightarrow\CC$ with compact support, where the 
convolution product is given by\begin{equation}\label{convol}
(f_1 f_2)(g)=\sum_{g_1\circ g_2=g} f_1(g_1) f_2(g_2).
\end{equation}
$C_c(\GG)$ is endowed with the involution $f\mapsto f^*$ defined by
$$f^*(g)=\overline{f(g^{-1})}.$$
For any $\chi\in X$, we define a $*$-representation of $C_c(\GG)$ on the Hilbert space $\ell^2(\GG_\chi)$ by
\begin{equation}
\label{regul}
\index{pi_chi@$\pi_\chi$ ($\chi\in X$)}
\forall f\in C_c(\GG),\:\forall\xi\in\ell^2(\GG_\chi),\:\:\:(\pi_\chi(f)\xi)(g)= \sum_{g_1\circ g_2=g} f(g_1) \xi (g_2).
\end{equation}
In other words, $\pi_\chi$ is the left regular representation on $\ell^2(\GG_\chi)$. Let us define a $C^*$-norm $\norm{\cdot}$ on $C_c(\GG)$ by $$\norm{f}\:=\:\sup_{\pi}\: \norm{\pi(f)},$$where $\pi$ runs over all 
$*$-representations of $C_c(\GG)$. The completion $C^*(\GG)$ of $C_c(\GG)$ under $\norm{\cdot}$ is a $C^*$-algebra, called the \emph{maximal $C^*$-algebra} of the groupoid $\GG$. For more details about groupoid $C^*$-algebras, see 
Renault \cite{Ren80}, Khoshkam and Skandalis \cite{KhoSka02}, or Connes \cite{Con94}, Chapter II, \S 5.

\begin{defi} We define the $C^*$-algebra \index{C_kinf@$C_{k,\infty}$} $C_{k,\infty}$ by letting $$C_{k,\infty}\:=\:C^*(\GG).$$
\end{defi}

By definition, any $*$-representation $\pi$ of $C_c(\GG)$ extends uniquely to a representation of $C_{k,\infty}$, which we will still note $\pi$.

\begin{lemma}
\label{xr19}
For any $*$-automorphism $\sigma$ of $C_c(\GG)$, there exists an unique extension of $\sigma$ to a $*$-automorphism of 
$C_{k,\infty}$.\end{lemma}
\begin{proof}
For any $*$-automorphism $\sigma$ of $C_c(\GG)$ and any $*$-representation $\pi$ of $C_c(\GG)$, note that $\pi\circ\sigma$ is a $*$-representation of $C_c(\GG)$. Thus, by definition of the norm $\norm{\cdot}$, $\sigma$ is an isometry: for all $f\in C_c(\GG)$, we have $\norm{\sigma(f)} = \norm{f}$. The result then follows easily.
\end{proof}

For any $g=(\chi,\ccc)\in\GG$, put $\N g=\N\ccc$, where $\N\ccc$ is the absolute norm of the fractional ideal $\ccc$, defined by \index{Nccc@$\N\ccc$ ($\ccc\in\FF_\OO$)} $\N\ccc=(\N\aaa)^{-1}\N\bbb$ for any $\aaa,\bbb\in\Ig_\OO$ such that $\ccc=\aaa^{-1}\bbb$.\\

Let us define a one parameter 
$*$-automorphism group $(\sigma_t)_{t\in\RR}$ of $C_c(\GG)$ by: \index{sigma_t@$\sigma_t$ ($t\in\RR$)} $$\forall t\in\RR,\:\:\forall f\in C_c(\GG),\:\:\forall 
g\in\GG,\:\:\:(\sigma_t(f))(g)=(\N g)^{it}f(g).$$

\begin{defi}
We will still note \index{sigma_t@$\sigma_t$ ($t\in\RR$)} $\sigma_t$ the unique extension (given by Lemma \ref{xr19}) of $\sigma_t$ to an automorphism of $C_{k,\infty}$.\end{defi}

It remains to check that the pair $(C_{k,\infty},\,(\sigma_t))$ is a $C^*$-dynamical system in the sense of \cite{BraRob}, \ie that the flow $(\sigma_t)$ is strongly continuous, \ie that for any $f\in C_{k,\infty}$, the map $t\mapsto \sigma_t(f)$ is continuous.

\begin{lemma}
The flow $(\sigma_t)$ on $C_{k,\infty}$ is strongly continuous.
\end{lemma}
\begin{proof}
Let $f\in C_{k,\infty}$. Let us show that the map $t\mapsto \sigma_t(f)$ is continuous. Let $\varepsilon >0$. It is enough to show that when $\abs{t}$ is small enough, we have $\norm{f-\sigma_t(f)}<\varepsilon$. Let $f'\in C_c(\GG)$ be such that $\norm{f-f'}<\varepsilon/3$. Like any $*$-automorphism, $\sigma_t$ is an isometry, so we have $\norm{\sigma_t(f) - \sigma_t(f')}=\norm{\sigma_t(f - f')}=\norm{f - f'}<\varepsilon/3$, so it is enough to show that when $\abs{t}$ is small enough, we have $\norm{f'-\sigma_t(f')}<\varepsilon/3$. For any $\ddd\in\FF_\OO$, define a function $f_\ddd'\in C_c(\GG)$ by
$$\forall (\chi,\ccc)\in\GG,\;\;\;f_\ddd'(\chi,\ccc)=\left\lbrace\begin{matrix} f'(\chi,\ccc) & \text{if} & \ccc=\ddd,\\ 0 &\text{if} & \ccc\neq\ddd. \end{matrix}\right.$$
Note that, as $f'$ has compact support, the set $\{\ddd\in\FF_\OO\vert f_\ddd'\neq 0\}$ is finite, and we have $f'=\sum_\ddd f_\ddd'$. For any $\ddd$ we have $\sigma_t (f_\ddd')=\N\ddd^{it} f_\ddd'$, so $$\norm{f'-\sigma_t(f')}\leqslant\sum_\ddd \norm{f_\ddd'-\sigma_t(f_\ddd')}\leqslant\sum_\ddd \abs{1-\N\ddd^{it}} \norm{f_\ddd'},$$
so it is now obvious that when $\abs{t}$ is small enough, this is smaller that $\varepsilon/3$.
\end{proof}

The resulting $C^*$-dynamical system $(C_{k,\infty},\,(\sigma_t))$ is the one that was announced in the introduction of this paper.

\section{Algebraic structure of $(C_{k,\infty},\,(\sigma_t))$}

\subsection{The rational subalgebra $\HH$}

In this subsection, we construct a sub-$*$-algebra $\HH$ which will play the rôle of the rational Hecke algebra in the Bost-Connes construction. Even though this algebra is called ``rational'', it is defined as an algebra over $\CC$. It is called rational because, as we will later see, it admits a presentation with rational coefficients.\\

For any $\aaa\in\Ig_\OO$, let \index{muaaa@$\mu_\aaa$ ($\aaa\in\Ig_\OO$)} $\mu_\aaa\in C_c(\GG)$ be defined by
$$\forall (\chi,\ccc)\in\GG,\:\:\: \mu_\aaa(\chi,\ccc)=1_{\ccc=\aaa}\,.$$
For any $\phi\in\Hay$, for any $\lambda\in\phi(\Ci)^\tor$, let us define a function \index{epsilambda@$e(\phi,\lambda)$ ($\phi\in\Hay$, $\lambda\in\phi(\Ci)^\tor$)} $e(\phi,\lambda)\in C_c(\GG)$ by
$$\forall (\chi,\ccc)\in\GG,\:\:\: e(\phi,\lambda)(\chi,\ccc)=1_{\ccc=1}\,1_{\chi\in X_\phi}\,\chi(\lambda).$$

\begin{defi}
Let \index{Hecke@$\HH$} $\HH$ denote the $*$-subalgebra of $C_c(\GG)$ generated by the $\mu_\aaa$, for all $\aaa\in\Ig_\OO$, and the 
$e(\phi,\lambda)$, for all $\lambda\in\phi(\Ci)^\tor$, for all $\phi\in\Hay$.\n
\end{defi}
We will later show (Proposition \ref{denseH}) that $\HH$ is dense in $C_{k,\infty}$. For now, we concentrate in checking 
several algebraic relations between the generators $\mu_\aaa$ and $e(\phi,\lambda)$ (see Proposition \ref{xr21}). We will 
later see (Proposition \ref{presentH}) that the relations of Proposition \ref{xr21} define a presentation of 
$\HH$.\n 

Recall that the inverse map in $\GG$ is given by
\begin{equation}
\label{invG}
(\chi,\ccc)^{-1}=(\chi^\ccc,\ccc^{-1}).
\end{equation}
The product law in $C_c(\GG)$, defined by equation (\ref{convol}), can be rewritten as:
\begin{equation}
\label{convolc2}
\forall f,g\in C_c(\GG),\:\:\forall(\chi,\ccc)\in\GG,\:\:\: (fg)(\chi,\ccc)
= \sum_{\ccc_2\in \FF_\chi} f(\chi^{\ccc_2},\,\ccc\ccc_2^{-1})\,g(\chi,\ccc_2).
\end{equation}
From equation (\ref{invG}), we check that for any $\aaa\in\Ig_\OO$, the adjoint $\mu_\aaa^*$ is given by
$$\forall (\chi,\ccc)\in\GG,\:\:\: \mu_\aaa^*(\chi,\ccc)=1_{\ccc=\aaa^{-1}}\,.$$
Using formula (\ref{convolc2}), we then check that, for all $f\in C_c(\GG)$, for all $(\chi,\ccc)\in\GG$, we have
\begin{eqnarray}
(\mu_\aaa f)(\chi,\ccc) & = & 1_{\ccc\aaa^{-1}\in\FF_\chi} \; f(\chi,\ccc\aaa^{-1}),\label{muaf}\\
(f \mu_\aaa)(\chi,\ccc) & = & f(\chi^\aaa,\ccc\aaa^{-1}),\label{fmua}\\
(\mu_\aaa^* f)(\chi,\ccc) & = & f(\chi,\ccc\aaa),\label{muaef}\\
(f \mu_\aaa^*)(\chi,\ccc) & = & 1_{\aaa^{-1}\in\FF_\chi} \; f(\chi^{\aaa^{-1}},\ccc\aaa).\nonumber
\end{eqnarray}
From that, we deduce that $C_c(\GG)$ is unital, with unit $\mu_1$ (where, as usual, $1$ denotes the principal ideal $(1)=\OO$)
$$\mu_1=1,$$
and we also deduce the formulas
\begin{eqnarray}
(\mu_\aaa f \mu_\bbb^*)(\chi,\ccc) & = & 1_{\ccc\aaa^{-1}\in\FF_\chi} \; 1_{\bbb^{-1}\in\FF_\chi}\; 
f(\chi^{\bbb^{-1}},\ccc\aaa^{-1}\bbb),\label{muafmube}\\(\mu_\bbb^* f \mu_\aaa)(\chi,\ccc) & = & 
f(\chi^\aaa,\ccc\aaa^{-1}\bbb),\nonumber\\(\mu_\aaa \mu_\bbb^*)(\chi,\ccc) & = & 
1_{\bbb^{-1}\in\FF_\chi}\; 1_{\ccc=\aaa\bbb^{-1}},\label{muamube}\\(\mu_\bbb^* \mu_\aaa)(\chi,\ccc) & = & 
1_{\ccc=\aaa\bbb^{-1}}.\label{mubemua}\end{eqnarray}
In particular, for $\bbb=\aaa$, equation (\ref{muamube}) gives
\begin{eqnarray}
(\mu_\aaa \mu_\aaa^*)(\chi,\ccc) & = & 1_{\aaa^{-1}\in\FF_\chi}\; 1_{\ccc=1}.\label{xr20}
\end{eqnarray}

\begin{prop}\label{xr21}
The functions $\mu_\aaa$, for $\aaa\in\Ig_\OO$, and $e(\phi,\lambda)$, for $\phi\in\Hay$ and $\lambda\in\phi(\Ci)^\tor$, satisfy the following relations:
\end{prop}
\noindent\begin{tabular}{lll}
(a$_1$) & $\mu_\aaa^*\mu_\aaa=\mu_1$ & $\forall \aaa\in \Ig_\OO$ \\
(a$_2$) & $\sum_\phi e(\phi,0) = \mu_1$ & where $\phi$ runs over $\Hay$\\
(b) & $\mu_\aaa\mu_\bbb=\mu_{\aaa\bbb}$ & $\forall \aaa,\bbb\in \Ig_\OO$ \\
(c) & $\mu_\aaa\mu_\bbb^*=\mu_\bbb^*\mu_\aaa$ & $\forall \aaa,\bbb\in \Ig_\OO\:$ relatively prime \\
(d$_1$) & $e(\phi,\lambda)^*=e(\phi,-\lambda)$ & $\forall \phi\in\Hay,\;\;\forall \lambda\in \phi(\Ci)^\tor$ \\
(d$_2$) & $e(\phi,\lambda_1)e(\phi,\lambda_2)=e(\phi,\lambda_1+\lambda_2)$ & $\forall \phi\in\Hay,\;\;\forall 
\lambda_1,\lambda_2\in \phi(\Ci)^\tor$\\(d$_3$) & $e(\phi^1,\lambda_1)e(\phi^2,\lambda_2)=0$ & $\forall 
\phi^1\neq\phi^2\in\Hay,\;\;\forall \lambda_i\in \phi^i(\Ci)^\tor$\\(e) & $e(\phi,\lambda)\mu_\aaa=\mu_\aaa 
e(\aaa*\phi,\phi_\aaa(\lambda))$ & $\forall \aaa\in \Ig_\OO,\;\;\forall \phi\in\Hay,\;\;\forall \lambda\in 
\phi(\Ci)^\tor$ \\(f) & $\mu_\aaa e(\phi,\lambda) \mu_\aaa^*=\frac{1}{\N \aaa}\sum_{(\aaa^{-1}*\phi)_\aaa(\mu)=\lambda} 
e(\aaa^{-1}*\phi,\mu)$ & $\forall \aaa\in \Ig_\OO,\;\; \forall \phi\in\Hay,\;\;\forall \lambda\in 
\phi(\Ci)^\tor$\end{tabular}\vskip 3 mm
\begin{proof} (a$_1$): Equation (\ref{mubemua}) applied with $\bbb=\aaa$ gives
$$(\mu_\aaa^* \mu_\aaa)(\chi,\ccc)  =  1_{c=1} = \mu_1 (\chi,\ccc).$$

\noindent (a$_2$): One checks directly that $\sum_\phi e(\phi,0) = \mu_1$.\n

\noindent (b): Equation (\ref{muaf}) applied with $f=\mu_\bbb$ gives
$$(\mu_\aaa \mu_\bbb)(\chi,\ccc)  = 1_{\ccc\aaa^{-1}\in\FF_\chi} \;  1_{\bbb=\ccc\aaa^{-1}}\; = 1_{\bbb\in\FF_\chi} 
\;  1_{\bbb=\ccc\aaa^{-1}}\;.$$As $\bbb$ is in $\Ig_\OO$, we always have $\bbb\in\FF_\chi$, so we find
$$(\mu_\aaa \mu_\bbb)(\chi,\ccc)  =1_{\bbb=\ccc\aaa^{-1}}\; = 1_{\aaa\bbb=\ccc}.$$
Thus, $\mu_\aaa \mu_\bbb=\mu_{\aaa\bbb}$.\n

\noindent (c): By equations (\ref{muamube}), (\ref{mubemua}), it is enough to show that for all $(\chi,\ccc)\in\GG$, 
we have $$1_{\bbb^{-1}\in\FF_\chi}\; 1_{\ccc=\aaa\bbb^{-1}}=1_{\ccc=\aaa\bbb^{-1}}.$$
If $\ccc\neq\aaa\bbb^{-1}$, then both sides are zero, so the equality holds. If $\ccc=\aaa\bbb^{-1}$, then we have 
$\aaa\bbb^{-1}\in\FF_\chi$. As $\aaa$ and $\bbb$ are relatively prime, Lemma \ref{xr16} then shows that 
$\bbb^{-1}\in\FF_\chi$, so the equality holds.\n

\noindent (d$_1$): For all 
$\chi\in X$, as $\chi$ is a character, we have
$$\forall\lambda\in\psi(\Ci)^\tor,\:\:\:\chi(-\lambda)=\overline{\chi(\lambda)}.$$
Relation (d$_1$) follows.\n

\noindent (d$_2$) and (d$_3$): from equation (\ref{convolc2}) and the formula 
$\chi(\lambda_1+\lambda_2)=\chi(\lambda_1)\chi(\lambda_2)$, one checks directly that for all $(\chi,\ccc)\in\GG$, letting $\psi$ be such that $\chi\in X_\psi$,
we have $$\left( e(\phi^1,\lambda_1) e(\phi^2,\lambda_2) \right) (\chi,\ccc)=1_{\ccc=1}\; 1_{\phi^1=\phi^2=\psi}\; 
\chi(\lambda_1+\lambda_2),$$which proves (d$_2$) and (d$_3$).\n

\noindent (e): By equation (\ref{muaf}) and the definition of $e(\phi,\lambda)$, we have , for any $\phi\in\Hay$, for 
any $\lambda\in\phi(\Ci)^\tor$, for any $(\chi,\ccc)\in\GG$, $$(\mu_\aaa 
e(\aaa*\phi,\phi_\aaa(\lambda)))(\chi,\ccc)  =  1_{\ccc\aaa^{-1}\in\FF_\chi} \; 1_{\ccc\aaa^{-1}=1} \; 
1_{\chi\in X_{\aaa*\phi}}\; \chi(\phi_\aaa(\lambda))= 1_{\ccc\aaa^{-1}=1} \; 1_{\chi\in X_{\aaa*\phi}}\; \chi^\aaa(\lambda) $$
so, by equation (\ref{fmua}), this is equal to $(e(\phi,\lambda)\mu_\aaa)(\chi,\ccc)$.\n

\noindent (f): For any 
$\phi\in\Hay$, for any $\lambda\in\phi(\Ci)^\tor$, for any $(\chi,\ccc)\in\GG$, we have\begin{eqnarray}
(\mu_\aaa e(\phi,\lambda) \mu_\aaa^*)(\chi,\ccc) & = & 1_{\ccc\aaa^{-1}\in\FF_\chi} \; 1_{\aaa^{-1}\in\FF_\chi}\; 
e(\phi,\lambda)(\chi^{\aaa^{-1}},\ccc)\;\;\;\text{by equation (\ref{muafmube})}\nonumber\\& = & 
1_{\ccc\aaa^{-1}\in\FF_\chi} \; 1_{\aaa^{-1}\in\FF_\chi}\; 1_{\ccc=1}\;1_{\chi\in X_{\aaa^{-1}*\phi}}\; \chi^{\aaa^{-1}} 
(\lambda)\nonumber\\& = & 1_{\aaa^{-1}\in\FF_\chi}\; 1_{\ccc=1}\;1_{\chi\in X_{\aaa^{-1}*\phi}}\; \chi^{\aaa^{-1}}(\lambda)\nonumber\\& 
= & 1_{\aaa^{-1}\in\FF_\chi}\; 1_{\ccc=1}\;1_{\chi\in X_{\aaa^{-1}*\phi}}\; (\N\aaa)^{-1} \sum_{(\aaa^{-1}*\phi)_\aaa(\mu)=\lambda} 
\chi(\mu)\;\;\;\text{by Lemma \ref{xr15}.}\nonumber \end{eqnarray}Let us first suppose that $\aaa^{-1}\in\FF_\chi$. We 
then have\begin{eqnarray}
(\mu_\aaa e(\phi,\lambda) \mu_\aaa^*)(\chi,\ccc) & = & 1_{\ccc=1}\; 1_{\chi\in X_{\aaa^{-1}*\phi}}\; (\N\aaa)^{-1} 
\sum_{(\aaa^{-1}*\phi)_\aaa(\mu)=\lambda} \chi(\mu)\nonumber\\& = & (\N\aaa)^{-1} 
\sum_{(\aaa^{-1}*\phi)_\aaa(\mu)=\lambda} e(\aaa^{-1}*\phi,\mu)(\chi,\ccc),\nonumber \end{eqnarray}
so we are done.\n
 
\noindent Let us now suppose that $\aaa^{-1}\not\in\FF_\chi$. We then have $1_{\aaa^{-1}\in\FF_\chi}=0$, and hence 
$(\mu_\aaa e(\phi,\lambda) \mu_\aaa^*)(\chi,\ccc)=0$. Thus, it is enough to show that 
$\sum_{(\aaa^{-1}*\phi)_\aaa(\mu)=\lambda} e(\aaa^{-1}*\phi,\mu)(\chi,\ccc)=0$. We 
have $$\sum_{(\aaa^{-1}*\phi)_\aaa(\mu)=\lambda} e(\aaa^{-1}*\phi,\mu)(\chi,\ccc)=1_{\ccc=1}\; 
1_{\chi\in X_{\aaa^{-1}*\phi}}\;\sum_{(\aaa^{-1}*\phi)_\aaa(\mu)=\lambda} \chi(\mu),$$so it is enough to show that if 
$\chi\in X_{\aaa^{-1}*\phi}$, then $\sum_{\psi_\aaa(\mu)=\lambda} \chi(\mu)=0$, where we have set $\psi=\aaa^{-1}*\phi$. Let $\mu_1\in\psi(\Ci)^\tor$ 
such that $\psi_\aaa(\mu_1)=\lambda$ (see Lemma \ref{idsurj}). We then have 
$$\sum_{\psi_\aaa(\mu)=\lambda} \chi(\mu)=\sum_{\psi_\aaa(\mu_0)=0} \chi(\mu_0+\mu_1),$$
so $$\sum_{\psi_\aaa(\mu)=\lambda} \chi(\mu)=\left(\sum_{\mu_0\in\psi[\aaa]} 
\chi(\mu_0)\right)\chi(\mu_1).$$
But, by Lemma \ref{xr15}, as $\aaa^{-1}\not\in\FF_\chi$, the restriction of $\chi$ to 
$\psi[\aaa]$ is a non-trivial character of $\psi[\aaa]$, so$$\sum_{\mu_0\in\psi[\aaa]} \chi(\mu_0) = 0,$$so 
$\sum_{\psi_\aaa(\mu)=\lambda} \chi(\mu)=0$, which completes the proof.\end{proof}

\subsection{Presentation of $\HH$}

The goal of this subsection is to show (Proposition \ref{presentH}) that the relations (a)-(f) of Proposition 
\ref{xr21} define a presentation of $\HH$ as a $*$-algebra.\n

The proof of the next lemma follows that of Proposition 18 in \cite{BosCon95}.

\begin{lemma}
\label{span}
Let $\widetilde{\HH}$ be a $*$-algebra with elements $\widetilde{\mu}_\aaa$, for $\aaa\in\Ig_\OO$, and 
$\widetilde{e}(\phi,\lambda)$, for $\phi\in\Hay$ and $\lambda\in\phi(\Ci)^\tor$, satisfying the relations (a)-(f) of 
Proposition \ref{xr21}. Let $S$ be the following subset of $\widetilde{\HH}$:
$$S\;=\;\left\{\widetilde{\mu}_\aaa 
\widetilde{e}(\phi,\lambda)\widetilde{\mu}_\bbb^*,\;\;\;\aaa,\bbb\in\Ig_\OO\;\text{relatively prime},\;\phi\in\Hay, \;\lambda\in\phi(\Ci)^\tor\right\}.$$
Then:
\begin{mylist}
\item The elements $\widetilde{\mu}_\aaa$, for $\aaa\in\Ig_\OO$, and 
$\widetilde{e}(\phi,\lambda)$, for $\phi\in\Hay$ and $\lambda\in\phi(\Ci)^\tor$, belong to the linear span of $S$. More specifically:
$$
\widetilde{\mu}_\aaa=\sum_{\phi\in\Hay}\widetilde{\mu}_\aaa \widetilde{e}(\phi,0) \widetilde{\mu}_1^*\;\;\;\text{and}\;\;\;
\widetilde{e}(\phi,\lambda)=\widetilde{\mu}_1 \widetilde{e}(\phi,\lambda) \widetilde{\mu}_1^*.
$$
\item Let $x_1,x_2\in S$. For $i=1,2$, write $x_i=\mu_{\aaa_i} e(\phi^i,\lambda_i)\mu_{\bbb_i}^*$. Let $\ddd=\aaa_2+\bbb_1$ be the gcd of $\aaa_2$ and $\bbb_1$. Let $\ccc$ be the gcd of $\ddd^{-1}\aaa_1\aaa_2$ and $\ddd^{-1}\bbb_1\bbb_2$. Set $\psi=\ccc^{-1}\ddd^{-1}\aaa_2*\phi^1$ and $\lambda'=\phi^1_{\ddd^{-1}\aaa_2}(\lambda_1)+\phi^2_{\ddd^{-1}\bbb_1}(\lambda_2)$. Then:
\noindent
\begin{eqnarray}
x_1x_2 & = & 1_{\aaa_2*\phi^1=\bbb_1*\phi^2}\:\widetilde{\mu}_{\ddd^{-1}\aaa_1\aaa_2}\widetilde{e}\left(\ddd^{-1}\aaa_2*\phi^1,\lambda'\right)\widetilde{\mu}_{\ddd^{-1}\bbb_1\bbb_2}^*\\
& = &
1_{\aaa_2*\phi^1=\bbb_1*\phi^2}\:\sum_{\psi_\ccc(\gamma)=\lambda'} \widetilde{\mu}_{\ccc^{-1}\ddd^{-1}\aaa_1\aaa_2}\widetilde{e}(\psi,\gamma)
\widetilde{\mu}_{\ccc^{-1}\ddd^{-1}\bbb_1\bbb_2}^*.\label{linearspan}
\end{eqnarray}
In particular, equation (\ref{linearspan}) shows that $x_1x_2$ belongs to the $\CC$-linear span of $S$.
\item If the elements $\widetilde{\mu}_\aaa$ and $\widetilde{e}(\phi,\lambda)$ generate $\widetilde\HH$ as a $*$-algebra, then the set $S$ generates $\widetilde\HH$ as a $\CC$-vector space.
\end{mylist}
\end{lemma}
\begin{proof}(1) easily follows from relations (a$_1$), (a$_2$) of Proposition \ref{xr21}.\\

\noindent (2): We have
$$x_1x_2=\mu_{\aaa_1} e(\phi^1,\lambda_1)\mu_{\bbb_1}^*\mu_{\aaa_2} e(\phi^2,\lambda_2)\mu_{\bbb_2}^*.$$
Using relations (a$_1$), (b) and (c) of Proposition \ref{xr21}, we 
find 
$$\widetilde{\mu}_{\bbb_1}^*\widetilde{\mu}_{\aaa_2}=\widetilde{\mu}_{\ddd^{-1}\bbb_1}^*\widetilde{\mu}_\ddd^*\widetilde
{\mu}_\ddd\widetilde{\mu}_{\ddd^{-1}\aaa_2}=\widetilde{\mu}_{\ddd^{-1}\bbb_1}^*\widetilde{\mu}_{\ddd^{-1}\aaa_2}.$$
Hence, we get
$$x_1x_2=\widetilde{\mu}_{\aaa_1} \widetilde{e}(\phi^1,\lambda_1)  
\widetilde{\mu}_{\ddd^{-1}\aaa_2}\widetilde{\mu}_{\ddd^{-1}\bbb_1}^*   
\widetilde{e}(\phi^2,\lambda_2)\widetilde{\mu}_{\bbb_2}^*.$$Using relations (e) and (d$_1$) of Proposition \ref{xr21}, we get $$x_1x_2=\widetilde{\mu}_{\aaa_1}\widetilde{\mu}_{\ddd^{-1}\aaa_2}\widetilde{e}\left(\ddd^{-1}\aaa_2*\phi^1,\phi^1_{\ddd^{-1}
\aaa_2}(\lambda_1)\right)\widetilde{e}\left(\ddd^{-1}\bbb_1*\phi^2, 
\phi^2_{\ddd^{-1}\bbb_1}(\lambda_2)\right)\widetilde{\mu}_{\ddd^{-1}\bbb_1}^*\widetilde{\mu}_{\bbb_2}^*.$$Hence, relation (b) of Proposition \ref{xr21} gives
$$x_1x_2=\widetilde{\mu}_{\ddd^{-1}\aaa_1\aaa_2}\widetilde{e}\left(\ddd^{-1}\aaa_2*\phi^1,\phi^1_{\ddd^{-1}\aaa_2}(\lambda_1)
\right)\widetilde{e}\left(\ddd^{-1}\bbb_1*\phi^2, 
\phi^2_{\ddd^{-1}\bbb_1}(\lambda_2)\right)\widetilde{\mu}_{\ddd^{-1}\bbb_1\bbb_2}^*.$$Thus, using relations (d$_2$) and 
(d$_3$) of Proposition \ref{xr21}, we get
$$x_1x_2=1_{\aaa_2*\phi^1=\bbb_1*\phi^2}\:\widetilde{\mu}_{\ddd^{-1}\aaa_1\aaa_2}\widetilde{e}\left(\ddd^{-1}\aaa_2*\phi^1,\phi^1_{\ddd^{-1}\aaa_2}(\lambda_1)+
\phi^2_{\ddd^{-1}\bbb_1}(\lambda_2)\right)\widetilde{\mu}_{\ddd^{-1}\bbb_1\bbb_2}^*.$$
By definition of $\psi$, $\lambda'$ and $\ccc$, and using relation (b) of Proposition \ref{xr21}, we obtain
$$x_1x_2=1_{\aaa_2*\phi^1=\bbb_1*\phi^2}\:\widetilde{\mu}_{\ccc^{-1}\ddd^{-1}\aaa_1\aaa_2}\left(\widetilde{\mu}_{\ccc}\widetilde{e}(\ccc*\psi,\lambda')\widetilde{\mu}_{\ccc}^*\right)\widetilde{\mu}_{\ccc^{-1}\ddd^{-1}\bbb_1\bbb_2}^*.$$
Relation (f) of Proposition \ref{xr21} then gives the result.\\

\noindent (3): The $\CC$-linear span of $S$ contains the generators $\widetilde\mu_\aaa$ and $\widetilde e(\phi,\lambda)$ by (1) and is stable under multiplication by (2). Moreover, it is obviously stable under the involution. Hence, it is equal to $\widetilde\HH$.
\end{proof}

\begin{lemma}\label{baseH}
The functions $\mu_\aaa e(\phi,\lambda) \mu_\bbb^*$, for  $\aaa,\bbb\in\Ig_\OO$ 
relatively prime, $\phi\in\Hay$, and $\lambda\in\phi(\Ci)^\tor$, form a basis of $\HH$ as a $\CC$-vector space.
\end{lemma}
\begin{proof}
By Lemma \ref{span} (3), they generate $\HH$ as a $\CC$-vector space. Thus, we only have to prove that they are linearly 
independent. Let us suppose that there exist $\alpha_1,\ldots,\alpha_n\in\CC$, 
$\aaa_0,\ldots,\aaa_n,\bbb_0,\ldots,\bbb_n\in\Ig_\OO$ with $\aaa_i$ relatively prime to $\bbb_i$ for each $i$, 
$\phi^1,\ldots,\phi^n\in\Hay$, and, for each $i$, $\lambda_i\in\phi(\Ci)^\tor$, such 
that$$\mu_{\aaa_0}e(\phi^0,\lambda_0)\mu_{\bbb_0}^*=\sum_{i=1}^n 
\alpha_i\mu_{\aaa_i}e(\phi^i,\lambda_i)\mu_{\bbb_i}^*.$$By equation (\ref{muafmube}) and the definition of 
$e(\phi^i,\lambda_i)$, we have $$\forall(\chi,\ccc)\in\GG,\;\;\;(\mu_{\aaa_i} e(\phi^i,\lambda_i) 
\mu_{\bbb_i}^*)(\chi,\ccc)=1_{\ccc{\aaa_i}^{-1}\in\FF_\chi}\; 1_{{\bbb_i}^{-1}\in\FF_\chi}\; 
1_{\ccc{\aaa_i}^{-1}{\bbb_i}=1}\; 1_{\chi\in X_{\bbb_i^{-1}*\phi^i}}\; \chi{\bbb_i}^{-1}(\lambda_i).$$Thus, the support of 
$\mu_{\aaa_i} e(\phi^i,\lambda_i) \mu_{\bbb_i}^*$ is included 
in$$\left\{g=(\chi,\ccc)\in\GG,\;\;\;\chi\in X_{\bbb_i^{-1}*\phi^i}\;\;\text{and}\;\;\ccc={\aaa_i}{\bbb_i}^{-1}\right\}.$$
Let $I$ denote the set of all $i\neq 0$ such that ${\aaa_i}{\bbb_i}^{-1}={\aaa_0}{\bbb_0}^{-1}$ and 
$\bbb_i^{-1}*\phi^i=\bbb_0^{-1}*\phi^0$. We thus have$$\mu_{\aaa_0}e(\phi^0,\lambda_0)\mu_{\bbb_0}^*=\sum_{i\in I} 
\alpha_i\mu_{\aaa_i}e(\phi^i,\lambda_i)\mu_{\bbb_i}^*.$$As $\aaa_i$ is relatively prime to $\bbb_i$, we see that for all 
$i\in I$, we have $\aaa_i=\aaa_0$ and $\bbb_i=\bbb_0$, so $\phi^i=\phi^0$. Hence, we 
get$$\mu_{\aaa_0}e(\phi^0,\lambda_0)\mu_{\bbb_0}^*=\sum_{i\in I} 
\alpha_i\mu_{\aaa_0}e(\phi^i,\lambda_i)\mu_{\bbb_0}^*.$$Hence, multiplying by $\mu_{\aaa_0}^*$ on the left and by 
$\mu_{\bbb_0}$ on the right, and using relation (a) of Proposition \ref{xr21}, we get$$e(\phi^0,\lambda_0)=\sum_{i\in 
I} \alpha_i e(\phi^0,\lambda_i).$$But the $e(\phi^0,\lambda)$, for $\lambda\in\phi^0(\Ci)^\tor$, are linearly 
independant (use, for instance, the isomorphism $C\left(X_{\phi^0}\times\{1\}\right) 
\simeq  C^*\left(\phi^0(\Ci)^\tor\right)$ as in Lemma \ref{ident}), so this is absurd.\end{proof}\begin{prop}\label{presentH}The relations (a)-(f) of Proposition \ref{xr21} define a presentation of $\HH$ as a $*$-algebra.\end{prop}\begin{proof}Let $\widetilde\HH$ be another $*$-algebra having elements 
$\widetilde{\mu}_\aaa$, for $\aaa\in\Ig_\OO$ and $\widetilde{e}(\phi,\lambda)$, for $\phi\in\Hay$, and 
$\lambda\in\phi(\Ci)^\tor$, satisfying the relations (a)-(f) of Proposition \ref{xr21}. We want to show that there 
exists a unique morphism $\sigma\;:\;\HH\rightarrow\widetilde\HH$ such that $\sigma\mu_\aaa=\widetilde{\mu}_\aaa$ and 
$\sigma e(\phi,\lambda)=\widetilde{e}(\phi,\lambda)$.\n 

The unicity is clear, by definition of $\HH$. Let us prove 
existence. By Lemma \ref{baseH}, we may define a $\CC$-linear map $\sigma\;:\;\HH\rightarrow\widetilde\HH$ by 
letting$$\sigma(\mu_\aaa e(\phi,\lambda) \mu_\bbb^*)=\widetilde{\mu}_\aaa 
\widetilde{e}(\phi,\lambda)\widetilde{\mu}_\bbb^*$$for all $\aaa,\bbb\in\Ig_\OO$ relatively prime, $\phi\in\Hay$, and 
$\lambda\in\phi(\Ci)^\tor$. Clearly, $\sigma(f^*)=\sigma(f)^*$. Moreover, Lemma \ref{span} shows that $\sigma\mu_\aaa=\widetilde{\mu}_\aaa$, that $\sigma e(\phi,\lambda)=\widetilde{e}(\phi,\lambda)$, and that$$\sigma(f_1 
f_2) = \sigma(f_1) \sigma (f_2),$$ which completes the proof.\end{proof}

\subsection{Presentation of $C_{k,\infty}$}
\label{subsec_pres_C}

The goal of this subsection is to 
show (Proposition \ref{presentC}) that the relations (a)-(f) of Proposition \ref{xr21} define a presentation of 
$C_{k,\infty}$ as a $C^*$-algebra.\n 

Let $\phi\in\Hay$. Let \index{C_psi@$C_\phi$ ($\phi\in\Hay$)} $C_\phi$ denote the subset of $C_c(\GG)$ of all functions whose support is a 
subset of $X_\phi\times\{1\}$.
\begin{lemma}
\label{prod}
Let $f_1,f_2\in C_\phi$. For all $g\in\GG$, we have
$$(f_1f_2)(g)=f_1(g)f_2(g).$$
\end{lemma}
\begin{proof}
Let $g=(\chi,\ccc)\in\GG$. By equation (\ref{convolc2}), we have
$$(f_1f_2)(\chi,\ccc)
= \sum_{\ccc_2\in \FF_\chi} f_1(\chi^{\ccc_2},\,\ccc\ccc_2^{-1})\,f_2(\chi,\ccc_2).$$
Thus, as $f_1,f_2\in C_\phi$, we can only have a nonzero term when $\ccc\ccc_2^{-1}=1$ and $\ccc_2=1$. If $\ccc\neq1$, 
then we get $(f_1f_2)(\chi,\ccc)=0$, as expected. If $\ccc=1$, then we get$$(f_1f_2)(\chi,1)= 
f_1(\chi,1)\,f_2(\chi,1),$$as expected.
\end{proof}

In particular, we see that for any $f_1,f_2\in C_\phi$, we have $f_1f_2\in C_\phi$. We also have $f_1^*\in C_\phi$. 
Thus, $C_\phi$ is a $*$-subalgebra of $C_c(\GG)$.\n

Let us define a norm $\norm{\cdot}_\phi$ on $C_\phi$ by:
$$\forall f\in C_\phi,\;\;\;\norm{f}_\phi=\sup_{g\in\GG}\abs{f(g)}.$$

\begin{lemma}
\label{ident}
$C_\phi$ is a $C^*$-algebra for the norm $\norm{\cdot}_\phi$. We have isomorphisms of $C^*$-algebras
$$
C_\phi \simeq C(X_\phi) \simeq  C^*(\phi(\Ci)^\tor).
$$
\end{lemma}
\begin{proof}
The identification $X_\phi\times\{1\}\simeq X_\phi$ gives a bijection
$C_\phi \simeq C(X_\phi)$. By Lemma \ref{prod}, this is a $*$-isomorphism. By definition of 
$\norm{\cdot}_\phi$, this is an isometry, so $\norm{\cdot}_\phi$ is a $C^*$-norm on $C_\phi$. The isomorphism 
$C(X_\phi) \simeq  C^*(\phi(\Ci)^\tor)$ is a classical result, see Davidson \cite{Dav91}, Proposition 
VII.1.1.\end{proof}

\begin{coro}
$C_\phi$ is a $C^*$-subalgebra of $C_{k,\infty}$.
\end{coro}
\begin{proof}
It is a classical result (see \cite{Dav91}, Theorem 1.5.5) that any injective $*$-morphism between two $C^*$-algebras is an 
isometry. Apply this to the inclusion map $\iota\;:\;C_\phi\rightarrow C_{k,\infty}$.
\end{proof}

\begin{lemma}
\label{densepsi}
The $e(\phi,\lambda)$, for $\lambda\in\phi(\Ci)^\tor$, generate a norm-dense $*$-subalgebra of 
$C_\phi$.\end{lemma}
\begin{proof}
By definition of the $e(\phi,\lambda)$, the isomorphism $C_\phi\simeq  C^*(\phi(\Ci)^\tor)$ given by Lemma \ref{ident}
identifies $e(\phi,\lambda)$ with $\lambda$. But, by definition of $C^*(\phi(\Ci)^\tor)$, the $\lambda$ generate a dense 
$*$-subalgebra of $C^*(\phi(\Ci)^\tor)$, so the result follows.\end{proof}

\begin{prop}
\label{denseH}
$\HH$ is dense in $C_{k,\infty}$, and any $*$-representation of $\HH$ extends uniquely to a representation of 
$C_{k,\infty}$.\end{prop}
\begin{proof}
Let us first prove density. As $C_c(\GG)$ is dense in $C_{k,\infty}$, it is enough to show that any $f\in C_c(\GG)$ can 
be approached by elements of $\HH$. Let $f\in C_c(\GG)$. As $f$ has compact support, there is a finite subset 
$\{\ccc_1,\ldots,\ccc_n\}\subset\FF_\OO$ such that for all $(\chi,\ccc)\in\GG$, if $\ccc\not\in\{\ccc_1,\ldots,\ccc_n\}$, 
then $f(\chi,\ccc)=0$. Let $f_i$ be defined 
by$$\forall(\chi,\ccc)\in\GG,\;\;\;f_i(\chi,\ccc)=1_{\ccc=\ccc_i}\;f(\chi,\ccc).$$We 
have$$f=f_1+\cdots+f_n.$$
It is thus enough to show that each of the $f_i$ can be approached by elements of $\HH$. Let $i\in\NN$ such that $1\leqslant i \leqslant n$. Write $\ccc_i=\aaa_i^{-1}\bbb_i$, with $\aaa_i,\bbb_i\in\Ig_\OO$ relatively prime. Let 
$f_i'=\mu_{\aaa_i} f_i \mu_{\bbb_i}^*$. We have $f_i=\mu_{\aaa_i}^* f_i' \mu_{\bbb_i}$, so it is enough to show that each of the $f_i'$ 
can be approached by elements of $\HH$. By equation (\ref{muafmube}), we have, for all $(\chi, \ccc)\in 
\GG$,$$f_i'(\chi,\ccc)  =  1_{\ccc\aaa_i^{-1}\in\FF_\chi} \; 1_{\bbb_i^{-1}\in\FF_\chi}\; 
f_i(\chi^{\bbb_i^{-1}},\ccc\ccc_i).$$Thus, the support of $f_i'$ is a subset of $X\times\{1\}$. For $\phi\in\Hay$, 
let $f_{i,\phi}'$ be defined 
by$$\forall(\chi,\ccc)\in\GG,\;\;\;f_{i,\phi}'(\chi,\ccc)=1_{\chi\in X_\phi}f_i'(\chi,\ccc).$$
We have
$$f_i'=\sum_{\phi\in\Hay} f_{i,\phi}',$$
so it is enough to show that each of the $f_{i,\phi}'$ can be approached by 
elements of $\HH$. We have $f_{i,\phi}'\in C_\phi$, so the result follows from Lemma \ref{densepsi}.\n 

Now, let us prove 
that any $*$-representation of $\HH$ extends uniquely to a representation of $C_{k,\infty}$. Unicity follows from the 
density of $\HH$ in $C_{k,\infty}$. Let us show existence. Let $\pi$ be a $*$-representation of $\HH$. By definition of $C_{k,\infty}$, it is enough to show that $\pi$ extends to a $*$-representation of $C_c(\GG)$. 
The construction we just made with the $f_i$, $f_i'$ and $f_{i,\phi}'$ shows that as a $*$-algebra, $C_c(\GG)$ is 
generated by the $C_\phi$, for $\phi\in\Hay$, and the $\mu_\aaa$, for $\aaa\in\Ig_\OO$. It is thus enough to show that 
the restriction of $\pi$ to the group algebra $\CC[\phi(\Ci)^\tor]$ extends to a representation of $C_\phi$. But this follows from Lemma 
\ref{ident}
\end{proof}

\begin{prop}
\label{presentC}
The relations (a)-(f) of Proposition \ref{xr21} define a presentation of $C_{k,\infty}$ as a $C^*$-algebra.\end{prop}\begin{proof}Let $\widetilde C$ be another $C^*$-algebra 
having elements $\widetilde{\mu}_\aaa$, for $\aaa\in\Ig_\OO$ and $\widetilde{e}(\phi,\lambda)$, for $\phi\in\Hay$, and 
$\lambda\in\phi(\Ci)^\tor$, satisfying the relations (a)-(f) of Proposition \ref{xr21}. We want to show that there 
exists a unique morphism $\sigma\;:\;C_{k,\infty}\rightarrow \widetilde C$ such that 
$\sigma\mu_\aaa=\widetilde{\mu}_\aaa$ and $\sigma e(\phi,\lambda)=\widetilde{e}(\phi,\lambda)$.\n 

Unicity follows from 
the density of $\HH$ in $C_{k,\infty}$, see Proposition \ref{denseH}. Let us prove existence.\n 

Let $\widetilde \HH$ 
denote the $*$-algebra generated by the $\widetilde{\mu}_\aaa$ and the $\widetilde{e}(\phi,\lambda)$. By the universal 
property of $\HH$ (Proposition \ref{presentH}), there exists a $*$-morphism $\sigma\;:\;\HH\rightarrow\widetilde \HH$ 
such that $\sigma\mu_\aaa=\widetilde{\mu}_\aaa$ and $\sigma e(\phi,\lambda)=\widetilde{e}(\phi,\lambda)$. Composing it 
with the inclusion $\widetilde\HH\rightarrow\widetilde C$ gives a $*$-representation of $\HH$. By Proposition 
\ref{denseH}, this representation extends to a $*$-morphism from $C_{k,\infty}$ into $\widetilde C$, so we 
are done.\end{proof}The flow $(\sigma_t)$ has a simple expression for this presentation: one checks directly 
that\begin{equation}\label{flowmu}\forall t\in\RR,\;\;\forall 
\aaa\in\Ig_\OO,\;\;\;\sigma_t(\mu_\aaa)=\N\aaa^{it}\mu_\aaa\end{equation}and\begin{equation}\label{flowe}
\forall t\in\RR,\;\;\forall \phi\in\Hay,\;\;\;\forall\lambda\in\phi(\Ci)^\tor,\;\;\;\sigma_t 
e(\phi,\lambda)=e(\phi,\lambda).\end{equation}

\subsection{Galois symmetry of $(C_{k,\infty},\,(\sigma_t))$}

Recall that an action of $\Gal(K/k)$ on $X$ has been defined by equation (\ref{galX}).\n

Let $\Gal(K/k)$ act by $*$-automorphisms on $C_c(\GG)$ by
\index{sigma_f@$\sigma f$ ($\sigma\in\Gal(K/k)$, $f\in C_{k,\infty}$)} $$\forall\sigma\in\Gal(K/k),\:\:\forall f\in C_c(\GG),\:\:\forall (\chi,\ccc)\in\GG,\:\:\:(\sigma f)(\chi,\ccc)=f(\sigma\chi,\ccc).$$

\begin{defi}
We will still note \index{sigma_f@$\sigma f$ ($\sigma\in\Gal(K/k)$, $f\in C_{k,\infty}$)} $(\sigma,f)\mapsto\sigma f$ the unique extension (given by Lemma \ref{xr19}) of this action to an 
action of $\Gal(K/k)$ on $C_{k,\infty}$.\end{defi}

One checks directly that the action of $\Gal(K/k)$ on the generators is given by
\begin{equation}
\label{sigmamu}
\forall\sigma\in\Gal(K/k),\;\;\forall \aaa\in \Ig_\OO,\;\;\;\sigma\mu_\aaa=\mu_\aaa
\end{equation}
and
\begin{equation}
\label{sigmaelambda}
\forall\sigma\in\Gal(K/k),\;\;\forall \phi\in\Hay,\;\;\;\forall\lambda\in\phi(\Ci)^\tor,\;\;\;\sigma 
(e(\phi,\lambda))=e(\sigma\phi,\sigma\lambda).\end{equation}

\begin{prop}
\label{galsym}
The group $\Gal(K/k)$, endowed with its profinite topology, is a topological symmetry group of $(C_{k,\infty},\,(\sigma_t))$. In other words, the action of $\Gal(K/k)$ on $C_{k,\infty}$ is faithful, continuous, and commutes with the flow $(\sigma_t)$, \ie
\begin{equation}
\label{commutact}
\forall\sigma\in\Gal(K/k),\:\:\forall t\in\RR,\:\:\forall f\in C_{k,\infty},\:\:\:\sigma(\sigma_t(f))=\sigma_t(\sigma
f).
\end{equation}
\end{prop}
\begin{proof}
By Lemma \ref{xr19}, it is enough to check equation (\ref{commutact}) for $f\in C_c(\GG)$, but it 
then follows from the definitions.\n

Let us check that the action of $\Gal(K/k)$ on $C_{k,\infty}$ is faithful. Let $\sigma\in\Gal(K/k)$ with $\sigma\neq 1$. Let 
$\phi\in\Hay$. If $\sigma\phi\neq\phi$, then it is clear that $\sigma$ acts non-trivially on $C_{k,\infty}$. If 
$\sigma\phi=\phi$, then, by definition of $H^+$, we have $\sigma\in\Gal(K/H^+)$. By definition of $K$, the action of 
$\Gal(K/H^+)$ on $\phi(\Ci)^\tor$ is faithful. Thus, there exists $\lambda\in\phi(\Ci)^\tor$ such that 
$\sigma\lambda\neq\lambda$, so $e(\phi,\sigma\lambda)\neq e(\phi,\lambda)$. Thus, by equation (\ref{sigmaelambda}), $\sigma e(\phi,\lambda)\neq e(\phi,\lambda)$, so the action of $\Gal(K/k)$ on $C_{k,\infty}$ is faithful.\\

Let us check that the action of $\Gal(K/k)$ on $C_{k,\infty}$ is continuous. Let $f\in C_{k,\infty}$ and $\varepsilon>0$.
By Proposition \ref{denseH}, the subalgebra $\HH$ is dense in $C_{k,\infty}$, so there exists $f_0\in \HH$ with $\norm{f-f_0} < \varepsilon /3$. Write $f_0$ in the basis provided by Lemma \ref{baseH}:
$$f_0=\sum_{i\in I} c_i \,\mu_{\aaa_i} \,e(\phi^i,\lambda_i) \,\mu^*_{\bbb_i},$$
where $I$ is a finite set and where, for all $i\in I$, we have $c_i\in\CC$, $\aaa_i,\bbb_i\in\Ig_\OO$ relatively prime, $\phi^i\in\Hay$, and $\lambda_i\in\phi^i(\Ci)^\tor$. Let $K_0$ be the extension of $k$ generated by the $\lambda_i$ and all their conjugates under $\Gal(K/k)$. Thus, $K_0/k$ is a finite Galois subextension of $K/k$. Let $V=\Gal(K/K_0)$. By definition of the profinite topology, $V$ is a neighborhood of 1 in $\Gal(K/k)$. For all $\sigma\in V$, we have $\sigma f_0 = f_0$. We have $\norm{\sigma f- f_0} = \norm{\sigma (f-f_0)}= \norm{f-f_0}< \varepsilon /3$, so we find $\norm{\sigma f-f}< 2\varepsilon/3$. Let $W$ denote the open ball of radius $\varepsilon/3$ centered at $f$. For all $f'\in W$, we have $\norm{\sigma f' - \sigma f}=\norm{\sigma(f'-f)}=\norm{f'-f}<\varepsilon/3$, whence $\norm{\sigma f'-f}<\varepsilon$, which completes the proof of the continuity.
\end{proof}

\subsection{The Galois-fixed subalgebra}\label{gisubsec}

In this subsection, we introduce two $C^*$-subalgebras of $C_{k,\infty}$, and it will turn out (Lemma \ref{galinvar}) that they are the same one.\\

The first one, denoted by \index{C_I_O@$C^*(\Ig_\OO)$} $C^*(\Ig_\OO)$, is the $C^*$-subalgebra of $C_{k,\infty}$ generated by the $\mu_\aaa$, for all $\aaa\in\Ig_\OO$.\\

The second one, denoted by \index{C_kinf_gal@$C_{k,\infty}^{\Gal(K/k)}$} $C_{k,\infty}^{\Gal(K/k)}$, is the subset of $C_{k,\infty}$ of all fixed points under the action of $\Gal(K/k)$. This is a $C^*$-subalgebra of $C_{k,\infty}$.\\

Let \index{PhiMaj_aaa@$\Phi(\aaa)$ ($\aaa\in\Ig_\OO$)} $$\Phi\::\:\Ig_\OO\rightarrow\NN$$
and \index{Moebius_aaa@$\M(\aaa)$ ($\aaa\in\Ig_\OO$)} $$\M\::\:\Ig_\OO\rightarrow\ZZ$$
denote the Euler totient and Möbius inversion functions respectively, \ie $\Phi$ and $\M$ are the multiplicative functions defined, for all primes $\pp$, for all $n\geqslant 0$, by
\index{PhiMaj_aaa@$\Phi(\aaa)$ ($\aaa\in\Ig_\OO$)} $$\Phi(\pp^n)=\N\pp^n - 1_{n\geqslant 1}\N\pp^{n-1}.$$
\label{arith}
and
\index{Moebius_aaa@$\M(\aaa)$ ($\aaa\in\Ig_\OO$)}$$\M(\pp^n)=1_{n=0}-1_{n=1}.$$
Note that we have, for all $\aaa\in\Ig_\OO$,
$$\Phi(\aaa)=\sum_{\bbb\mid\aaa}\M(\bbb^{-1}\aaa)\N\bbb.$$

\begin{lemma}
\label{genphia}
For all $\phi\in\Hay$, for all $\aaa\in\Ig_\OO$, the $\OO$-module $\phi[\aaa]$ has exactly $\Phi(\aaa)$ generators.
\end{lemma}
\begin{proof}
Let $\aaa=\prod_i \pp_i^{n_i}$ be the factorization of $\aaa$. As the $\pp_i$ are relatively prime, we have $\aaa=\bigcap_i \pp_i^{n_i}$, so 
by equations (\ref{xr6}), and (\ref{xr7}), we have 
$$\psi[\aaa]=\bigoplus_i \psi[\pp_i^{n_i}],$$
so it is enough to do the proof when $\aaa$ is a prime power, but it is then easy.
\end{proof}

The proof of the next lemma has been inspired by that of Proposition 21 (b) in \cite{BosCon95} and of Proposition 4.1 (3) in \cite{HarLei97}.

\begin{lemma}\label{galinvar}
The two subalgebras  $C^*(\Ig_\OO)$ and $C_{k,\infty}^{\Gal(K/k)}$ of $C_{k,\infty}$ are the same:
$$C^*(\Ig_\OO)=C_{k,\infty}^{\Gal(K/k)}.$$
\end{lemma}
\begin{defi}\label{xr22}
We let \index{C_1@$C_1$} $C_1$ denote this $C^*$-algebra:
$$C_1=C^*(\Ig_\OO)=C_{k,\infty}^{\Gal(K/k)}.$$
This notation will be justified in subsection \ref{uniqkms}, where $C_1$ will be viewed as a spectral subspace of $C_{k,\infty}$ for the action of $\Gal(K/k)$.
\end{defi}
\begin{proof}
One inclusion is clear: $C_{k,\infty}^{\Gal(K/k)}$ contains $C^*(\Ig_\OO)$. Let us check the other inclusion. The Galois group $\Gal(K/k)$ is endowed with its profinite topology, so it is a compact abelian group. Let \index{dsigma@$d\sigma$} $d\sigma$ be the 
normalized Haar measure on it. Let us consider the map $\E$ defined by:
\index{E_gras@$\E$}
\begin{eqnarray}
\E\;:\;C_{k,\infty} & \longrightarrow & C_{k,\infty}^{\Gal(K/k)} \nonumber \\
x & \longmapsto & \int_{\Gal(K/k)} \sigma(x) d\sigma.\label{defE}
\end{eqnarray}
By Proposition \ref{denseH}, $\HH$ is dense in $C_{k,\infty}$, so $\E(\HH)$ is dense in $C_{k,\infty}^{\Gal(K/k)}$. But, by 
Lemma \ref{span}, $\HH$ is the linear span of the $\mu_\aaa e(\phi,\lambda)\mu_\bbb^*$, for $\aaa,\bbb\in\Ig_\OO$, 
$\phi\in\Hay$, and $\lambda\in\phi(\Ci)^\tor$. Thus, $\E(\HH)$ is the linear span of the $\mu_\aaa 
\E(e(\phi,\lambda))\mu_\bbb^*$. Hence, it is enough to show that for all $\phi\in\Hay$, for all 
$\lambda\in\phi(\Ci)^\tor$, the element $\E(e(\phi,\lambda))$ belongs to $C^*(\Ig_\OO)$.\n 

So let $\phi\in\Hay$ and 
$\lambda\in\phi(\Ci)^\tor$.\n 

Let us first assume that $\lambda=0$. The group $\Gal(H^+/k)$ acts transitively on $\Hay$ 
(see Theorem \ref{hayes1}). By Galois theory, the restriction map $\Gal(K/k)\rightarrow\Gal(H^+/k)$ is 
surjective. Hence, $\Gal(K/k)$ acts transitively on $\Hay$. Thus, by relation (a$_2$) in Proposition \ref{xr21}, we 
get 
\begin{equation}
\label{Ezero}
\E(e(\phi,\lambda))=1/h(\sgn),
\end{equation}
where $h(\sgn)$ is the cardinal of $\Hay$. So the proof is complete.\n 

Let us now assume that $\lambda\neq 0$. Let
$$\aaa=\ann_\OO(\lambda)=\left\lbrace a\in\OO,\;\;\;\phi_a(\lambda)=0\right\rbrace.$$
We have 
$\lambda\in\phi[\aaa]$ and, for all $\bbb\neq\aaa$ such that $\bbb\mid\aaa$, $\lambda\not\in\phi[\bbb]$. So $\lambda$ is a generator 
of the $\OO$-module $\phi[\aaa]$. Let $K_{\aaa}$ denote the extension of $H^+$ generated by the elements of 
$\phi[\aaa]$. By \cite{Hay92}, Theorem 16.2, $\Gal(K_{\aaa}/k)$ acts transitively on the set $\XX_{\aaa}$ defined 
by
$$\XX_{\aaa}=\left\{ (\psi,\mu),\;\;\;\psi\in\Hay,\;\;\;\mu\;\;\text{is a generator of}\;\; \psi[\aaa] \right\}.$$
By Galois theory, the map $\Gal(K/k)\rightarrow\Gal(K_{\aaa}/k)$ is surjective, so $\Gal(K/k)$ also acts 
transitively on $\XX_{\aaa}$. Thus, $\E(e(\phi,\lambda))$ only depends on $\aaa$. We therefore 
note
$$\E(e(\aaa^{-1}))=\E(e(\phi,\lambda)).$$
Relation (f) of Proposition \ref{xr21}, gives, for all $\psi\in\Hay$, for all $\bbb\in\Ig_\OO$,
$$\mu_{\bbb} e(\psi,0) \mu_{\bbb}^*=\frac{1}{\N \bbb}\sum_{\mu\in(\bbb^{-1}*\psi)[\bbb]} e(\bbb^{-1}*\psi,\mu).$$
Thus, equation (a$_2$) of Proposition \ref{xr21} gives
$$\mu_{\bbb} \mu_{\bbb}^*=\frac{1}{\N \bbb}\sum_{\psi\in\Hay}\;\;\sum_{\mu\in(\bbb^{-1}*\psi)[\bbb]} e(\bbb^{-1}*\psi,\mu).$$
Applying $\E$ to this equality and using Lemma \ref{genphia}, we get
$$\N \bbb\, \mu_{\bbb} \mu_{\bbb}^*\,=\,h(\sgn)\sum_{\ccc\mid\bbb} \Phi(\ccc) \E(e(\ccc^{-1})),$$
where $h(\sgn)$ is the cardinal of $\Hay$. Doing a Möbius inversion, we then find
$$h(\sgn)\;\Phi(\bbb)\E(e(\bbb^{-1}))=\sum_{\ccc\mid\bbb} \M(\ccc^{-1}\bbb)\,\N \ccc\, \mu_{\ccc} \mu_{\ccc}^*.$$
Thus, for all $\bbb\in\Ig_\OO$, we get the following explicit expression of 
$\E(e(\bbb^{-1}))$ as an element of $C^*(\Ig_\OO)$:
\begin{equation}
\label{expl_E}\E(e(\bbb^{-1}))=\frac{\sum_{\ccc\mid\bbb} \M(\ccc^{-1}\bbb)\N \ccc\, \mu_{\ccc} \mu_{\ccc}^*}{h(\sgn)\;\Phi(\bbb)}=\frac{\sum_{\ccc\mid\bbb} \M(\ccc^{-1}\bbb)\N \ccc\, \mu_{\ccc} \mu_{\ccc}^*}{h(\sgn)\;\sum_{\ccc\mid\bbb} \M(\ccc^{-1}\bbb)\N \ccc}\cdot\qedhere
\end{equation}
\end{proof}

\begin{prop}\label{presentinvar} $C_1$ is isomorphic to the universal $C^*$-algebra generated 
by elements $\widetilde\mu_\aaa$, for $\aaa\in\Ig_\OO$, subject to the relations (a$_1$), (b) and (c) of Proposition 
\ref{xr21}.\end{prop}\begin{proof}This follows directly from Proposition \ref{presentC} and Lemma 
\ref{galinvar}.\end{proof}

\subsection{Admissible characters}
\label{admsubsec}

Some ideas in this subsection have been inspired by \cite{HarLei97}, \S 5. Our main goal here is to prove Proposition \ref{adminj}, which will be useful for the classification of extremal KMS$_\beta$ states at low temperature.

\begin{lemma}
Let $\chi\in X$. Let $\phi\in\Hay$ be such that $\chi\in X_\phi$. The following conditions are equivalent:
\begin{mylist}
\item For any maximal ideal $\pp\in\Ig_\OO$, the restriction of $\chi$ to $\phi[\pp]$ is non-trivial.
\item For any $\bbb\in\Ig_\OO$ different from $1$, the restriction of $\chi$ to $\phi[\bbb]$ is non-trivial
\item $\FF_\chi=\Ig_\OO$.
\end{mylist}
\end{lemma}
\begin{proof}
(2)$\Rightarrow$(1) is trivial. (1)$\Rightarrow$(2): as $\bbb\neq 1$, there exists a maximal ideal $\pp$ dividing 
$\bbb$. By Equation (\ref{xr5}), we then have $\phi[\pp]\subset\phi[\bbb]$, so the result follows. (2)$\Rightarrow$(3): Let $\ccc\in\FF_x$. Write $\ccc=\bbb^{-1}\aaa$ with $\aaa,\bbb\in\Ig_\OO$ relatively prime. By Lemma \ref{xr16}, we have $\bbb^{-1}\in\FF_x$. Thus, by Lemma \ref{xr15}, the restriction of $\chi$ to $\phi[\bbb]$ is trivial, so $\bbb=1$, so $\ccc\in\Ig_\OO$. (3)$\Rightarrow$(2): Let $\bbb\in\Ig_\OO$ with $\bbb\neq 1$. We have $\bbb^{-1}\not\in\FF_x$, so the result follows by Lemma \ref{xr15}.
\end{proof}

\begin{defi}$\text{}$
A character $\chi\in X$ is said to be \emph{admissible} if it satisfies the above equivalent conditions. Let \index{X_adm@$X^\adm$} $X^\adm$ denote the topological subspace of $X$ of admissible elements.\end{defi}


Recall that $A_f$ is the ring of finite adèles of $k$ with respect to $\OO$. Thus, $A_f$ is the restricted product of the $k_\pp$ 
with respect to the $\OO_\pp$, where $\pp$ runs over all finite places of $k$.\n

The following lemma is well-known.

\begin{lemma}\label{kpap}
Let $\aaa\in\Ig_\OO$. The diagonal map
 $\iota\,:\,k\hookrightarrow A_f$
induces an $\OO$-module isomorphism
 $$k/\aaa\xrightarrow{\sim}\bigoplus_{\pp}k_\pp/\aaa_\pp,$$
 where $\pp$ runs over all finite places of $k$, $k_\pp$ is the completion of $k$ at $\pp$, and $\aaa_\pp$ is the 
closure of $\aaa$ in $k_\pp$.\end{lemma}
\begin{proof}
Let $R=\prod_{\pp}\aaa_\pp\subset A_f$. This contains $\iota(\aaa)$. Hence $\iota$ induces a map
$$k/\aaa\to A_f/R.$$
This map is an $\OO$-module morphism. It is injective because $\iota^{-1}(R)=\aaa$. By the strong approximation theorem 
(Theorem \ref{strongapprox}), the range of $\iota$ is dense in $A_f$. But by definition of the restricted product, $R$ 
is an open subset of $A_f$. Hence $\iota$ induces a surjection modulo $R$. Thus, $\iota$ induces an isomorphism of 
$\OO$-modules $k/\aaa\simeq A_f/R$. But $A_f/R=\bigoplus_{\pp}k_\pp/\aaa_\pp$, so the result follows.\end{proof}
\begin{lemma}
\label{hl}
For any ideal $\aaa\in\Ig_\OO$, for any finite place $\pp$ of $k$, there exists a character $\chi$ of $k_\pp/\aaa_\pp$ 
whose restriction to $\pp^{-1}\aaa_\pp/\aaa_\pp$ is nontrivial.\end{lemma}
\begin{proof}
Let $\F_\pp$ denote the residue field of $\OO_\pp$. This is a finite extension of $\F_p$. The ring $\OO_\pp$ is 
principal (as is any local ring of any Dedekind ring), so its maximal ideal $\pp\OO_\pp$ is equal to $u\OO_\pp$ for some 
$u\in\OO_\pp$. Now $\aaa_\pp$ is also an ideal of $\OO_\pp$, so it is equal to $u^v\OO_\pp$ for some $v\geqslant 0$. 
Hence we have $\pp^{-1}\aaa_\pp/\aaa_\pp=u^{v-1}\OO_\pp/u^{v}\OO_\pp$. But we have $k_\pp=\F_\pp((u))$ and 
$\OO_\pp=\F_\pp[[u]]$, so we can define a character $\chi$ on $k_\pp/\aaa_\pp$ by letting$$\chi\left(\sum_{k\in\ZZ} a_k 
u^k\right) = \exp\left(2i\pi\Tr^{\F_\pp}_{\F_p}(a_{v-1})/p\right).$$The restriction of $\chi$ to 
$\pp^{-1}\aaa_\pp/\aaa_\pp$ is non-trivial since we have 
$\chi(u^{v-1})=\exp(2i\pi/p).$\end{proof}\begin{lemma}\label{xr23}For any ideal $\aaa\in\Ig_\OO$, there exists a 
character $\chi$ of $k/\aaa$ whose restriction to $\pp^{-1}\aaa_\pp/\aaa_\pp$, for any finite place $\pp$ of $k$, is 
nontrivial.\end{lemma}\begin{proof}
Use Lemma \ref{kpap} to identify $k/\aaa$ with $\bigoplus_\pp k_\pp/\aaa_\pp$. For all $\pp$, let $\chi_\pp$ be a 
character of $k_\pp/\aaa_\pp$ as given by the preceding lemma. Let $\chi=\prod_\pp\chi_\pp$. Then $\chi$ is a 
character of $k/\aaa$ that has the required property.\end{proof}
\begin{lemma}
For any $\phi\in\Hay$, there exists an admissible character $\chi\in X_\phi$. In 
particular, $X^\adm$ is non-empty.\end{lemma}
\begin{proof}
Let $L$ denote the lattice corresponding to $\phi$. Write $L=\xi\aaa$ with $\xi\in\Ci^*$ and $\aaa\in\Ig_\OO$. Let 
$\chi_0$ be a character of $k/\aaa$ as given by Lemma \ref{xr23}. Define a character $\chi$ of $\phi(\Ci)^\tor$ by: 
$$\chi(\lambda)=\chi_0(e_L^{-1}(\lambda)/\xi).$$Then $\chi$ is admissible.\end{proof}

\begin{lemma}
\label{admidinj}
For any $\chi\in X^\adm$, the map $\Ig_\OO\rightarrow X$, $\aaa\mapsto \chi^\aaa$, is injective.
\end{lemma}
\begin{proof}
By definition of admissibility and equation (\ref{xr14}), we have $\FF_{\chi^\aaa}=\aaa^{-1}\Ig_\OO$, so the result follows.
\end{proof}

\begin{lemma}
For any $\chi\in X^\adm$, for any $\sigma\in\Gal(K/k)$, we have $\sigma\chi\in X^\adm$.
\end{lemma}
\begin{proof}
The actions of $\Gal(K/k)$ and of $\Ig_\OO$ on $X$ commute with one another. Hence, $\FF_{\sigma\chi}=\FF_\chi=\Ig_\OO$. Hence, $\sigma\chi$ is admissible.
\end{proof}

\begin{prop}
\label{adminj}
For any $\chi\in X^\adm$, the map $\Gal(K/k)\rightarrow X^\adm$, $\sigma\mapsto \sigma \chi$, is injective.
\end{prop}
\begin{proof}
Let $\phi\in\Hay$ such that $\chi\in X_\phi$. Let $1\neq\sigma\in\Gal(K/k)$. Suppose that $\sigma \chi=\chi$. We have $\sigma\chi\in X_{\sigma^{-1}\phi}$, so $\sigma^{-1}\phi=\phi$. Thus, $\sigma\phi=\phi$, so by 
definition of $H^+$, we see that $\sigma\in\Gal(K/H^+)$. Also, $\sigma$ induces a map 
$$\sigma\;:\;\phi(\Ci)^\tor\rightarrow\phi(\Ci)^\tor.$$For any $\lambda\in\phi(\Ci)^\tor$, for any $a\in\OO$, we have 
$\phi_a(\sigma\lambda)=(\sigma\phi_a)(\sigma\lambda)=\sigma(\phi_a(\lambda)),$ so $\sigma$ is an $\OO$-module 
automorphism of $\phi(\Ci)^\tor$. Let $L$ denote the lattice corresponding to $\phi$. Write $L=\xi\aaa$ with $\xi\in\Ci^*$ 
and $\aaa\in\Ig_\OO$. Thus, we have $\OO$-module isomorphisms\begin{equation}\label{isoms}k/\aaa\xrightarrow{\xi} kL/L 
\xrightarrow{e_L} \phi(\Ci)^\tor,\end{equation}which we use to identify $k/\aaa$ with $\phi(\Ci)^\tor$, as 
$\OO$-modules. Thus, $\sigma$ is seen as an $\OO$-module automorphism of $k/\aaa$. Use Lemma \ref{kpap} to identify $k/\aaa$ with $\bigoplus_\pp 
k_\pp/\aaa_\pp$. For any finite place $\pp$ of $k$, writing $k_\pp$ as a field of Laurent series as in the proof of Lemma \ref{hl}, one sees that as $k_\pp/\aaa_\pp\simeq k_\pp/\OO_\pp$ as $\OO_\pp$-modules, hence as $\OO$-modules. Hence, $\End_\OO(k_\pp/\aaa_\pp)=\OO_\pp$ acting by multiplication. Hence,
$$\End_\OO(k/\aaa)=\prod_\pp\OO_\pp.$$
View $\sigma$ as an element of $\End_\OO(k/\aaa)$, and write $\sigma=\prod_\pp \sigma_\pp$ with $\sigma_\pp\in\OO_\pp$ for all $\pp$.\n

By definition of $K$, the action of $\Gal(K/H^+)$ on $\phi(\Ci)^\tor$ 
is faithful. Thus, as an $\OO$-module automorphism of $\phi(\Ci)^\tor$, we have $\sigma\neq 1$. Thus, there exists a 
$\pp$ such that $\sigma_\pp\neq 1$, so 
$\sigma_\pp-1\in\OO_\pp-\{0\}$. As $\chi$ is admissible, there exists $\lambda\in\phi[\pp]$ such that $\chi(\lambda)\neq 1$. 
View $\lambda$ as an element of 
$\pp^{-1}\aaa_\pp/\aaa_\pp$. Let $\widetilde\lambda\in\pp^{-1}\aaa_\pp\subset k_\pp$ be a 
representative of $\lambda$. Let $\widetilde\mu=(\sigma_\pp-1)^{-1}\widetilde\lambda\in k_\pp$. Let $\mu$ denote the class 
of $\widetilde\mu$ in $k_\pp/\aaa_\pp$. We have $(\sigma_\pp-1)\mu=\lambda$, so $(\sigma-1)\mu=\lambda$, so $$\chi((\sigma-1)\mu)\neq 1,$$ so
$$\chi(\sigma\mu)\neq\chi(\mu),$$
which is absurd, as $\sigma\chi=\chi$.
\end{proof}

\subsection{Irreducibility of regular representations at admissible characters}
\label{subsec_irred}

The goal of this subsection is to show that the regular representations of $\GG$ associated to admissible characters are irreducible. This will be used to classify extremal KMS states at low temperature.\n

Recall that for any $\chi\in X$, we defined the regular representation $\pi_\chi$ of $C_c(\GG)$ by equation (\ref{regul}). By definition of $C_{k,\infty}$, $\pi_\chi$ extends uniquely to a representation of $C_{k,\infty}$.\n

Recall that $X^\adm$ is the subset of $X$ of admissible elements.

\begin{lemma}
\label{irredrep}
For all $\chi\in X^\adm$, the regular representation $\pi_\chi$ of $C_{k,\infty}$  is irreducible.
\end{lemma}
\begin{proof}
Let $\chi\in X^\adm$. Let $\phi$ be such that $\chi\in X_\phi$. The representation $\pi_\chi$ is a map $C_{k,\infty}\rightarrow B\ell^2(\GG_\chi)$. Identify $\FF_\chi$ with $\GG_\chi$ through the map $\ccc\mapsto (\chi,\ccc)$. As $\chi$ is admissible, we have $\FF_\chi=\Ig_\OO$. Thus, $\GG_\chi$ is identified with $\Ig_\OO$. Let $A\in B\ell^2(\Ig_\OO)$ such that
$$\forall f\in C_{k,\infty},\;\;\;\pi_\chi(f)A=A\pi_\chi(f).$$
Let us show that $A$ is a scalar multiple of the identity. For that, let us first prove that $A$ is diagonal. Let $(\varepsilon_\ccc)_{\ccc\in\Ig_\OO}$ be the standard orthonormal basis of $\ell^2(\Ig_\OO)$: in other words, for all $\ccc,\aaa\in\Ig_\OO$, $\varepsilon_\ccc(\aaa)=1_{\aaa=\ccc}$. Let $(a_{\ccc,\ddd})$ be the matrix representing $A$ in this basis. Thus, we have
$$\forall \ccc\in\Ig_\OO,\;\;\;A\varepsilon_\ddd=\sum_\ccc a_{\ccc,\ddd} \varepsilon_\ccc.$$
Using equation (\ref{regul}), we check that
\begin{equation}
\label{regulmu}
\forall \aaa\in\Ig_\OO,\:\:\forall\bbb\in\Ig_\OO\:\:\:\pi_\chi(\mu_\aaa)\varepsilon_\bbb=\varepsilon_{\aaa\bbb}
\end{equation}
and
$$\forall \psi\in\Hay,\:\:\forall \lambda\in\psi(\Ci)^\tor,\:\:\forall\bbb\in\Ig_\OO,\:\:\:\pi_\chi(e(\psi,\lambda))\varepsilon_\bbb=1_{\psi=\bbb^{-1}*\phi}\;\chi^\bbb(\lambda)\varepsilon_\bbb.$$
Now let $\underline\lambda=(\lambda_\psi)_{\psi\in\Hay}$ be a family with $\forall\psi\in\Hay$, $\lambda_\psi\in\psi(\Ci)^\tor$. Let
$$e(\underline\lambda)=\sum_{\psi\in\Hay} e(\psi,\lambda_\psi).$$
We have:
$$\forall\bbb\in\Ig_\OO,\:\:\: \pi_\chi(e(\underline\lambda))\varepsilon_\bbb=\chi^\bbb(\lambda_{\bbb^{-1}*\phi})\varepsilon_\bbb.$$
Thus, for all $\bbb\in\Ig_\OO$, we get
\begin{eqnarray}
A\pi_\chi(e(\underline\lambda))\varepsilon_\bbb & = & \sum_{\aaa\in\Ig_\OO} a_{\aaa,\bbb} \chi^\bbb(\lambda_{\bbb^{-1}*\phi}) \varepsilon_\aaa,\nonumber\\
\pi_\chi(e(\underline\lambda))A\varepsilon_\bbb & = & \sum_{\aaa\in\Ig_\OO} a_{\aaa,\bbb} \chi^\aaa(\lambda_{\aaa^{-1}*\phi}) \varepsilon_\aaa.\nonumber
\end{eqnarray}
Thus, for all $\aaa,\bbb\in\Ig_\OO$ such that $a_{\aaa,\bbb}\neq 0$, for all $\underline \lambda$, we get
\begin{equation}
\label{xr24}
\chi^\bbb(\lambda_{\bbb^{-1}*\phi})=\chi^\aaa(\lambda_{\aaa^{-1}*\phi}).
\end{equation}
If $\bbb^{-1}*\phi\neq\aaa^{-1}*\phi$, as $\chi$ is admissible, we can obviously choose $\underline\lambda$ to make equation (\ref{xr24}) fail. Thus, we have $\bbb^{-1}*\phi=\aaa^{-1}*\phi$. By letting $\underline\lambda$ vary, we see that $\chi^\bbb$ and $\chi^\aaa$ are the same character of $(\aaa^{-1}*\phi)(\Ci)^\tor$. Thus,  $\chi^\aaa=\chi^\bbb$. Thus, as $\chi$ is admissible, by Lemma \ref{admidinj}, we find $\aaa=\bbb$. Thus, $(a_{\ccc,\ddd})$ is a diagonal matrix. Lastly, using the equality $A\pi_\chi(\mu_\aaa)=\pi_\chi(\mu_\aaa) A$ for all $\aaa\in\Ig_\OO$, one sees that the diagonal entries $(a_{\ccc,\ccc})$ are all equal, so that $(a_{\ccc,\ddd})$ is a scalar multiple of the identity matrix. Thus $\pi_\chi$ is an irreductible representation.
\end{proof}

\subsection{A lemma on the action of $\Gal(K/k)$ on $\HH$}

In this subsection, we prove an important lemma which we will use in subsections \ref{uniqkms} and \ref{typeiii}.

\begin{defi}
\label{xr25} Let $F$ be a set of finite places of $k$. An ideal $\ccc\in\Ig_\OO$ is said to be $F$-\emph{localized} if all its prime divisors belong to $F$.
\end{defi}

\begin{defi}Let $\ddd\in\Ig_\OO$. Let \index{F_d@$F_\ddd$ ($\ddd\in\Ig_\OO$)} $F_\ddd$ be the set of all places of $k$ dividing $\ddd$. We define \index{Hecke_d@$\HH[\ddd]$} $\HH[\ddd]$ to be the $*$-algebra generated by the $\mu_\aaa$, for all $F_\ddd$-localized ideals $\aaa\in\Ig_\OO$, and the $e(\phi,\lambda)$, for all $\phi\in\Hay$ and $\lambda\in\phi[\ddd]$.
\end{defi}

Note that, for any $\ddd\in\Ig_\OO$, $\Gal(K/K_\ddd)$ acts trivially on $\HH[\ddd]$. Thus, the action of $\Gal(K/k)$ on $\HH[\ddd]$ gives an action of the quotient group $\Gal(K_\ddd/k)=\Gal(K/k)/\Gal(K/K_\ddd)$ on $\HH[\ddd]$ (remember that the field $K_\ddd$ was defined in Definition \ref{Kc}).

\begin{lemma}
\label{xr26}
Let $\ddd\in\Ig_\OO$. Let $\pp$ be a maximal ideal of $\Ig_\OO$ not dividing $\ddd$. Let $\sigma_\pp=(\pp,K_\ddd/k)\in\Gal(K_\ddd/k)$ be the Artin automorphism of $K_\ddd$ associated to $\pp$. For all $x\in \HH[\ddd]$, we have
\begin{equation}
 \label{xr27}
x\mu_\pp=\mu_\pp\sigma_\pp(x).
\end{equation}
\end{lemma}
\begin{proof}
Let $A$ denote the subset of $\HH[\ddd]$ of all elements $x$ such that equation (\ref{xr27}) holds. Obviously, $A$ is a sub-$\CC$-algebra of $\HH[\ddd]$. But $\HH[\ddd]$ is generated as a $\CC$-algebra by the $\mu_\aaa$, the $\mu_\aaa^*$, and the $e(\phi,\lambda)$, for all $F_\ddd$-localized ideals $\aaa\in\Ig_\OO$,  all $\phi\in\Hay$, and all $\lambda\in\phi[\ddd]$. Indeed, by relation (d$_1$) of Proposition \ref{xr21}, we have $e(\phi,\lambda)^*=e(\phi,-\lambda)$. Hence, in order to prove that $A=\HH[\ddd]$, it is enough to check that $\mu_\aaa\in A$, $\mu_\aaa^*\in A$, and $e(\phi,\lambda)\in A$ for any $F_\ddd$-localized ideal $\aaa\in\Ig_\OO$, any $\phi\in\Hay$, and any $\lambda\in\phi[\ddd]$.\\

Let $\aaa\in\Ig_\OO$ be a $F_\ddd$-localized ideal. By relation (b) of Proposition \ref{xr21}, we have $\mu_\aaa\mu_\pp=\mu_\pp\mu_\aaa=\mu_\pp\sigma_\pp(\mu_\aaa)$, so $\mu_\aaa\in A$. As $\aaa$ is $F_\ddd$-localized and $\pp$ does not divide $\ddd$, relation (c) of Proposition \ref{xr21} gives $\mu_\aaa^*\mu_\pp=\mu_\pp\mu_\aaa^*=\mu_\pp\sigma_\pp(\mu_\aaa^*)$, so $\mu_\aaa^*\in A$.\\

Now let $\phi\in\Hay$ and $\lambda\in\phi[\ddd]$. We have
\begin{eqnarray*}
e(\phi,\lambda)\mu_\pp & = & \mu_\pp e(\pp*\phi,\phi_\pp(\lambda))\;\;\;\text{by relation (e) of Proposition \ref{xr21}}\\
 & = & \mu_\pp e(\sigma_\pp\phi,\phi_\pp(\lambda))\;\;\;\text{by Theorem \ref{hayes1}}\\
 & = & \mu_\pp e(\sigma_\pp\phi, \sigma_\pp(\lambda))\;\;\;\text{by Theorem \ref{hayes2}, as $\pp\nmid\ddd$}\\
 & = & \mu_\pp\,\sigma_\pp(e(\phi,\lambda)),
\end{eqnarray*}
so $e(\phi,\lambda)\in A$, which completes the proof.
\end{proof}

\section{KMS$_\beta$ equilibrium states of $(C_{k,\infty},\,(\sigma_t))$}

\subsection{The Galois-invariant KMS$_\beta$ state at any temperature}

The goal of this subsection is to construct (Proposition \ref{existkms}), for any $\beta\in\RR^*_+$, a Galois-invariant KMS$_\beta$ state $\varphi_\beta$ of $(C_{k,\infty},\,(\sigma_t))$. We will also show (Proposition \ref{xr29}) that $\varphi_\beta$ is the only Galois-invariant KMS$_\beta$ state of $(C_{k,\infty},(\sigma_t))$.\\

Proposition \ref{presentinvar} shows that $C_1$ is isomorphic to the infinite tensor product
$$C_1=\bigotimes_\pp \tau_\pp,$$
where $\pp$ runs over the finite places of $k$ and where, for each $\pp$, $\tau_\pp$ \index{tau_p@$\tau_\pp$ ($\pp$ a finite place)} is the (Toeplitz) $C^*$-algebra generated by $\mu_\pp$. Note that the $\tau_\pp$ are nuclear.\\

Let $\beta\in\RR^*_+$. For each $\pp$, define a state $\varphi_{\beta,\pp}$ on $\tau_\pp$ by
$$\forall n,m\geqslant 0,\;\;\;\varphi_{\beta,\pp}(\mu_\pp^n\mu_\pp^{*m})=1_{n=m}\N\pp^{-n\beta}.$$
Define a state $\varphi_\beta$ \index{phi_var_beta@$\varphi_\beta$ ($\beta\in\RR^*_+$)} on $C_1$ by
$$\varphi_\beta=\bigotimes_\pp\varphi_{\beta,\pp}.$$
Note that we have
\begin{equation}
\label{xr28}
\forall\aaa,\bbb\in\Ig_\OO,\;\;\;\varphi_{\beta}(\mu_\aaa\mu_\bbb^{*})=1_{\aaa=\bbb}\N\aaa^{-\beta}.
\end{equation}
Recall that the map $\E\::\:C_{k,\infty}\rightarrow C_1$ was defined in equation (\ref{defE}).

\begin{defi}
We extend $\varphi_\beta$ to a state on $C_{k,\infty}$ by letting
\index{phi_var_beta@$\varphi_\beta$ ($\beta\in\RR^*_+$)}
$$\forall f\in C_{k,\infty},\;\;\;\varphi_\beta (f)= \varphi_\beta \left( \E ( f) \right).$$
\end{defi}

\begin{prop}
\label{existkms}
For any $\beta\in\RR^*_+$, the state $\varphi_\beta$ on $C_{k,\infty}$ is a KMS$_\beta$ state of $\left(C_{k,\infty},(\sigma_t)\right)$. In particular, the state $\varphi_\beta$ on $C_1$ is a KMS$_\beta$ state of $\left(C_1,(\sigma_t)\right)$.
\end{prop}
\begin{proof}
For any $f_1,f_2\in C_{k,\infty}$, we look for a bounded holomorphic function $F_{\beta,f_1,f_2}$ on the strip $0 < \im z < \beta$ realizing the KMS$_\beta$ property for the state $\varphi_\beta$ and the pair $(f_1,f_2)$.\n

As $\HH$ is a dense $(\sigma_t)$-invariant $*$-subalgebra of $C_{k,\infty}$, by \cite{BraRob}, \S 5.3.1, it is enough to do that for $f_1,f_2\in\HH$.\\

In Lemma \ref{baseH}, we found a basis of $\HH$ as a $\CC$-vector space. Obviously, it is enough to check the KMS$_\beta$ condition in the case when $f_1$ and $f_2$ are elements of that basis.\n

Thus, write $f_1=\mu_{\aaa_1} e(\psi^1,\lambda_1) \mu_{\bbb_1}^*$ and $f_2=\mu_{\aaa_2} e(\psi^2,\lambda_2) \mu_{\bbb_2}^*$ with $\aaa_i,\bbb_i\in\Ig_\OO$ relatively prime, with $\psi^i\in\Hay$ and with $\lambda_i\in\psi^i(\Ci)^\tor$.\n

By Lemma \ref{span} (2),
 $$f_1f_2\:=\:1_{\aaa_2*\psi^1=\bbb_1*\psi^2}\: \mu_{\ddd^{-1}\aaa_1\aaa_2} e( \ddd^{-1}\aaa_2*\psi^1, \lambda')\mu_{\ddd^{-1}\bbb_1\bbb_2}^*,$$
where $\ddd$ is the gcd of $\aaa_2$ and $\bbb_1$, and where $\lambda'=\psi^1_{\ddd^{-1}\aaa_2}(\lambda_1)+\psi^2_{\ddd^{-1}\bbb_1}(\lambda_2).$  We thus have
$$\E(f_1f_2)=1_{\aaa_2*\psi^1=\bbb_1*\psi^2}\:\mu_{\ddd^{-1}\aaa_1\aaa_2} \E (e( \ddd^{-1}\aaa_2*\psi^1, \lambda')) \mu_{\ddd^{-1}\bbb_1\bbb_2}^*.$$
Let $\ccc=\ann_\OO(\lambda')$. Using equation (\ref{expl_E}), we deduce
$$\E(f_1f_2)=1_{\aaa_2*\psi^1=\bbb_1*\psi^2}\;\mu_{\ddd^{-1}\aaa_1\aaa_2}\;
\frac{\sum_{\fff\mid\ccc} \M(\fff^{-1}\ccc)\N\fff\, \mu_{\fff} \mu_{\fff}^*}{h(\sgn)\;\Phi(\ccc)}\;
\mu_{\ddd^{-1}\bbb_1\bbb_2}^*,$$
where $h(\sgn)$ is the cardinal of $\Hay$. Using the formula for $\varphi_\beta (\mu_\aaa\mu_\bbb^*)$ given in equation (\ref{xr28}), we then get:
$$\varphi_\beta(f_1f_2)=\frac{1_{\aaa_2*\psi^1=\bbb_1*\psi^2}}{h(\sgn)\;\Phi(\ccc)}\;\sum_{\fff\mid\ccc} \M(\fff^{-1}\ccc)\,\N\fff\;
1_{\ddd^{-1}\aaa_1\aaa_2\fff=\ddd^{-1}\bbb_1\bbb_2\fff}\;\N(\ddd^{-1}\aaa_1\aaa_2\fff)^{-\beta}.$$
Now, the condition $\ddd^{-1}\aaa_1\aaa_2\fff=\ddd^{-1}\bbb_1\bbb_2\fff$ is equivalent to $\aaa_1\aaa_2=\bbb_1\bbb_2$ and, as $\aaa_i$ is relatively prime to $\bbb_i$, this is equivalent to $\aaa_1=\bbb_2$ and $\aaa_2=\bbb_1$. We thus get:
$$\varphi_\beta(f_1f_2)=\frac{1_{\aaa_1=\bbb_2}\;1_{\aaa_2=\bbb_1}\;1_{\psi^1=\psi^2}}{h(\sgn)\;\Phi(\ccc)}\;\sum_{\fff\mid\ccc} \M(\fff^{-1}\ccc)\,\N\fff\;
\N(\ddd^{-1}\aaa_1\aaa_2\fff)^{-\beta}.$$
Now, if $\aaa_2=\bbb_1$, then $\ddd=\aaa_2=\bbb_1$, and so $\ccc=\ann_\OO(\lambda_1+\lambda_2)$. Summing this up, we have
\begin{equation}
\label{phiffprime}
\varphi_\beta(f_1f_2)=\frac{1_{\aaa_1=\bbb_2}\;1_{\aaa_2=\bbb_1}\;1_{\psi^1=\psi^2}}{h(\sgn)\;\Phi(\ccc)}\;\sum_{\fff\mid\ccc} \M(\fff^{-1}\ccc)\,\N\fff\;
\N(\aaa_1\fff)^{-\beta}\;\;\;\;\text{where}\;\;\;\;\ccc=\ann_\OO(\lambda_1+\lambda_2).
\end{equation}
Swapping $f_1$ with $f_2$ amounts to swapping $1$ with $2$ in the indices, so we get
$$\varphi_\beta(f_2f_1)=\frac{1_{\aaa_1=\bbb_2}\;1_{\aaa_2=\bbb_1}\;1_{\psi^1=\psi^2}}{h(\sgn)\;\Phi(\ccc)}\;\sum_{\fff\mid\ccc} \M(\fff^{-1}\ccc)\,\N\fff\;
\N(\aaa_2\fff)^{-\beta}\;\;\;\;\text{where}\;\;\;\;\ccc=\ann_\OO(\lambda_1+\lambda_2).$$
Thus, we find
$$\varphi_\beta(f_2f_1)=\left(\frac{\N\aaa_2}{\N\aaa_1}\right)^{-\beta} \varphi_\beta(f_1f_2).$$
We already know that both sides vanish unless $\aaa_1=\bbb_2$, so we get
$$\varphi_\beta(f_2f_1)=\left(\frac{\N\aaa_2}{\N\bbb_2}\right)^{-\beta} \varphi_\beta(f_1f_2).$$
Now, we have for all $t\in\RR$,
$$\sigma_t (f_2) = \sigma_t (\mu_{\aaa_2} e(\psi^2,\lambda_2) \mu_{\bbb_2}^*) = \N\aaa_2^{it}\mu_{\aaa_2}e(\psi^2,\lambda_2) \N\bbb_2^{-it}\mu_{\bbb_2}^*=\left(\frac{\N\aaa_2}{\N\bbb_2}\right)^{it} f_2.$$
Thus, letting
$$F_{\beta,f_1,f_2}(z)=\left(\frac{\N\aaa_2}{\N\bbb_2}\right)^{iz} \varphi_\beta(f_1f_2)$$
defines a bounded holomorphic function $F_{\beta,f_1,f_2}$ on the strip, realizing the KMS$_\beta$ property for the state $\varphi_\beta$ and the pair $(f_1,f_2)$.
\end{proof}

\begin{prop}
\label{xr29}
Let $\beta\in\RR^*_+$.
\begin{mylist}
\item The state $\varphi_\beta$ on $C_1$ is the only KMS$_\beta$ state of $(C_1,(\sigma_t))$.
\item The state $\varphi_\beta$ on $C_{k,\infty}$ is the only Galois-invariant KMS$_\beta$ state of $(C_{k,\infty},(\sigma_t))$.
\end{mylist}
\end{prop}
\begin{proof}
Clearly, the two statements are equivalent. Let us prove (1). Let $\varphi$ be a KMS$_\beta$ state of $(C_1,(\sigma_t))$. Let us show that
$$\varphi=\varphi_\beta.$$
Let $\aaa,\bbb\in\Ig_\OO$. We have
\begin{eqnarray}
\varphi(\mu_\aaa\mu_\bbb^*) & = & \varphi(\mu_\bbb^*\sigma_{i\beta}(\mu_\aaa)) \nonumber\\
 & = & \N\aaa^{-\beta}\varphi(\mu_\bbb^*\mu_\aaa). \label{xr30}
\end{eqnarray}
Let us first work in the case when $\aaa\neq\bbb$. Let us prove that $\varphi(\mu_\aaa\mu_\bbb^*)=0$. As $\varphi(\mu_\aaa\mu_\bbb^*)=\overline{\varphi(\mu_\bbb\mu_\aaa^*)}$, we may swap $\aaa$ and $\bbb$, and therefore we may assume without loss of generality that $\aaa\nmid\bbb$. Let $\ddd=\aaa+\bbb$ denote the gcd of $\aaa$ and $\bbb$. We have
\begin{equation}
\label{xr31}
\mu_\bbb^*\mu_\aaa=\mu_{\ddd^{-1}\bbb}^*\mu_\ddd^*\mu_\ddd\mu_{\ddd^{-1}\aaa}=
\mu_{\ddd^{-1}\bbb}^*\mu_{\ddd^{-1}\aaa}.
\end{equation}
As $\ddd^{-1}\aaa$ and $\ddd^{-1}\bbb$ are relatively prime, we have 
\begin{equation}
\label{xr32}
\mu_{\ddd^{-1}\bbb}^*\mu_{\ddd^{-1}\aaa}=
\mu_{\ddd^{-1}\aaa}\mu_{\ddd^{-1}\bbb}^*.
\end{equation}
 Thus, equation (\ref{xr30}) applied to $\ddd^{-1}\aaa$ and $\ddd^{-1}\bbb$ gives
\begin{equation}
\label{xr33}
\varphi(\mu_{\ddd^{-1}\aaa}\mu_{\ddd^{-1}\bbb}^*)=\N(\ddd^{-1}\aaa)^{-\beta}\varphi_\beta(\mu_{\ddd^{-1}\aaa}\mu_{\ddd^{-1}\bbb}^*).
\end{equation}
As $\aaa\nmid\bbb$, we have $\ddd^{-1}\aaa\neq 1$, so equation (\ref{xr33}) gives
$$\varphi(\mu_{\ddd^{-1}\aaa}\mu_{\ddd^{-1}\bbb}^*)=0.$$
Hence, equation (\ref{xr32}) gives $\varphi(\mu_{\ddd^{-1}\bbb}^*\mu_{\ddd^{-1}\aaa})=0$, so equation (\ref{xr31}) gives $\varphi(\mu_\bbb^*\mu_\aaa)=0$, so equation (\ref{xr30}) gives $\varphi(\mu_\aaa\mu_\bbb^*)=0$. Hence, we have proven that
$$\aaa\neq\bbb\;\Rightarrow\;\varphi(\mu_\aaa\mu_\bbb^*)=0=\varphi_\beta(\mu_\aaa\mu_\bbb^*).$$
In the case when $\aaa=\bbb$, equation (\ref{xr30}) gives
\begin{eqnarray*}
\varphi(\mu_\aaa\mu_\aaa^*)  & = & \N\aaa^{-\beta}\varphi(\mu_\aaa^*\mu_\aaa)\\
& = & \N\aaa^{-\beta}\\
& = & \varphi_\beta(\mu_\aaa\mu_\aaa^*).
\end{eqnarray*}
Thus, we have proven that
$$\forall\aaa,\bbb\in\Ig_\OO,\;\;\;\varphi(\mu_\aaa\mu_\bbb^*)=\varphi_\beta(\mu_\aaa\mu_\bbb^*).$$
As the linear span of the $\mu_\aaa\mu_\bbb^*$ is the $*$-algebra generated by the $\mu_\aaa$, it is dense in $C_1$ (by Definition \ref{xr22}), so we get $\varphi=\varphi_\beta$.
\end{proof}

\subsection{Action of $\Gal(K/k)$ on extremal KMS$_\beta$ states}

As usual, $\Gal(K/k)$ is endowed with its profinite topology. It acts on the set of states by $(\sigma,\varphi)\mapsto\varphi\circ\sigma$. Obviously, the KMS$_\beta$ condition and factoriality are preserved by this action. Hence, the sets $K_\beta$ and $\EE(K_\beta)$ are invariant under the action of $\Gal(K/k)$.\n

The proof of the next proposition comes from that of Theorem 25 in \cite{BosCon95}.

\begin{prop}
\label{transkms}
For any $\beta\in\RR^*_+$, the action of $\Gal(K/k)$ on $\EE(K_\beta)$ is transitive.
\end{prop}
\begin{proof}
The main ingredient is that the Galois-fixed subalgebra has a unique KMS$_\beta$ state (cf. Proposition \ref{xr29}). As in the proof of Lemma \ref{galinvar}, let $d\sigma$ be the 
normalized Haar measure on $\Gal(K/k)$, and let $\E$ denote the map defined in equation (\ref{defE}).\n

Let $\varphi_1$,$\varphi_2\in\EE(K_\beta)$. Then $\varphi_1\circ \E$ and $\varphi_2\circ \E$ are Galois-invariant elements of $K_\beta$. Thus, by Proposition \ref{xr29}, they are equal:
$$\varphi_1\circ \E = \varphi_2\circ \E.$$
But we have, for $i=1,2$, 
\begin{equation}
\label{decextr}
\varphi_i\circ \E = \int_{\Gal(K/k)} \varphi_i\circ\sigma\;d\sigma.
\end{equation}
Equation (\ref{decextr}) gives two decompositions of the same state as a barycenter of extremal KMS$_\beta$ states, but such a decomposition is unique (cf. \cite{BraRob}, II, Theorem 5.3.30), so the orbits of $\varphi_1$ and of $\varphi_2$ under $\Gal(K/k)$ are the same one. 
\end{proof}

Let $S$ denote the space of all states of $C_{k,\infty}$, endowed with the weak$^*$ topology. Recall that the weak$^*$ topology on $S$ is the one for which a basis of open neighborhoods of a state $\varphi_0$ is given by the
\begin{equation}
\label{xr34}
B(\varphi_0;x_1,\ldots,x_n;\varepsilon)\:=\:\left\{\varphi\in S\:;\;\;\;\forall i,\;\abs{\varphi(x_i)-\varphi_0(x_i)}<\varepsilon \right \}
\end{equation}
for all $n\geqslant 1$, $x_1,\ldots,x_n\in C_{k,\infty}$, and $\varepsilon>0$.

\begin{lemma}
\label{contkms}
The action of $\Gal(K/k)$ on $S$, given by $(\sigma,\varphi)\mapsto \varphi\circ\sigma$, is continuous.
\end{lemma}
\begin{proof}
Let $\varphi_0\in S$, $n\leqslant 1$, $x_1,\ldots,x_n\in C_{k,\infty}$, and $\varepsilon>0$. Let $U=B(\varphi_0;x_1,\ldots,x_n;\varepsilon)$, as defined in equation (\ref{xr34}). Let us find an open set $V\subset\Gal(K/k)$ and an open set $W\subset S$ such that
\begin{equation}
\label{xr35}
\forall \sigma\in V,\;\forall \varphi\in W,\;\;\;\varphi\circ\sigma\in U.
\end{equation}
Let us take $W=B(\varphi_0\in S;x_1,\ldots,x_n;\varepsilon/2)$. By Proposition \ref{galsym}, for any $i$, $1\leqslant i\leqslant n$, the map
\begin{eqnarray}
\Gal(K/k) & \rightarrow & C_{k,\infty} \nonumber \\
\sigma & \mapsto & \sigma(x_i) \nonumber
\end{eqnarray}
is continuous, so the finite intersection
$$V=\bigcap_{i=1}^n \left\{\sigma\in\Gal(K/k),\;\;\;\norm{\sigma(x_i)-x_i}<\varepsilon/2\right\}$$
is an open neighborhood of $1$ in $\Gal(K/k)$. Hence, all we have to do is to check equation (\ref{xr35}). Let $\sigma\in V$ and $\varphi\in W$. Let $i$, $1\leqslant i \leqslant n$. We have $\norm{\sigma(x_i)-x_i}<\varepsilon/2$, so, as $\varphi$ is a state, $\abs{\varphi(\sigma(x_i))-\varphi(x_i)}<\varepsilon/2$. On the other hand, as $\varphi\in W$, we have $\abs{\varphi(x_i)-\varphi_0(x_i)}<\varepsilon/2$. Thus, $\abs{\varphi(\sigma(x_i))-\varphi_0(x_i)}<\varepsilon$, so $\varphi\circ\sigma\in U$.
\end{proof}

\subsection{Extremal KMS$_\beta$ states at low temperature $1/\beta<1$ and special values}

Recall that $X^\adm$ is the subspace of $X$ of admissible elements, that $\EE(K_\beta)$ is endowed with the weak$^*$ topology, and that $\Gal(K/k)$ is endowed with its profinite topology. In this subsection, for any $\beta>1$, we will construct a homeomorphism $X^\adm\rightarrow\EE(K_\beta)$,  $\chi\mapsto\varphi_{\beta,\chi}$, commuting with the actions of $\Gal(K/k)$, and we will show that both $\EE(K_\beta)$ (for $\beta>1$) and $X^\adm$ are principal homogeneous spaces under $\Gal(K/k)$. Moreover, we will compute the values of $\varphi_{\beta,\chi}$ at certain elements of $\HH$ and relate them to special values of partial zeta functions of $k$.\n

For any $\chi\in X^\adm$, as at the beginning of the proof of Lemma \ref{irredrep}, let us make the identification $\GG_\chi=\FF_\chi=\Ig_\OO$, so that $\pi_\chi$ is seen as a representation in $\ell^2(\Ig_\OO)$. Let $(\varepsilon_\aaa)_{\aaa\in\Ig_\OO}$ be the standard orthonormal basis of $\ell^2(\Ig_\OO)$.

\begin{defi}
Let \index{H@$H$} $H$ be the unbounded operator on $\ell^2(\Ig_\OO)$ defined by
$$\forall\aaa\in\Ig_\OO,\:\:\:H\varepsilon_\aaa=(\log\N\aaa)\varepsilon_\aaa.$$
\end{defi}

\begin{lemma}\label{hamilflow}
For any $\chi\in X^\adm$, for all $t\in\RR$, for all $f\in C_{k,\infty}$, we have
$$\pi_\chi(\sigma_t(f))=e^{itH}\,\pi_\chi(f)\;e^{-itH}.$$
\end{lemma}
\begin{proof}
By Lemma \ref{xr19}, it is enough to do the proof in the case when $f\in C_c(\GG)$. It is then a straightforward computation.
\end{proof}

The function $\beta\mapsto\Tr(e^{-\beta H})$ is trivially computed:

\begin{lemma}\label{xr36}
For all $\beta>1$, we have $\Tr(e^{-\beta H})=\zeta_{k,\infty}(\beta)$.
\end{lemma}
\begin{proof}
$\Tr(e^{-\beta H})=\sum_{\aaa\in\Ig_\OO} e^{-\beta\log\N\aaa}=\sum_{\aaa\in\Ig_\OO} \N\aaa^{-\beta}=\zeta_{k,\infty}(\beta).$
\end{proof}

\begin{defi}
For any $\chi\in X^\adm$, for any $\beta>1$, we define a linear functional $\varphi_{\beta,\chi}$ on $C_{k,\infty}$ by:
\index{phi_var_betax@$\varphi_{\beta,\chi}$ ($\beta\in\RR^*_+$, $\chi\in X^\adm$)}
$$\varphi_{\beta,\chi}(f)=\zeta_{k,\infty}(\beta)^{-1}\Tr\left(\pi_\chi(f)\; e^{-\beta H}\right).$$
\end{defi}

Let  $(\varepsilon_\aaa)_{\aaa\in\Ig_\OO}$ denote the standard basis of $\ell^2(\Ig_\OO)$.

\begin{lemma}
\label{xr38}
For any $\chi\in X^\adm$, for any $\beta>1$, $\varphi_{\beta,\chi}$ is a KMS$_\beta$ state of the $C^*$-dynamical system $(C_{k,\infty},(\sigma_t))$.
\end{lemma}
\begin{proof}
By Lemma \ref{xr36}, we have $\varphi_{\beta,\chi}(1)=1$. We also have, for any $f\in C_{k,\infty}$,
$$\varphi_{\beta,\chi}(ff^*)=\zeta_{k,\infty}(\beta)^{-1}\Tr\left(\pi_\chi(f^*)\;e^{-\beta H}\,\pi_\chi(f)\right)\geqslant 0,$$
so $\varphi_{\beta,\chi}$ is a state on $C_{k,\infty}$. For any $f,f'\in C_{k,\infty}$, let us define a bounded continuous function $F_{\beta,\chi,f,f'}$ on the strip $\{z\in\CC\:\vert\:0\leqslant \im z\leqslant\beta\}$ by:
$$F_{\beta,\chi,f,f'}(z)=\zeta_{k,\infty}(\beta)^{-1}\Tr\left(e^{-\beta H}\, \pi_\chi(f)\; e^{izH}\, \pi_\chi(f') \; e^{-izH}\right).$$
One checks that the restriction of $F_{\beta,\chi,f,f'}$ to $\{z\in\CC\:\vert\:0< \im z<\beta\}$ is holomorphic. By Lemma \ref{hamilflow}, we have, for all $t\in\RR$,
$$F_{\beta,\chi,f,f'}(t)=\varphi_{\beta,\chi}(f\sigma_t(f'))\:\:\:\:\text{and}\:\:\:\:F_{\beta,\chi,f,f'}(t+i\beta)=\varphi_{\beta,\chi}(\sigma_t(f')f).$$
So $\varphi_{\beta,\chi}$ is a KMS$_\beta$ state of $(C_{k,\infty},(\sigma_t))$.
\end{proof}

\begin{lemma}\label{xr37}
For any $\chi\in X^\adm$, for any $\beta>1$, for any $\sigma\in\Gal(K/k)$, we have
$$\varphi_{\beta,\sigma \chi}=\varphi_{\beta,\chi}\circ\sigma.$$
\end{lemma}
\begin{proof}
By definition of $\varphi_{\beta,\sigma \chi}$, it is enough to check that $\pi_{\sigma \chi}(f)=\pi_{\chi}(\sigma f).$ By Proposition \ref{denseH}, it is enough to prove it when $f$ is one of the $e(\psi,\lambda)$ or one of the $\mu_\aaa$. The result then follows from equations (\ref{sigmamu}), (\ref{sigmaelambda}).
\end{proof}

\begin{lemma}
\label{gns}
For any $\chi\in X^\adm$, for any $\beta>1$, the GNS representation of $\varphi_{\beta,\chi}$ is ($\pi_{\beta,\chi}, \Omega_{\beta,\chi}$), where $\pi_{\beta,\chi}\::\:C_{k,\infty}\rightarrow B(\ell^2(\Ig_\OO)\otimes\ell^2(\Ig_\OO))$ is given by
$$\pi_{\beta,\chi}(f)(\xi\otimes\eta)=\pi_\chi(f)\xi\otimes\eta,$$
and the cyclic vector $\Omega_{\beta,\chi}\in\ell^2(\Ig_\OO)\otimes\ell^2(\Ig_\OO)$ is given by
$$\Omega_{\beta,\chi}=\zeta_{k,\infty}(\beta)^{-1/2}\sum_{\aaa\in\Ig_\OO}\N\aaa^{-\beta/2}\varepsilon_\aaa\otimes\varepsilon_\aaa.$$
\end{lemma}
\begin{proof}
We obviously have
$$\varphi_{\beta,\chi}(f)=\langle\pi_{\beta,\chi}(f)\Omega_{\beta,\chi},\Omega_{\beta,\chi}\rangle.$$
Hence, we only have to show that $\Omega_{\beta,\chi}$ is a cyclic vector for $\pi_{\beta,\chi}$. For any maximal ideal $\pp$ of $\OO$ and any $n\geqslant 0$, using equation(\ref{regulmu}), we find
$$\forall \aaa\in\Ig_\OO,\;\;\;\pi_\chi(\mu^*_{\pp^n})\varepsilon_\aaa=1_{\pp^n\mid\aaa}\;\varepsilon_{\pp^{-n}\aaa}$$
and hence
$$\forall \aaa\in\Ig_\OO,\;\;\;\pi_\chi(\mu_{\pp^n}\mu^*_{\pp^n})\varepsilon_\aaa=1_{\pp^n\mid\aaa}\;\varepsilon_\aaa.$$
Thus, if we let $\nu_{\pp^n}=\mu_{\pp^n}\mu^*_{\pp^n}-\mu_{\pp^{n+1}}\mu^*_{\pp^{n+1}},$ we get
$$\forall \aaa\in\Ig_\OO,\;\;\;\pi_\chi(\nu_{\pp^n})\varepsilon_\aaa=1_{\pp^n\mid\aaa\;\text{and}\;\pp^{n+1 }\nmid\aaa}\;\;\varepsilon_\aaa.$$
Now let $\bbb\in\Ig_\OO$. Let us show that $\varepsilon_\bbb\otimes\varepsilon_\bbb$ is in the closure of $\pi_{\beta,\chi}(C_{k,\infty})(\Omega_{\beta,\chi})$. Write $\bbb=\prod_\pp \pp ^{n_\pp}$ with $n_\pp\geqslant 0$. For $T>0$, let $P_T$ denote the set of all maximal ideals $\pp$ with $\N\pp<T$. The family $(P_T)_T$ is a growing family of finite sets whose union is the set of all maximal ideals of $\OO$. For all $T$, let $\nu_T\in C_{k,\infty}$ be defined by
$$\nu_T=\prod_{\pp\in P_T} \nu_{\pp^{n_\pp}}.$$
We have $$\pi_{\beta,\chi}(\nu_T)(\Omega_{\beta,\chi})=\zeta_{k,\infty}(\beta)^{-1/2}\sum_{\aaa\in Q_T}\N\aaa^{-\beta/2}\varepsilon_\aaa\otimes\varepsilon_\aaa,$$
where $Q_T$ is the set of all $\aaa\in\Ig_\OO$ such that for all $\pp\in P_T$, the $\pp$-adic valuations of $\aaa$ and $\bbb$ are equal. As the series $\sum_{\aaa} \N\aaa^{-\beta}$ is convergent, we see that
$$\pi_{\beta,\chi}(\nu_T)(\Omega_{\beta,\chi}) \xrightarrow{T\rightarrow+\infty} \zeta_{k,\infty}(\beta)^{-1/2}\N\bbb^{-\beta/2}\varepsilon_\bbb\otimes\varepsilon_\bbb.$$
Thus, we have shown that $\varepsilon_\bbb\otimes\varepsilon_\bbb$ is in the closure of $\pi_{\beta,\chi}(C_{k,\infty})(\Omega_{\beta,\chi})$. Applying the $\pi_{\beta,\chi}(\mu_\aaa)$ and the $\pi_{\beta,\chi}(\mu_\aaa^*)$ to that shows that for all $\bbb_1,\bbb_2\in\Ig_\OO$, the element $\varepsilon_{\bbb_1}\otimes\varepsilon_{\bbb_2}$ is in the closure of $\pi_{\beta,\chi}(C_{k,\infty})(\Omega_{\beta,\chi})$.
\end{proof}

\begin{prop}
\label{factorIi}
For any $\chi\in X^\adm$, for any $\beta>1$, the state $\varphi_{\beta,\chi}$ is factorial (hence extremal) of type $\I_{\infty}$.
\end{prop}
\begin{proof}
Let $A$ denote the weak closure of $\pi_{\beta,\chi}(C_{k,\infty})$ in $B(\ell^2(\Ig_\OO)\otimes\ell^2(\Ig_\OO))$. By Lemma \ref{irredrep}, the representation $\pi_\chi$ is irreducible. Thus, inside $B\ell^2(\Ig_\OO)$, we have
$\pi_\chi(C_{k,\infty})'=\CC.$ Using Takesaki \cite{Tak}, I, Chapter IV, Proposition 1.6 (i), we deduce that inside $B(\ell^2(\Ig_\OO)\otimes\ell^2(\Ig_\OO))$, we have $$\pi_{\beta,\chi}(C_{k,\infty})'=\CC\otimes B\ell^2(\Ig_\OO).$$
Thus, using \cite{Tak}, I, Chapter IV, Proposition 1.6 (ii), we deduce
$$A=\pi_{\beta,\chi}(C_{k,\infty})''=B\ell^2(\Ig_\OO)\otimes\CC.$$
In particular, we have $A\simeq B\ell^2(\Ig_\OO)$, so $A$ is a factor of type $\I_\infty$.
\end{proof}

\begin{lemma}
\label{injkms}
For any $\beta>1$, the map $X^\adm\rightarrow \EE(K_\beta)$, $\chi\mapsto \varphi_{\beta,\chi}$ is injective.
\end{lemma}
\begin{proof}
We reuse the notations of the proof of the previous lemma. Let us extend  $\varphi_{\beta,\chi}$ to a state $\widetilde\varphi_{\beta,\chi}$ on the von Neumann algebra $A=B\ell^2(\Ig_\OO)\otimes\CC$ by:
$$\forall a\in B\ell^2(\Ig_\OO),\;\;\;\widetilde\varphi_{\beta,\chi}(a\otimes1)=\langle a(\Omega_{\beta,\chi}),\Omega_{\beta,\chi}\rangle.$$
For any $\widehat \beta > 0$ , we have $e^{-\widehat \beta H}\in B\ell^2(\Ig_\OO)$. We have, for all $\psi\in\Hay$, for all $\lambda\in\psi(\Ci)^\tor$:
$$\zeta_{k,\infty}(\beta)\lim_{\widehat\beta\rightarrow+\infty} \widetilde\varphi_{\beta,\chi}(\pi_\chi(e(\psi,\lambda))e^{-\widehat \beta H} \otimes 1)=\langle\pi_\chi(e(\psi,\lambda))(\varepsilon_1),\varepsilon_1\rangle=1_{\chi\in X_\psi}\;\chi(\lambda).$$
Thus, $\chi$ is uniquely determined.
\end{proof}

We can now prove the main result classifying extremal KMS$_\beta$ states at low temperature. Recall that $\Gal(K/k)$ is endowed with its profinite topology, and $\EE(K_\beta)$ is endowed with the weak$^*$ topology.

\begin{theo}
\label{propkms}
For any $\beta>1$, the topological space $\EE(K_\beta)$ is principal homogeneous under $\Gal(K/k)$. 
\end{theo}
\begin{proof}
We must show that for any $\varphi\in\EE(K_\beta)$, the map $\Gal(K/k)\rightarrow\EE(K_\beta)$, $\sigma\mapsto\varphi\circ\sigma$ is a homeomorphism. We already know that it is surjective (Proposition \ref{transkms}) and continuous (Lemma \ref{contkms}). Thus, as $\Gal(K/k)$ is compact, it only remains to show that it is injective. Let $\varphi\in\EE(K_\beta)$ and $\sigma\in\Gal(K/k)$ such that $\varphi\circ\sigma=\varphi$. We have to show that $\sigma=1$. Let $\chi\in X^\adm$. By Proposition \ref{factorIi}, we have $\varphi_{\beta,\chi}\in \EE(K_\beta)$. By Proposition \ref{transkms}, there exists $\tau\in\Gal(K/k)$ such that $\varphi=\varphi_{\beta,\chi}\circ\tau$. By Lemma \ref{xr37}, we have $\varphi=\varphi_{\beta,\tau \chi}$ and $\varphi\circ\sigma=\varphi_{\beta,\sigma\tau \chi}$, so $\varphi_{\beta,\sigma\tau \chi}=\varphi_{\beta,\tau \chi}$. By Lemma \ref{injkms}, we deduce $\sigma\tau \chi=\tau \chi$. By Proposition \ref{adminj}, we find $\sigma\tau=\tau$, so $\sigma=1$.
\end{proof}

\begin{theo}
\label{xr39}
For any $\beta>1$, the map $X^\adm\rightarrow\EE(K_\beta)$, $\chi\mapsto\varphi_{\beta,\chi}$ is a homeomorphism.
\end{theo}
\begin{proof}
It is injective by Lemma \ref{injkms}. Let us check surjectivity: let $\varphi\in\EE(K_\beta)$. Let $\chi_0\in X^\adm$. By Proposition \ref{transkms} and Lemma \ref{xr37}, there exists $\sigma\in\Gal(K/k)$ such that $\varphi=\varphi_{\beta,\sigma\chi_0}$. Thus, the map $\chi\mapsto\varphi_{\beta,\chi}$ is bijective. One checks that it is continuous. By definition of an admissible character, $X^\adm$ is a closed subspace of $X$. Thus, $X^\adm$ is compact, so the considered map is a homeomorphism.
\end{proof}

\paragraph*{Relations between certain special values of KMS$_\beta$ states and of partial zeta functions.} Let us now compute the values of the states $\varphi_{\beta,\chi}$ on some of the generators $e(\phi,\lambda)$. Let \index{A_plus@$A_+$} $A_+$ denote the subset of $\Ig_\OO$ of all ideals $\aaa$ such that $\sigma_\aaa=1$, where $\sigma_\aaa=(\aaa,H^+/k)\in\Gal(H^+/k)$ is the Artin automorphism of $H^+$ associated to $\aaa$. For any $\ccc\in\Ig_\OO$ and any $\sigma\in\Gal(K_\ccc/k)$, let \index{A_c_sigma@$A_{\ccc,\sigma}$ ($\ccc\in\Ig_\OO$, $\sigma\in\Gal(K_\ccc/H^+)$)} $A_{\ccc,\sigma}$ denote the subset of $A_+$ of all ideals $\aaa$ prime to $\ccc$ and such that $\sigma_\aaa=\sigma$, where $\sigma_\aaa=(\aaa,K_\ccc/k)\in\Gal(K_\ccc/k)$ is the Artin automorphism of $K_\ccc$ associated to $\aaa$. Note that $A_+$ and the $A_{\ccc,\sigma}$ are generalized ideal classes of $\OO$.\\

Let \index{zeta_k_inf_A@$\zeta_{k,\infty}^+$} $\zeta_{k,\infty}^+$ and \index{zeta_k_inf_c_sigma@$\zeta_{k,\infty}^{\ccc,\sigma}$} $\zeta_{k,\infty}^{\ccc,\sigma}$ (for any $\ccc\in\Ig_\OO$ and $\sigma\in\Gal(K_\ccc/H^+)$) be the partial zeta functions associated to $A_+$ and $A_{\ccc,\sigma}$ respectively:
\begin{eqnarray*}
\zeta_{k,\infty}^+(\beta) & = & \sum_{\aaa\in A_+}\N\aaa^{-\beta}, \\
\zeta_{k,\infty}^{\ccc,\sigma}(\beta) & = & \sum_{\aaa\in A_{\ccc,\sigma}}\N\aaa^{-\beta}.
\end{eqnarray*}

\begin{theo}
\label{xr40}
Let $\beta>1$, $\phi\in\Hay$, and $\chi\in X^\adm\cap X_\phi$.
\begin{mylist}
\item We have
$$\varphi_{\beta,\chi}(e(\phi,0))=\frac{\zeta_{k,\infty}^+(\beta)}{\zeta_{k,\infty}(\beta)}\cdot$$
\item For any maximal ideal $\pp$ of $\OO$, for any $\lambda\in\phi[\pp]$, we have
$$\varphi_{\beta,\chi}(e(\phi,\lambda))=\zeta_{k,\infty}(\beta)^{-1}\left(\N\pp^{-\beta}\zeta_{k,\infty}^+(\beta)\;+\sum_{\sigma\in\Gal(K_\pp/H^+)}\chi(\sigma\lambda)\,\zeta_{k,\infty}^{\pp,\sigma}(\beta)\right).$$
\end{mylist}

\end{theo}
\begin{proof}
Let us first prove (1). By definition, $A_+$ is the subset of $\Ig_\OO$ of all ideals $\aaa$ such that $\sigma_\aaa=1$, where $\sigma_\aaa=(\aaa,H^+/k)\in\Gal(H^+/k)$ is the Artin automorphism of $H^+$ associated to $\aaa$. Hence, by Theorem \ref{hayes1}, 
we have
\begin{eqnarray*}
A_+ & =& \left\{\aaa\in\Ig_\OO,\;\;\;\aaa*\phi=\phi\right\}\\
&=& \left\{\aaa\in\Ig_\OO,\;\;\;\aaa^{-1}*\phi=\phi\right\}\\
&=& \left\{\aaa\in\Ig_\OO,\;\;\;\chi^\aaa\in X_\phi\right\}.
\end{eqnarray*}
Thus, by definition of $\varphi_{\beta,\chi}$, for any $\lambda\in\phi(\Ci)^\tor$, we have
\begin{equation}
\label{xr42}
\varphi_{\beta,\chi}(e(\phi,\lambda)) = \zeta_{k,\infty}(\beta)^{-1}\sum_{\aaa\in\Ig_\OO}1_{\chi^\aaa\in X_\phi}\,\chi^\aaa(\lambda)\N\aaa^{-\beta} = \zeta_{k,\infty}(\beta)^{-1}\sum_{\aaa\in A_+}\chi(\phi_\aaa(\lambda))\N\aaa^{-\beta}.
\end{equation}
Applying this equality to $\lambda=0$, we get (1).\\

Let us now prove (2). Let $\aaa\in A_+$. In the case when $\pp\mid\aaa$, we have $\phi_\aaa(\lambda)=0$, so
$$\chi(\phi_\aaa(\lambda))=1.$$
In the case when $\pp\nmid\aaa$, by Theorem \ref{hayes2}, we have $\phi_\aaa(\lambda)=\sigma_\aaa(\lambda)$. Hence, equation (\ref{xr42}) gives
\begin{eqnarray*}
\zeta_{k,\infty}(\beta)\varphi_{\beta,\chi}(e(\phi,\lambda))  & = & \sum_{\aaa\in A_+,\:\pp\mid\aaa}\chi(\phi_\aaa(\lambda))\N\aaa^{-\beta}+\sum_{\aaa\in A_+,\:\pp\nmid\aaa}\chi(\phi_\aaa(\lambda))\N\aaa^{-\beta}\\
& = & \sum_{\aaa\in A_+,\:\pp\mid\aaa}\N\aaa^{-\beta}+\sum_{\aaa\in A_+,\:\pp\nmid\aaa}\chi(\sigma_\aaa(\lambda))\N\aaa^{-\beta}\\
& = & \sum_{\aaa\in A_+}\N(\pp\aaa)^{-\beta}+\sum_{\sigma\in\Gal(K_\pp/H^+)}\:\sum_{\aaa\in A_{\pp,\sigma}}\chi(\sigma\lambda)\,\N\aaa^{-\beta}\\
& = & \N\pp^{-\beta}\zeta_{k,\infty}^+(\beta)+\sum_{\sigma\in\Gal(K_\pp/H^+)}\chi(\sigma\lambda)\,\zeta_{k,\infty}^{\pp,\sigma}(\beta),
\end{eqnarray*}
which proves (2).
\end{proof}

\subsection{Unicity of the KMS$_\beta$ state at high temperature $1/\beta\geqslant 1$}
\label{uniqkms}

Recall that in Proposition \ref{existkms}, for any $\beta\in\RR^*_+$, we found a Galois-invariant KMS$_\beta$ state $\varphi_\beta$ of $\left(C_{k,\infty},(\sigma_t)\right)$.\n

In this subsection, we will prove (Theorem \ref{xr61}) that when $\beta\leqslant 1$, there is no other KMS$_\beta$ state of $\left(C_{k,\infty},(\sigma_t)\right)$. In other words,
$$\beta\leqslant 1\;\;\Rightarrow\;\;K_\beta=\{\varphi_\beta\}.$$
Most of the ideas here come from \cite{BosCon95}, \S 7.\n

\noindent Let $\beta\in\RR^*_+$ be such that $\beta\leqslant 1$, and let $\psi$ be a KMS$_\beta$ state of $\left(C_{k,\infty},(\sigma_t)\right)$. We must show that $\psi=\varphi_\beta$.\\

\noindent Let $\widehat{\Gal(K/k)}$ be the dual group of $\Gal(K/k)$. As $\Gal(K/k)$ is profinite, $\widehat{\Gal(K/k)}$ is discrete.\\

\noindent Let $F$ be a nonempty finite set of finite places of $k$. Recall from Definition \ref{xr25} that an ideal $\aaa\in\Ig_\OO$ is $F$-localized if all its prime divisors belong to $F$. We will also need to define what it means to be \emph{$F$-localized} for an element of $\Gal(K/k)$, and for an element of $C(X)$.\\

Let us first define what it means to be $F$-localized for an element of $\Gal(K/k)$. We have
$$K=\lim_{\ccc\rightarrow} K_\ccc,$$
so
$$\Gal(K/k)=\lim_{\leftarrow\ccc} \Gal(K_\ccc/k),$$
so
$$\widehat{\Gal(K/k)}=\lim_{\ccc\rightarrow} \widehat{\Gal(K_\ccc/k)}.$$
This means that for any character $\nu$ of $\Gal(K/k)$, there exists $\ccc\in\Ig_\OO$ such that $\nu$ factors through the projection $\Gal(K/k)\rightarrow\Gal(K_\ccc/k)$. 

\begin{defi}
A character $\nu$ of $\Gal(K/k)$ is said to be $F$-\emph{localized} if there exists a $F$-localized ideal $\ccc\in\Ig_\OO$ such that $\nu$ factors through the projection $\Gal(K/k)\rightarrow\Gal(K_\ccc/k)$.
\end{defi}
\noindent Thus, any $\nu\in\Gal(K/k)$ is $F$-localized for some $F$.\\

Let \index{K_F@$K_F$} $K_F$ denote the extension of $H^+$ generated by the elements of the $\phi[F]$, for $\phi\in\Hay$. In other words,
$$K_F=\lim_{\ccc\rightarrow}K_\ccc,$$
where $\ccc$ runs over $\Ig_\OO$. Thus, a character $\nu$ of $\Gal(K/k)$ is $F$-localized if, and only if it factors through the surjection
$$\Gal(K/k)\rightarrow\Gal(K_F/k).$$

Let us now define what it means to be $F$-localized for an element of $C(X)$. For any $\phi\in\Hay$, let \index{phi_F@$\phi[F]$ ($\phi$ a Drinfel'd module)} $\phi[F]$ denote the following subgroup of $\phi(\Ci)^\tor$:
$$\phi[F]=\bigcup_{\ccc\;\text{is $F$-loc.}}\phi[\ccc],$$
where $\ccc$ runs over all the $F$-localized ideals in $\Ig_\OO$. Let $X_{\phi,F}$ denote the dual group of $\phi[F]$. The restriction-to-$\phi[F]$ map is a surjective morphism
$$X_{\phi}\rightarrow X_{\phi,F}.$$
Let \index{X_F@$X_F$}$X_F$ denote the (disjoint) union of the $X_{\phi,F}$, for all $\phi\in\Hay$. The restriction maps $X_\phi\rightarrow X_{\phi,F}$ give a surjection
$$X\rightarrow X_F.$$
This gives an injective morphism of $C^*$-algebras
$$C(X_F)\hookrightarrow C(X).$$
Thus, we see $C(X_F)$ as a $C^*$-subalgebra of $C(X)$.

\begin{defi}
An element $f\in C(X)$ is said to be $F$-\emph{localized} if it belongs to $C(X_F)$. In other words, $f$ is $F$-localized if, seen as a function $f\::\:X\rightarrow \CC$, it factors through the map $X\rightarrow X_F$.
\end{defi}

\begin{lemma}
\label{xr43}
The $C^*$-algebra $C(X_F)$ is generated by the $e(\phi,\lambda)$, for all $\phi\in\Hay$ and all $\lambda\in\phi[F]$.
\end{lemma}
\begin{proof}
That can be checked like Lemma \ref{densepsi}.
\end{proof}

For any character $\nu$ of $\Gal(K/k)$, let \index{C_chi@$C_\nu$ ($\nu$ a character of $\Gal(K/k)$)}\index{C_1@$C_1$}$C_\nu$ be the following spectral subspace of $C_{k,\infty}$:
$$C_\nu=\left\lbrace f\in C_{k,\infty},\;\;\;\forall \sigma\in \Gal(K/k),\;\;\;\sigma f=\nu(\sigma)f\right\rbrace.$$
Thus, when $\nu=1$ is the trivial character, the corresponding subspace $C_1$ is the Galois-fixed subalgebra computed in \ref{galinvar}.

\begin{lemma}
\label{xr44}
The following subspace is dense in $C_{k,\infty}$:
 $$\bigoplus_{\nu\in\widehat{\Gal(K/k)}} C_\nu.$$
\end{lemma}
\begin{proof}
As $\Gal(K/k)$ is a compact abelian group of $*$-automorphisms of $C_{k,\infty}$, this follows from a result found in Pedersen \cite{Ped79}, \S\S 8.1.4 and 8.1.10 page 349.
\end{proof}

\begin{lemma}
\label{xr45}
The states $\psi$ and $\varphi_\beta$ agree on $C_1$.
\end{lemma}
\begin{proof}
We saw in Proposition \ref{xr29} that $(C_1,\sigma_t)$ has only one KMS$_\beta$ state. Thus, as $\psi$ and $\varphi_\beta$ are KMS$_\beta$, they must agree on $C_1$.
\end{proof}

\begin{lemma}
\label{xr46}
Suppose that for any $\nu\in\widehat{\Gal(K/k)}$ with $\nu\neq 1$, the state $\psi$ vanishes on the spectral subspace $C_\nu$. Then $\psi=\varphi_\beta$.
\end{lemma}
\begin{proof}
By Lemma \ref{xr44}, in order to show that $\psi$ and $\varphi_\beta$ are equal, it is enough to show that they agree on $C_\nu$, for all $\nu$. We already know that $\psi$ and $\varphi_\beta$ agree on $C_1$. As $\varphi_\beta$ is $\Gal(K/k)$-invariant, it is easy to see that it vanishes on $C_\nu$ for any nontrivial $\nu$, so we deduce that $\psi=\varphi_\beta$.
\end{proof}

Thus, in order to prove that $\psi=\varphi_\beta$, it is enough to prove that $\psi$ vanishes on each of the spectral subspaces $C_\nu$, for $\nu\neq 1$. The following lemma, which is inspired by Lemma 27 (c) in \cite{BosCon95}, will be useful to prove that.

\begin{lemma}
\label{xr47}
Let $\nu\in\widehat{\Gal(K/k)}$ with $\nu\neq 1$.  Let $F$ be a nonempty finite set of finite places of $k$ such that $\nu$ is $F$-localized. Suppose that for any $F$-localized partial isometry $V\in C(X)\cap C_\nu$, we have
$$\forall x\in C_1,\;\;\;\psi(Vx)=0.$$
Then $\psi$ vanishes on the spectral subspace $C_\nu$.
\end{lemma}
\begin{proof}
From Theorems \ref{propkms} and \ref{xr39}, we know that $X^\adm$ is principal homogeneous under $\Gal(K/k)$. Thus, by choosing a base point $\chi_0\in X^\adm$, we can identify $\Gal(K/k)$ with $X^\adm$ through the map $\sigma\mapsto\sigma\chi_0$. Let $\fff\in\Ig_\OO$ be defined by
$$\fff=\prod_{\pp\in F}\pp.$$
For any $n\geqslant 1$, let $V_n\in C(X)$ be defined as follows. Let $\chi\in X$. If $\FF_\chi$ is of the form $\aaa^{-1}\Ig_\OO$ with $\aaa\mid\fff^n$, write $\chi=\sigma\chi_0^\aaa$ with $\sigma\in\Gal(K/k)$ and put $V_n(\chi)=\nu(\sigma)$. Otherwise, put $V_n(\chi)=0$. Note that $V_n$ is a partial isometry and belongs to $C_\nu$. Moreover, $V_n$ is $F$-localized because $\nu$ is.\\

 For any $\chi\in X$, we have
$\abs{V_n(\chi)}=1$ if $\FF_\chi$ is of the form $\aaa^{-1}\Ig_\OO$ for some $\aaa\mid\fff^n$, and $\abs{V_n(\chi)}=0$ otherwise. As $\abs{V_n}$ takes values in $\{0,1\}$, we have $\abs{V_n}=\abs{V_n}^2=V_nV_n^*$. Thus, by Lemmas \ref{xr16} and \ref{xr17}, for any $\chi\in X$, we have $$V_nV_n^*(\chi)\;=\;1_{\FF_\chi\,\subset\,\fff^{-n}\Ig_\OO}\;=\;\prod_{\pp\in F}1_{\fff^{-n}\pp^{-1}\not\in\FF_\chi}.$$
Hence, by equation (\ref{xr20}), we get
$$V_nV_n^*\;=\;\prod_{\pp\in F}\left(1-\mu_{\fff^{n}\pp}\mu_{\fff^{n}\pp}^*\right).$$
As $\psi$ and $\varphi_\beta$ agree on $C_1$ (Lemma \ref{xr45}) and $F$ is finite, we get
\begin{equation}
\label{xr48}
\psi\left(V_nV_n^*\right)\;=\;\varphi_\beta\left(\prod_{\pp\in F}\left(1-\mu_{\fff^{n}\pp}\mu_{\fff^{n}\pp}^*\right)\right)\xrightarrow{n\rightarrow\infty}1.
\end{equation}

Now let $x\in C_\nu$. We want to prove that $\psi(x)=0$.
For any $n\geqslant 1$, let $P_n=1-V_nV_n^*$. The Schwarz inequality gives
\begin{equation}
\label{xr49}
\abs{\psi(P_n x)}^2\leqslant\psi(P_n)\psi(xx^*).
\end{equation}
By equation (\ref{xr48}), we have $\psi(P_n)\xrightarrow{n\rightarrow\infty}0$, so equation (\ref{xr49}) gives $\psi(P_n x)\xrightarrow{n\rightarrow\infty}0$, so
\begin{equation}
\label{xr50}
\psi(V_nV_n^*x)\xrightarrow{n\rightarrow\infty}\psi(x).
\end{equation}
For any $n\geqslant 1$, as $x\in C_\nu$ and $V_n^*\in C_{\nu^{-1}}$, we have $V_n^*x\in C_1$. Hence, by assumption, we have
$$\psi(V_nV_n^*x)=0.$$
Together with equation (\ref{xr50}), this gives $\psi(x)=0$, which completes the proof of Lemma \ref{xr47}.
\end{proof}

Thus, in order to prove that $\psi=\varphi_\beta$, it is enough to prove the following lemma:

\begin{lemma}\label{xr51}
Let $\nu\in\widehat{\Gal(K/k)}$ with $\nu\neq 1$. Let $F$ be a nonempty finite set of finite places of $k$ such that $\nu$ is $F$-localized. For any $F$-localized partial isometry $V\in C(X)\cap C_\nu$, we have
$$\forall x\in C_1,\;\;\;\psi(Vx)=0.$$
\end{lemma}
\begin{proof}
This proof is directly inspired by the proof of Lemma 27 (b) of \cite{BosCon95}. It will make use of Lemmas \ref{Emup_lemma}, \ref{thetaext}, \ref{alphaout}, \ref{alphatwist}, \ref{xr58} and \ref{xr59}, and will only be completed on page \pageref{xr60}.\\

Let \index{V@$V$} $V\in C(X)$ be a $F$-localized partial isometry such that $V\in C_\nu$.\\

Let \index{E@$E$} $E=V^*V=VV^*$ (the algebra $C(X)$ is commutative). Note that $E$ is a projection and belongs to $C_1$. Let \index{C_1E@$C_{1,E}$ (reduced $C^*$-algebra)}$$C_{1,E}=E C_1 E=\left\lbrace f\in C_1,\;\;\;f=fE=Ef\right\rbrace$$ denote the reduced algebra. As $V$ is fixed by the flow $(\sigma_t)$ and $\psi$ and $\varphi_\beta$ are KMS$_\beta$ states for the flow $(\sigma_t)$, we see that $V$ belongs to the centralizer of $\psi$ and of $\varphi_\beta$.\\

 Let $\alpha$ denote the following automorphism of $C_{1,E}$:
\index{alpha(f)@$\alpha$}
$$\forall f\in C_{1,E},\;\;\;\alpha(f)=VfV^*.$$
Let \index{M@$M$} $M$ be the weak closure of $C_{k,\infty}$ in the GNS representation of $\varphi_\beta$. Let us extend the state $\varphi_\beta$ to a normal state \index{phi_var_beta_tilde@$\widetilde\varphi_\beta$ ($\beta\in\RR^*_+$)}$\widetilde\varphi_\beta$ on $M$. Let \index{M_1@$M_1$}$M_1\subset M$ denote the weak closure of $C_1$ in the GNS representation of $\varphi_\beta$.\\

As $V$ belongs to the centralizer of $\varphi_\beta$, for all $f\in C_{1,E}$, we have $\varphi_\beta(\alpha(f))=\varphi_\beta(f)$. Thus, $\alpha$ preserves $\varphi_\beta$, so it extends to an automorphism of the reduced algebra \index{M_1_E@$M_{1,E}$}$M_{1,E}$ preserving $\widetilde\varphi_\beta$.\\

Let $\ccc\in\Ig_\OO$ be a $F$-localized ideal such that $\nu$ factors through $\Gal(K_\ccc/k)$. 

\begin{lemma}
\label{Emup_lemma}
Let $\pp$ be a finite place of $k$ with $\pp\not\in F$. We have
\begin{equation}
\label{Emup}
E\mu_\pp\in C_{1,E}\;\;\;\text{and}\;\;\;\forall\pp\not\in F,\;\;\;\alpha(E\mu_\pp) =  \nu(\sigma_\pp)E\mu_\pp,
\end{equation}
where $\sigma_\pp=(\pp,K_F/k)\in\Gal(K_F/k)$ is the Artin automorphism of $K_F$ associated to $\pp$.
\end{lemma}
\begin{proof}
Let $\HH_F$ denote the $*$-algebra generated by the $e(\phi,\lambda)$, for all $\phi\in\Hay$ and all $\lambda\in\phi[F]$. By Proposition \ref{xr43}, we know that $\HH_F$ is norm-dense in $C(X_F)$. So let $(V_n)_{n\in\NN}$ be a sequence of elements of $\HH_F$ converging to $V$ in the norm topology.  Obviously, we have
$$\HH_F\subset\bigcup_{\ddd\in\Ig_\OO,\:\ddd\:F-\text{loc.}}\HH[\ddd],$$
where $\ddd$ runs over the $F$-localized elements of $\Ig_\OO$. Thus, for any $n\in\NN$, there exists a $F$-localized $\ddd_n\in\Ig_\OO$ such that $V_n\in\HH[\ddd_n]$. As $\pp\not\in F$ and $\ddd_n$ if $F$-localized, we have $\pp\nmid\ddd_n$, so Lemma \ref{xr26} gives $V_n \mu_\pp  =   \mu_\pp \sigma_\pp(V_n)$. Now view $\sigma_\pp$ as an automorphism of the $C^*$-algebra $C(X_F)$. In particular, it is continuous. Hence, we obtain
\begin{equation}
\label{Vmup}
\forall\pp\not\in F,\;\;\;V \mu_\pp = \mu_\pp \sigma_\pp(V)=\nu(\sigma_\pp) \mu_\pp V,
\end{equation}
and the result follows.
\end{proof}

\paragraph*{The ITPFI structure of $M_1$.} For any $\pp$, recall that \index{tau_p@$\tau_\pp$ ($\pp$ a finite place)}$\tau_\pp$ is the (Toeplitz) $C^*$-algebra generated by $\mu_\pp$, and that $\varphi_{\beta,\pp}$ is the restriction of $\varphi_\beta$ to $\tau_\pp$. Let $(\varepsilon_n)_{n\geqslant 0}$ denote the standard orthonormal basis of $\ell^2(\NN)$. Let $\pi_{\beta,\pp}$ be the following representation of $\tau_\pp$:
\begin{eqnarray}
\pi_{\beta,\pp}: \tau_\pp & \rightarrow & B(\ell^2(\NN)\otimes\ell^2(\NN)) \nonumber\\
\mu_\pp & \mapsto & \left(\varepsilon_n\otimes\varepsilon_m\mapsto \varepsilon_{n+1}\otimes\varepsilon_m.\right) \nonumber
\end{eqnarray}
Let $\Omega_{\beta,\pp}\in \ell^2(\NN)\otimes\ell^2(\NN)$ be the following vector:
$$\Omega_{\beta,\pp}=\sqrt{1-\N\pp^{-\beta}}\sum_{n\geqslant 0} \N\pp^{-n\beta/2}\varepsilon_n\otimes \varepsilon_n.$$
It is easy to check that the pair $\left(\pi_{\beta,\pp},\Omega_{\beta,\pp}\right)$ is the GNS representation of $\varphi_{\beta,\pp}$. Let \index{M_1p@$M_{1,\pp}$ ($\pp$ a finite place)} $M_{1,\pp}$ denote the weak closure of $\tau_\pp$ in the representation $\pi_{\beta,\pp}$. One checks that
\begin{equation}
\label{M1p}
M_{1,\pp}=B\ell^2(\NN)\otimes\CC.
\end{equation}
In particular $M_{1,\pp}$ is a type $\I_\infty$ factor. Let $\widetilde\varphi_{\beta,\pp}$ be the unique extension of $\varphi_{\beta,\pp}$ to a normal linear functional on $M_{1,\pp}$. Alternatively, $\widetilde\varphi_{\beta,\pp}$ is the restriction of $\widetilde\varphi_{\beta}$ to $M_{1,\pp}$. Note that the eigenvalue list of $\widetilde\varphi_{\beta,\pp}$ is the sequence \label{evlist} $$\left(\left(1-\N\pp^{-\beta}\right)\N\pp^{-n\beta}\right)_{n\geqslant 0}.$$
We have
\begin{equation}
\label{itpfi}
\left(M_1,\widetilde\varphi_\beta\right)\;=\;\bigotimes_{\pp} \left(M_{1,\pp},\widetilde\varphi_{\beta,\pp}\right),
\end{equation}
where $\pp$ runs over the finite places of $k$. Recall from \cite{Tak}, III, Chapter XIV, Corollary 1.10, that any ITPFI is a factor. In particular, \label{M1factor} $M_1$ is a factor. We will check later (Lemma \ref{xr65}) that it is of type $\III_{q^{-\beta}}$.\\

For any $\lambda\in\CC$ with $\abs{\lambda}=1$, let $\rho_{\pp,\lambda}$ denote the $*$-automorphism of $\tau_\pp$ such that \index{rho_pp_lambda@$\rho_{\pp,\lambda}$} $\rho_{\pp,\lambda}(\mu_\pp)=\lambda\mu_\pp$. As $\rho_{\pp,\lambda}$ preserves $\varphi_{\beta,\pp}$, it extends to an automorphism of $M_{1,\pp}$. Let \index{theta@$\theta_{F,\nu}$}$\theta_{F,\nu}$ be the following automorphism of $M_1$:
$$\theta_{F,\nu}=\left(\otimes_{\pp\in F}\; \id_{M_{1,\pp}}\right)\; \otimes\; \left(\otimes_{\pp\not\in F}\; \rho_{\pp,\nu(\sigma_\pp)}\right).$$

\begin{lemma}
\label{thetaext}
$\theta_{F,\nu}$ is an outer automorphism of $M_1$.
\end{lemma}
\begin{proof}
Suppose that $\theta_{F,\nu}$ is inner. Lemma 1.3.8 (b) from Connes \cite{Con73} states that there exists a sequence $(u_\pp)$, where, for any finite place $\pp$ of $k$, $u_\pp$ is an unitary of $M_{1,\pp}$ with
\begin{equation}
\label{xr52}
\forall x\in M_{1,\pp},\;\;\;\theta_{F,\nu}(x) = u_\pp x u_\pp^*
\end{equation}
and such that
\begin{equation}
\label{xr53}
\sum_{\pp}\left(1-\abs{\varphi_{\beta,\pp}(u_\pp)}\right) < \infty.
\end{equation}
As $M_{1,\pp}$ is a factor, equation (\ref{xr52}) determines $u_\pp$ up to multiplication by a $z\in\CC$ with $\abs{z}=1$. By definition of $\theta_{F,\nu}$, when $\pp\in F$, one can take $u_\pp=1$. When $\pp\not\in F$, one can take $u_\pp\in M_{1,\pp}=B\ell^2(\NN)$ to be the diagonal matrix with eigenvalue list $\left(\nu(\sigma_\pp^n)\right)_{n\in\NN}$. Using the expression of the GNS representation of $\varphi_{\beta,\pp}$ that we saw above, we get
\begin{eqnarray}
\varphi_{\beta,\pp}(u_\pp) & = & \left(1-\N\pp^{-\beta}\right) \sum_{n\in\NN} \nu(\sigma_\pp)^n \N\pp^{-n\beta} \nonumber\\
& = & \frac{1-\N\pp^{-\beta}}{1 - \nu(\sigma_\pp)\N\pp^{-\beta}}\cdot\nonumber
\end{eqnarray}
This is equal to $1$ whenever $\nu(\sigma_\pp)=1$, so, letting
$$Y=\left\{\pp\:;\;\;\;\pp\;\text{is a finite place of}\;k,\;\pp\not\in F,\;\text{and}\;\nu(\sigma_\pp)\neq 1\right\},$$
 Equation (\ref{xr53}) gives
\begin{equation}
\label{xr54}
\sum_{\pp\in Y}\left(1-\abs{\frac{1-\N\pp^{-\beta}}{1 - \nu(\sigma_\pp)\N\pp^{-\beta}}}\right) < \infty.
\end{equation}
Recall that we let $\ccc$ be a $F$-localized ideal such that $\nu$ factors through $\Gal(K_\ccc/k)$. This is a finite group, so the range of $\nu$ is finite, so there exists a $\gamma$ with $0<\gamma<1$ such that for any $\sigma\in\Gal(K/k)$ with $\nu(\sigma)\neq 1$, we have $\re \nu(\sigma) \leqslant \gamma$.
Let $\pp\in Y$. We have $\abs{1 - \nu(\sigma_\pp)\N\pp^{-\beta}}\geqslant 1 - \gamma \N\pp^{-\beta}$.
Thus, we find
$$1-\abs{\frac{1-\N\pp^{-\beta}}{1 - \nu(\sigma_\pp)\N\pp^{-\beta}}}\geqslant 1-\frac{1-\N\pp^{-\beta}}{1 - \gamma \N\pp^{-\beta}}=\frac{(1-\gamma)\N\pp^{-\beta}}{1 - \gamma \N\pp^{-\beta}}\geqslant (1-\gamma)\N\pp^{-\beta}.$$
As $1-\gamma>0$, together with Equation (\ref{xr54}), this gives
\begin{equation}
\label{xr55}
\sum_{\pp\in Y} \N\pp^{-\beta} < \infty.
\end{equation}
As $\beta \leqslant 1$, this implies that for all $s\geqslant 1$, we have
\begin{equation}
\label{xr56}
\sum_{\pp\in Y} \N\pp^{-s} \leqslant \sum_{\pp\in Y} \N\pp^{-1} < \infty.
\end{equation}
The \v{C}ebotarev density theorem (Theorem \ref{cebotarev}) states that for any $\sigma\in\Gal(K_\ccc/k)$, the following set $P_\sigma$ of places of $k$,
$$P_\sigma=\{\pp\:;\;\;\;\pp\;\text{does not ramify in}\;K_\ccc\;\text{and}\;\sigma_\pp=\sigma\},$$
has a positive Dirichlet density:
$$d(P_\sigma)>0.$$
Up to a finite set of places of $k$, we have $$Y=\bigcup_{\nu(\sigma)\neq 1} P_\sigma,$$
where $\sigma$ runs over $\Gal(K_\ccc/k)$. Hence, as $\nu\neq 1$, we have $d(Y)>0$, so $$\lim_{s\rightarrow 1_+}\sum_{\pp\in Y} \N\pp^{-s}=\infty,$$
contradicting Equation (\ref{xr56}).
\end{proof}

Define two subfactors $M_1^{F\pm}$ of $M_1$ as follows:
\index{M_1Fplusminus@$M_1^{F\pm}$} $$M_1^{F+}\;=\;\bigotimes_{\pp\in F} M_{1,\pp},$$
 $$M_1^{F-}\;=\;\bigotimes_{\pp\not\in F}
\left(M_{1,\pp},\widetilde\varphi_{\beta,\pp}\right),$$
where $\pp$ runs over the finite places of $k$. We thus have
$$M_1 = M_1^{F+} \otimes M_1^{F-}.$$
As the projection $E$ is $F$-localized, we have $E\in M_1^{F+}$, so letting \index{N@$N$} $$N=\left(M_1^{F+}\right)_E$$and using \cite{Tak}, I, Chapter IV, Proposition 1.9, we get
\begin{equation}
\label{xr57}
M_{1,E} = N \otimes M_1^{F-}.
\end{equation}

\begin{lemma}
\label{alphaout}
$\alpha$ is an outer automorphism of $M_{1,E}$.
\end{lemma}
\begin{proof}
Suppose that $\alpha$ is an inner automorphism of $M_{1,E}$. Let $\tau=\alpha^{-1}\circ\theta_{F,\nu}\in\Aut\left(M_{1,E}\right)$. By construction, $\tau$ induces the identity on $\CC E\otimes M_1^{F-}$. As $N\otimes \CC$ is of type $\I_\infty$, the restriction of $\tau$ to $N\otimes \CC$ is inner. By equation (\ref{xr57}), we get that $\tau$ is an inner automorphism of $M_{1,E}$. Hence, $\theta_{F,\nu}$ restricts to an inner automorphism of $M_{1,E}$. As $M_1$ is a factor, using Lemma 1.5.2 of \cite{Con73}, we deduce that $\theta_{F,\nu}$ is an inner automorphism of $M_1$, contradicting Lemma \ref{thetaext}.
\end{proof}

As we already noted, $V$ belongs to the centralizer of $\psi$. Define a linear functional \index{L@$L$ (linear functional)}$L$ on $C_{1,E}$ by
$$\forall x\in C_{1,E},\;\;\;L(x)=\psi(Vx)=\psi(xV).$$
We want to prove by contradiction that $L$ is zero. Thus, suppose that $L$ is nonzero.\\

The Schwarz inequality gives
$$\forall x\in C_{1,E},\;\;\;\abs{L(x)}^2\leqslant \psi(E) \psi(x^* x).$$
By Lemma \ref{xr45}, the states $\psi$ and $\varphi_\beta$ agree on $C_{1,E}$, so $\psi(x^*x)= \varphi_\beta(x^* x)$. Thus, $L$ extends to a normal linear functional on $M_{1,E}$, which we still note $L$.\\

As $\varphi_\beta$ is KMS$_\beta$ on $C_{k,\infty}$, by \cite{BraRob}, II, Corollary 5.3.4, there exists a unique extension of $(\sigma_t)$ to an ultraweakly continuous flow \index{sigma_t_tilde@$\widetilde\sigma_t$ ($t\in\RR$)}$(\widetilde\sigma_t)$ on $M$ for which $\widetilde\varphi_\beta$ is KMS$_\beta$.

\begin{lemma}
\label{alphatwist}
The linear functional $L$ satisfies the $\alpha$-twisted KMS$_\beta$ condition for the flow $(\widetilde\sigma_t)$ on $M_{1,E}$. In other words, for any $x,y\in M_{1,E}$, there exists a bounded continuous function $F_{x,y}$ on the strip $0\leqslant \im z \leqslant \beta$, holomorphic on the interior of the strip, such that for any $t\in\RR$, we have
\begin{equation}
\label{alphatweq}
F_{x,y}(t)=L(x\sigma_t(y))\;\;\;\text{and}\;\;\;F_{x,y}(t+i\beta)=L(\sigma_t(y)\alpha(x)).
\end{equation}
\end{lemma}
\begin{proof}
In the case when $x,y\in C_{1,E}$, this can be easily checked by applying the KMS$_\beta$ condition for $\psi$ to the pair $(Vx,y)$. As both $(\widetilde\sigma_t)$ and $L$ are ultraweakly continuous, the result follows.
\end{proof}

\begin{lemma}
\label{xr58}
There exists a nonzero $\widetilde\sigma_t$-invariant $w\in M_{1,E}$ such that 
$$\forall x\in M_{1,E},\;\;\; L(x)=\widetilde\varphi_\beta(wx).$$
\end{lemma}
\begin{proof}
Let $L = \abs{L}u$ be the polar decomposition of $L$ (see \cite{Tak}, I, Chapter III, Theorem 4.2). In particular, $u\in M_{1,E}$ is a partial isometry, $\abs{L}$ is a positive normal linear functional on $M_{1,E}$, and
$$\forall x\in M_{1,E},\;\;\;L(x) = \abs{L} (ux).$$
We want to apply Connes' Radon-Nikod\'ym theorem to $\abs{L}$ and $\widetilde\varphi_\beta$, seen as finite normal faithful weights on $M_{1,E}$.\\

As $L$ is $\widetilde\sigma_t$-invariant, by unicity of the polar decomposition, $u$ and $\abs{L}$ are $\widetilde\sigma_t$-invariant. As $\widetilde\varphi_\beta$ is KMS$_\beta$ for the flow $(\widetilde\sigma_t)$, we deduce that $\abs{L}$ is $\sigma_t^{\widetilde\varphi_\beta}$-invariant, where $(\sigma_t^{\widetilde\varphi_\beta})$ is the modular automorphism group associated to the finite faithful normal weight $\widetilde\varphi_\beta$ on $M_{1,E}$.\\

Connes' Radon-Nikod\'ym theorem (\cite{Con73}, Lemme 1.2.3 (b)) then states that there exists a positive $\widetilde\sigma_t$-invariant $h\in M_{1,E}$ such that 
$$\forall x\in M_{1,E},\;\;\; \abs{L}(x)=\widetilde\varphi_\beta(hx).$$
Letting $w=hu$, we get
$$\forall x\in M_{1,E},\;\;\; L(x)=\widetilde\varphi_\beta(wx),$$
and $w$ is nonzero by our assumption that $L$ is nonzero. It is $\widetilde\sigma_t$-invariant because both $u$ and $h$ are.
\end{proof}

\begin{lemma}
\label{xr59}
Let $w$ be given by Lemma \ref{xr58}. Then
$$\forall x\in M_{1,E},\;\;\; \alpha(x)w=wx.$$
\end{lemma}
\begin{proof}
Let $x,y\in M_{1,E}$. Let $F^L_{x,y}$ be the function given by Lemma \ref{alphatwist}, such that for any $t\in\RR$, we have
$$F^L_{x,y}(t)=L(x\widetilde\sigma_t(y))\;\;\;\text{and}\;\;\;F^L_{x,y}(t+i\beta)=L(\widetilde\sigma_t(y)\alpha(x)).$$
By definition of $w$, we get
$$F^L_{x,y}(t)=\widetilde\varphi_\beta(wx\widetilde\sigma_t(y))\;\;\;
\text{and}\;\;\;F^L_{x,y}(t+i\beta)=\widetilde\varphi_\beta(\widetilde\sigma_t(wy)\alpha(x)).$$
Now let $F^{\widetilde\varphi_\beta}_{\alpha(x),wy}$ be the function given by the KMS$_\beta$ property of $\widetilde\varphi_\beta$ applied to the pair $(\alpha(x),wy)$, so that
$$F^{\widetilde\varphi_\beta}_{\alpha(x),wy}(t)=\widetilde\varphi_\beta(\alpha(x)\widetilde\sigma_t(wy))\;\;\;
\text{and}\;\;\;F^{\widetilde\varphi_\beta}_{\alpha(x),wy}(t+i\beta)=\widetilde\varphi_\beta(\widetilde\sigma_t(wy)\alpha(x)).$$
Let $G=F^{\widetilde\varphi_\beta}_{\alpha(x),wy}-F^L_{x,y}$. Note that $G$ vanishes on $\RR+i\beta$. Therefore, one can extend $G$ to a holomorphic function on the broader strip $0<\im z<2\beta$ by letting
$$\forall z\in\CC\;\;\text{with}\;\;\beta<\im z<2\beta,\;\;\;G(z)=\overline{G(\bar z + 2i\beta)}.$$ As $G$ vanishes on $\RR+i\beta$ and is holomorphic on an open set containing $\RR+i\beta$, it vanishes everywhere, so $F^{\widetilde\varphi_\beta}_{\alpha(x),wy}=F^L_{x,y}$. In particular, evaluating that at $0$, we get
$$\widetilde\varphi_\beta(wxy)=\widetilde\varphi_\beta(\alpha(x)wy).$$
As this holds for all $y\in M_{1,E}$ and the state $\widetilde\varphi_\beta$ is faithful on $M_{1,E}$, we get
\[wx=\alpha(x)w.\qedhere \]
\end{proof}

We already know by Lemma \ref{alphaout} that $\alpha$ is outer. Thus, Proposition 4.1.16 of Sunder \cite{Sun87} shows that
$$\{y\in M_{1,E}\,;\;\;\;\forall x\in M_{1,E},\;\alpha(x)y=yx\}=\{0\}.$$
Together with Lemma \ref{xr59}, this shows that $w=0$. But by construction (cf. Lemma \ref{xr58}), $w$ is nonzero, so we get a contradiction. Thus, our assumption that $L$ is nonzero was false. Thus, $L$ is zero, so
$$\forall x\in C_{1,E},\;\;\;\psi(Vx)=0.$$

Now let $x\in C_1$. We have $ExE\in C_{1,E}$, so $\psi(VExE)=0$. As $E=V^*V=VV^*$ is a projection and belongs to the centralizer of $\psi$, we get $\psi(VExE)=\psi(EVEx)=\psi(Vx)$, so $\psi(Vx)=0$, which proves Lemma \ref{xr51}.\label{xr60}
\end{proof}

From that, we can deduce the main result of this subsection. Recall that we have assumed $0<\beta\leqslant 1$.

\begin{theo}
\label{xr61}
The $C^*$-system $(C_{k,\infty},\sigma_t)$ has exactly one KMS$_\beta$ state, $\varphi_\beta$.
\end{theo}
\begin{proof}
That follows from Lemmas \ref{xr46}, \ref{xr47}, and \ref{xr51}.
\end{proof}

\begin{coro}
\label{xr62}
The state $\varphi_\beta$ of $C_{k,\infty}$ is a factor state, \ie the von Neumann algebra $M$ is a factor.
\end{coro}
\begin{proof}
That follows from Theorem \ref{xr61} and \cite{BraRob}, II, Theorem 5.3.30 (3).
\end{proof}

\subsection{The type $\III_{q^{-\beta}}$ of the KMS$_\beta$ state at high temperature $1/\beta\geqslant 1$}
\label{typeiii}

{\bf Erratum:} This article wrongly claims that at high temperature $1/\beta\geqslant 1$, the unique KMS$_\beta$ state is of type $\III_{q^{-\beta}}$, where $q$ is the cardinal of the constant subfield of $k$. It has been shown by Neshveyev and Rustad \cite{NesRus12} that the correct type is $\III_{q^{-\beta d_\infty}}$ where $d_\infty$ is the degree of the place $\infty$. {\bf Lemma \ref{xr71} and the results that depend on it, especially Theorem \ref{xr80}, are wrong.} More specific errata have been inserted below.\\

Let us go on with the notations of the preceding subsection. In particular, we assume $\beta\leqslant 1$. The goal of this subsection is to prove (Theorem \ref{xr80}) that the state $\varphi_\beta$ on $C_{k,\infty}$ is of type $\III_{q^{-\beta}}$. In other words, we want to show that the factor $M$ is of type $\III_{q^{-\beta}}$. Before doing that, we will show that the subfactor $M_1$ is of type $\III_{q^{-\beta}}$.\\

Recall from equation (\ref{itpfi}) that $M_1$ is the following infinite tensor product, where $\pp$ runs over the finite places of $k$:
$$\left(M_1,\widetilde\varphi_\beta\right)\;=\;\bigotimes_{\pp} \left(M_{1,\pp},\widetilde\varphi_{\beta,\pp}\right).$$
Here, each of the $M_{1,\pp}$ is a type $\I_\infty$ factor, so the usual methods (cf. Araki and Woods \cite{AraWoo68}, \cite{Con94}) allowing to compute asymptotic ratio sets cannot be applied directly to $M_1$. Instead, we will first find an integer $\tau\in\NN$ and, for each $\pp$, a projection $e_\pp\in M_{1,\pp}$ such that the reduced factor $M_{1,\pp,e_\pp}$ is of type $\I_\tau$, and such that the infinite tensor product $e=\otimes_\pp e_\pp$ is a nonzero projection in $M_1$.\\

Let $\tau\in\NN$ be such that $\tau>1/\beta$. For any finite place $\pp$ of $k$, let $e_\pp=1-\mu_\pp^\tau\mu_\pp^{*\tau}\in M_{1,\pp}$.\\

Recall from equation (\ref{M1p}) that $M_{1,\pp}$ is naturally identified with $B\ell^2(\NN)$. Under this identification, the projection $e_\pp$ is the diagonal matrix whose $\tau$ first diagonal entries are $1$, and whose other entries are $0$. Thus, the reduced subfactor $M_{1,\pp,e_\pp}$ is of type $\I_\tau$. Note that
\begin{equation}
\label{xr63}
\widetilde\varphi_{\beta,\pp}(e_\pp)=1-\N\pp^{-\tau\beta}.
\end{equation}
Any decreasing sequence of projections in a von Neumann algebra converges weakly to a projection, so we can define a projection $e\in M_1$ by
$$e=\prod_\pp e_\pp=\bigotimes_\pp e_\pp.$$
By definition of $\tau$, we have $\tau\beta>1$, so
$$\widetilde\varphi_\beta(e)=\prod_\pp\widetilde\varphi_{\beta,\pp}(e_\pp)=\prod_\pp\left(1-\N\pp^{-\tau\beta}\right)=\frac{1}{\zeta_{k,\infty}(\tau\beta)}\neq 0.$$
In particular, $e\neq 0$. Let us define a state \index{phi_var_beta_tilde_e@$\widetilde\varphi_{\beta,e}$}$\widetilde\varphi_{\beta,e}$ on $M_{1,e}$ by
$$\forall x\in M_{1,e},\;\widetilde\varphi_{\beta,e}(x)
=\zeta_{k,\infty}(\tau\beta)\,\widetilde\varphi_\beta(x)=\frac{\widetilde\varphi_\beta(x)}{\widetilde\varphi_\beta(e)}\cdot$$
For any $\pp$, let us define a state \index{phi_var_beta_tilde_p_e_p@$\widetilde\varphi_{\beta,\pp,e_\pp}$}$\widetilde\varphi_{\beta,\pp,e_\pp}$ on $M_{1,\pp,e_\pp}$ by
$$\forall x\in M_{1,\pp,e_\pp},\;\widetilde\varphi_{\beta,\pp,e_\pp}(x) =\left(1-\N\pp^{-\tau\beta}\right)^{-1}\widetilde\varphi_{\beta,\pp}(x)=\frac{\widetilde\varphi_{\beta,\pp}(x)}{\widetilde\varphi_{\beta,\pp}(e_\pp)}\cdot$$
Coming back to the definition of the infinite tensor product (\cite{Tak}, III, Chapter XIV, \S 1), and using the expression (\ref{itpfi}) of $M_1$ as the infinite tensor product of the $(M_{1,\pp},\widetilde\varphi_{\beta,\pp})$, one can check that
\begin{equation}
\label{itpfired}
\left(M_{1,e},\;\widetilde\varphi_{\beta,e}\right)\;=\;\bigotimes_{\pp} \left(M_{1,\pp,e_\pp},\;\widetilde\varphi_{\beta,\pp,e_\pp}\right).
\end{equation}

Let $(\sigma^{\widetilde\varphi_\beta}_t)$ denote the modular flow of $\widetilde\varphi_\beta$. As $\widetilde\varphi_\beta$ is KMS$_\beta$ for the flow $\left(\widetilde\sigma_t\right)$, we have
\begin{equation}
\label{xr64}
\forall t\in\RR,\;\;\;\sigma^{\widetilde\varphi_\beta}_t=\widetilde\sigma_{\beta t}.
\end{equation}

\begin{lemma}
\label{xr65}
The factor $M_{1}$ is of type $\III_{q^{-\beta}}$.
\end{lemma}
\begin{proof}
Let us first prove that $q^{-\beta}$ belongs to the asymptotic ratio set $r_\infty(M_{1,e})$. We want to apply the criterion given on page 465 of \cite{Con94} to the ITPFI in equation (\ref{itpfired}). The eigenvalue list of $\widetilde\varphi_{\beta,\pp,e_\pp}$ is $(\lambda_{\pp,a})_{a=0,\ldots,\tau-1}$, where
$$\lambda_{\pp,a}=\frac{\left(1-\N\pp^{-\beta}\right)\N\pp^{-a\beta}}{1-\N\pp^{-\tau\beta}}\cdot$$
Let $r$ be such that $0 < r < 1$. For any $n\in\NN$, let
$$r(n)=\lfloor rq^n/n \rfloor.$$
By equation (\ref{riemann3}), there exists a $n_0\geqslant 1$ such that for any $n\geqslant n_0$, there exist (at least) $r(n)$ distinct finite places $$\pp_n^1,\ldots,\pp_n^{r(n)}$$ of $k$ such that $$\forall i\in\{1,\ldots,r(n)\},\;\;\;\N\pp_n^i=q^n.$$
For any $n\geqslant n_0$, let $I_n$ be the following set of places of $k$:
$$I_n=\{\pp_{2n}^1,\ldots,\pp_{2n}^{r(2n)},\pp_{2n+1}^1,\ldots,\pp_{2n+1}^{r(2n)}\}.$$
Let $X(I_n)=\{0,\ldots,\tau-1\}^{I_n}$ be the set of all applications from $I_n$ to $\{0,\ldots,\tau-1\}$. Define a measure $\lambda$ on $X(I_n)$ by
$$\lambda(\{f\})=\prod_{i=1}^{r(2n)} \lambda_{\pp_{2n}^i,f(\pp_{2n}^i)}\;\lambda_{\pp_{2n+1}^i,f(\pp_{2n+1}^i)}.$$
For any $i\in\{1,\ldots,r(2n)\}$, define elements $k^{1,i}_n$ and $k^{2,i}_n$ of $X(I_n)$ by
$$\forall \pp\in I_n,\;\;k^{1,i}_n(\pp)=1_{\pp=\pp_{2n}^i}
\;\;\;\text{and}\;\;\;k^{2,i}_n(\pp)=1_{\pp=\pp_{2n+1}^i}.$$
Put $K^1_n=\{k^{1,1}_n,\ldots,k^{1,r(2n)}_n\}$ and $K^2_n=\{k^{2,1}_n,\ldots,k^{2,r(2n)}_n\}$. For any $i\in\{1,\ldots,r(2n)\}$, we have $\lambda(\{k^{1,i}_n\})=\lambda(\{k^{1,1}_n\})$, so
\begin{eqnarray*}
\lambda(K^1_n) & = & r(2n)\, \lambda(\{k^{1,1}_n\}) \\
 & = & r\,\cdot \left(\frac{1-q^{-2n\beta}}{1-q^{-2n\tau\beta}}\right)^{\lfloor rq^{2n}/(2n) \rfloor}
\left(\frac{1-q^{-(2n+1)\beta}}{1-q^{-(2n+1)\tau\beta}}\right)^{\lfloor rq^{2n}/(2n) \rfloor}\left(q^{1-\beta}\right)^{2n}/(2n).
\end{eqnarray*}
Since $\beta\leqslant 1$, one checks easily that
\begin{equation}
\label{xr66}
\sum_{n\geqslant n_0} \lambda(K^1_n) = \infty.
\end{equation}
Let $\phi_n\::\:K^1_n\rightarrow K^2_n$ be the bijection defined by $\phi_n(k^{1,i}_n)=k^{2,i}_n$. For any $i\in\{1,\ldots,2n\}$, we have
$$\lambda\left(\{\phi_n(k^{1,i}_n)\}\right)/\lambda\left(\{k^{1,i}_n\}\right)=\frac{\lambda_{\pp_{2n}^i,0}\lambda_{\pp_{2n+1}^i,1}}
{\lambda_{\pp_{2n}^i,1}\lambda_{\pp_{2n+1}^i,0}}=\frac{q^{-(2n+1)\beta}}{q^{-2n\beta}}=q^{-\beta}.$$
Together with equation (\ref{xr66}), this allows to apply the criterion given on page 465 of \cite{Con94}, and we get
$$q^{-\beta}\in r_\infty(M_{1,e}).$$
Hence, by \cite{Con73}, Théorème 3.6.1, we have $q^{-\beta}\in S(M_{1,e})$. Hence, by \cite{Con73}, Corollaire 3.2.8 (b), we have $$q^{-\beta}\in S(M_1).$$
In particular, this shows that $S(M_1)\neq\{0,1\}$, so, by \cite{Con73}, Théorème 3.4.1, one gets that $S(M_1)\cap\RR^*_+$ is the orthogonal of $T(M_1)$ for the duality $(s,t)\mapsto s^{it}$. By construction, $\widetilde\sigma_{2\pi/\log q}=1$, so equation (\ref{xr64}) gives
\begin{equation}
\label{xr67}
\sigma^{\widetilde\varphi_\beta}_{2\pi/(\beta\log q)}=1.
\end{equation}
Thus,
$$T(M_1)\supset\frac{2\pi}{\beta\log q}\,\ZZ.$$
Hence, by orthogonality, we get
$$S(M_1)\cap\RR^*_+ \subset q^{\beta\ZZ}.$$
As we already know that $q^{-\beta}\in S(M_1)$, we get
$$S(M_1)\cap\RR^*_+ = q^{\beta\ZZ}.$$
Thus, $M_1$ is of type $\III_{q^{-\beta}}$.
\end{proof}

Let $M_{1,\widetilde\varphi_\beta}$ denote the centralizer of $\widetilde\varphi_\beta$ in $M_1$. We will only use Lemma \ref{xr65} through the following corollary:

\begin{coro}
\label{xr68}
The centralizer $M_{1,\widetilde\varphi_\beta}$ is a factor, of type $\II_1$.
\end{coro}
\begin{proof}
Lemma \ref{xr65} and equation (\ref{xr67}) allow to apply \cite{Con73}, Théorème 4.2.6, and we obtain that $M_{1,\widetilde\varphi_\beta}$ is a factor. Note that $\widetilde\varphi_\beta$ is a finite faithful normal trace on $M_{1,\widetilde\varphi_\beta}$. Hence, the type of $M_{1,\widetilde\varphi_\beta}$ can only be either $\II_1$ or $\I_n$ with $n\in\NN^*$. Let $\pp$ be a finite place of $k$. For any $n\geqslant 1$, set $x_n=\mu_\pp^n\mu_\pp^{*n}$. Note that the $x_n$ are fixed by the flow $(\sigma_t)$, hence by equation (\ref{xr64}) they belong to $M_{1,\widetilde\varphi_\beta}$. Equation (\ref{xr20}) shows that the $x_n$ are linearly independent over $\CC$. Thus, $M_{1,\widetilde\varphi_\beta}$ is infinite-dimensional over $\CC$, so its type cannot be $\I_n$ with $n\in\NN^*$. Hence, it must be $\II_1$.
\end{proof}

Our next goal is to prove (Lemma \ref{xr76}) that the centralizer $M_{\widetilde\varphi_\beta}$ of $\widetilde\varphi_\beta$ in $M$ is also a factor.

\begin{defi}For any $\ddd\in\Ig_\OO$, let \index{M_d@$M[\ddd]$ ($\ddd\in\Ig_\OO$)} $M[\ddd]$ denote the weak closure of $\HH[\ddd]$ in $M$.
\end{defi}

\begin{lemma}
\label{xr69}
Let $\ddd\in\Ig_\OO$. Let $\pp$ be a maximal ideal of $\Ig_\OO$ not dividing $\ddd$. Let $\sigma_\pp=(\pp,K_\ddd/k)\in\Gal(K_\ddd/k)$ be the Artin automorphism of $K_\ddd$ associated to $\pp$. Then:
\begin{mylist}
\item The automorphism $\sigma_\pp$ of $\HH[\ddd]$ extends uniquely to an ultraweakly continuous automorphism of $M[\ddd]$.
\item For all $x\in M[\ddd]$, we have
\begin{equation}
 \label{xr70}
x\mu_\pp=\mu_\pp\sigma_\pp(x).
\end{equation}
\end{mylist}
\end{lemma}
\begin{proof}
Let us first prove (1). Unicity is clear because, by the von Neumann density theorem, $\HH[\ddd]$ is ultraweakly dense in $M[\ddd]$. Let $\sigma\in\Gal(K/k)$ be such that $\sigma\vert_{K_\ddd}=\sigma_\pp$. As $\varphi_\beta\circ\sigma=\varphi_\beta$ on $C_{k,\infty}$, we know that $\sigma$ extends to an ultraweakly continuous automorphism of $M$, which we still note $\sigma$. The required extension of $\sigma_\pp$ is then obtained by taking the restriction of $\sigma$ to $M[\ddd]$.\\

Let us now check (2). By density, it is enough to check equation (\ref{xr70}) when $x\in\HH[\ddd]$. It then follows from Lemma \ref{xr26}.
\end{proof}

\begin{lemma}
\label{xr71}
Let $\ddd\in\Ig_\OO$. Let $M[\ddd]_{\widetilde\varphi_\beta}$ denote the centralizer of $\widetilde\varphi_\beta$ in $M[\ddd]$. Let $Z(M[\ddd]_{\widetilde\varphi_\beta})$ denote the center of $M[\ddd]_{\widetilde\varphi_\beta}$. Then:
$$Z(M[\ddd]_{\widetilde\varphi_\beta})\subset M_1.$$
\end{lemma}

{\bf Erratum: the statement and proof of this Lemma are wrong, as shown to us by Neshveyev.}

\begin{proof}
Let $x$ belong to $Z(M[\ddd]_{\widetilde\varphi_\beta})$. As $x$ belongs to $M[\ddd]$, it is fixed by $\Gal(K/K_\ddd)$. Let $\sigma\in\Gal(K_\ddd/k)=\Gal(K/k)/\Gal(K/K_\ddd)$. By Corollary \ref{xr2} ({\bf Erratum:} that corollary is wrong, as mentioned in the Erratum next to it), there exist finite places $\pp,\qq$ of $k$ not dividing $\ddd$, such that $\N\pp=\N\qq$, $\sigma_{\pp}=\sigma$ and $\sigma_\qq=1$. As $\N\pp=\N\qq$, we have
$$\forall t\in\RR,\;\;\sigma_t(\mu_\pp\mu_\qq^*)=\N\pp^{it}\N\qq^{-it}\mu_\pp\mu_\qq^*=\mu_\pp\mu_\qq^*.$$
Hence, by equation (\ref{xr64}), $\mu_\pp\mu_\qq^*\in M[\ddd]_{\widetilde\varphi_\beta}$. Thus, as $x$ belongs to the center of $M[\ddd]_{\widetilde\varphi_\beta}$, we have
\begin{equation}
\label{xr72}
x\mu_\pp\mu_\qq^*=\mu_\pp\mu_\qq^*x
\end{equation}
({\bf Erratum: } The previous statement is wrong, as it incorrectly assumes that $\mu_\pp$ and $\mu_\qq$ belong to $M[\ddd]$). On the other hand, by Lemma \ref{xr69} (2), we have
\begin{equation}
\label{xr73}
x\mu_\pp\mu_\qq^* = \mu_\pp\sigma_\pp(x)\mu_\qq^*,
\end{equation}
and we also compute
\begin{eqnarray}
\mu_\pp\mu_\qq^*x & = & \mu_\pp(x^*\mu_\qq)^*\nonumber\\
& = & \mu_\pp(\mu_\qq \sigma_\qq (x^*))^* \;\;\;\text{by Lemma \ref{xr69} (2)}\nonumber\\
& = & \mu_\pp\sigma_\qq (x)\mu_\qq^*.\label{xr74}
\end{eqnarray}
Combining equations (\ref{xr72}), (\ref{xr73}), and (\ref{xr74}), we get
\begin{equation}
 \label{xr75}
\mu_\pp\sigma_\pp(x)\mu_\qq^*=\mu_\pp\sigma_\qq (x)\mu_\qq^*.
\end{equation}
Multiplying both sides of equation (\ref{xr75}) by $\mu_\pp^*$ on the left and by $\mu_\qq$ on the right, and applying relation ($a_1$) of Proposition \ref{xr21}, we get
$$\sigma_\pp(x)=\sigma_\qq (x).$$
As $\sigma_\pp=\sigma$ and $\sigma_\qq=1$, we get
$$\sigma(x)=x.$$
Thus, $x\in M_1$.
\end{proof}

Let $M_{\widetilde\varphi_\beta}$ denote the centralizer of $\widetilde\varphi_\beta$ in $M$. 

\begin{lemma}
\label{xr76}
The centralizer $M_{\widetilde\varphi_\beta}$ is a factor, of type $\II_1$.
\end{lemma}
\begin{proof}
Note that $\widetilde\varphi_\beta$ is a finite, faithful, normal, positive, normalized trace on $M_{\widetilde\varphi_\beta}$. Let $\tr$ be another such trace on $M_{\widetilde\varphi_\beta}$.  Let us prove that $\tr=\widetilde\varphi_\beta$. Let $\ddd\in\Ig_\OO$.
By Connes' Radon-Nikod\'ym theorem, \cite{Con73}, Lemme 1.2.3 (b), there exists a positive element $h$ of $M[\ddd]_{\widetilde\varphi_\beta}$ such that
$$\forall x\in M[\ddd]_{\widetilde\varphi_\beta},\;\;\;\tr(x)=\widetilde\varphi_\beta(hx).$$
As $\widetilde\varphi_\beta$ and $\tr$ are faithful traces, one easily checks that $h$ belongs to the center $Z(M[\ddd]_{\widetilde\varphi_\beta})$. Thus, by Lemma \ref{xr71}, $h\in M_1$. Hence, the restriction of $\tr$ to $M[\ddd]_{\widetilde\varphi_\beta}$ is $\Gal(K/k)$-invariant, so
\begin{equation}
 \label{xr77}
\forall x\in M[\ddd]_{\widetilde\varphi_\beta},\;\;\;\tr(x) = \tr(\E(x)).
\end{equation}
As $(\sigma_t^{\widetilde\varphi_\beta})$ is $(2\pi/\log q)$-periodic, we have a normal conditional expectation
\begin{eqnarray*}
 E_{\widetilde\varphi_\beta}\;:\;M & \rightarrow & M_{\widetilde\varphi_\beta} \\
 x & \mapsto & \frac{\log q}{2\pi} \int_0^{2\pi/\log q} \sigma_t^{\widetilde\varphi_\beta}(x) dt.
\end{eqnarray*}
As $\HH$ is norm-dense in $C_{k,\infty}$ (see Proposition \ref{denseH}), it is ultraweakly dense in $M$, and it follows that $E_{\widetilde\varphi_\beta}(\HH)$ is ultraweakly dense in $E_{\widetilde\varphi_\beta}(M)$. We have
\begin{eqnarray*}
E_{\widetilde\varphi_\beta}(\HH) & = & E_{\widetilde\varphi_\beta}\left(\bigcup_{\ddd\in\Ig_\OO} \HH[\ddd]\right) \;\subset\; E_{\widetilde\varphi_\beta}\left(\bigcup_{\ddd\in\Ig_\OO} M[\ddd]\right)\\
& \subset & \bigcup_{\ddd\in\Ig_\OO} E_{\widetilde\varphi_\beta}\left(M[\ddd]\right) \;\subset\; \bigcup_{\ddd\in\Ig_\OO} M[\ddd]_{\widetilde\varphi_\beta}.
\end{eqnarray*}
Thus, $\bigcup_{\ddd\in\Ig_\OO} M[\ddd]_{\widetilde\varphi_\beta}$ is ultraweakly dense in $E_{\widetilde\varphi_\beta}(M)=M_{\widetilde\varphi_\beta}$. Thus, equation (\ref{xr77}) gives
\begin{equation}
\label{xr78}
\forall x\in M_{\widetilde\varphi_\beta},\;\;\;\tr(x) = \tr(\E(x)).
\end{equation}
We know by Corollary \ref{xr68} that $M_{1,\widetilde\varphi_\beta}$ is a type $\II_1$ factor. Hence, by Jones \cite{Jon03}, Corollary 7.1.19, we know that $\tr$ and $\widetilde\varphi_\beta$ agree on $M_{1,\widetilde\varphi_\beta}$. Thus, by equation (\ref{xr78}), we deduce that $\tr$ and $\widetilde\varphi_\beta$ agree on $M_{\widetilde\varphi_\beta}$. Hence, by \cite{Jon03}, Corollary 7.1.20, we deduce that $M_{\widetilde\varphi_\beta}$ is a factor, and the same argument that we made for $M_{1,\widetilde\varphi_\beta}$ shows that $M_{\widetilde\varphi_\beta}$ is also of type $\II_1$.
\end{proof}

\begin{coro}
\label{xr79}
We have $S(M)\neq\{0,1\}$. In other words, the factor $M$ is not of type $\III_0$.
\end{coro}
\begin{proof}
Suppose that $S(M)=\{0,1\}$. Then, by \cite{Con73}, Corollaire 3.2.7 (b), the center of $M_{\widetilde\varphi_\beta}$ has no minimal nonzero projection. Hence, by Lemma \ref{xr76}, one deduces that $\CC$ has no minimal nonzero projection, which is absurd.
\end{proof}

At last, we can prove the main result of this subsection. Recall that we have assumed $0<\beta\leqslant 1$.

\begin{theo}
\label{xr80}
The state $\varphi_\beta$ on $C_{k,\infty}$ is of type $\III_{q^{-\beta}}$. In other words, the factor $M$ is of type $\III_{q^{-\beta}}$.
\end{theo}

{\bf Erratum:} This is the wrong result mention in previous errata. It has been shown by Neshveyev and Rustad \cite{NesRus12} that the correct type is $\III_{q^{-\beta d_\infty}}$ where $d_\infty$ is the degree of the place $\infty$. The original statements have been kept for reference only.

\begin{proof}
By Corollary \ref{xr79} and \cite{Con73}, Théorème 3.4.1, the set $S(M)\cap\RR^*_+$ is the orthogonal of $T(M)$ for the duality $(s,t)\mapsto s^{it}$. Hence, it is enough to prove that
$$T(M)=\frac{2\pi}{\beta\log q}\cdot\ZZ.$$
As $\widetilde\sigma_{2\pi/\log q}=1$, equation (\ref{xr64}) gives
$$\frac{2\pi}{\beta\log q}\in T(M),$$
which proves one inclusion. Let us prove the other one. Let $t_0\in\RR$ be such that $t_0/\beta\in T(M)$. Thus, by equation (\ref{xr64}), $\widetilde\sigma_{t_0}$ is an inner automorphism of $M$. Let $u$ be an unitary of $M$ such that
$$\forall x\in M,\;\;\;\widetilde\sigma_{t_0}(x)=uxu^*.$$
For any $t\in\RR$ and $x\in M$, we have
$$\widetilde\sigma_t(u)\widetilde\sigma_t(x)\widetilde\sigma_t(u)^*=\widetilde\sigma_{t_0+t}(uxu^*)=u\widetilde\sigma_t(x)u^*,$$
so the unitaries $u$ and $\sigma_t(u)$ implement the same inner automorphism of the factor $M$, so there exists some $z_t\in\CC$ with $\abs{z_t}=1$ and $\sigma_t(u)=z_t u$. The map $t\mapsto z_t$ is a character of $\RR$, so there exists $\theta\in\RR$ such that
$$\forall t\in\RR,\;\;\;z_t=e^{i\theta t}.$$
The KMS$_\beta$ property of the state $\widetilde\varphi_\beta$ for the flow $(\widetilde\sigma_t)$ applied to the pair $(u^*,u)$ gives a bounded continuous function $F$ on the strip $0\leqslant \im z \leqslant \beta$, holomorphic on the interior of the strip, such that
$$\forall t\in\RR,\;\;\;F(t)=\widetilde\varphi_\beta(u^*\sigma_t(u))\;\;\;\text{and}\;\;\;F(t+i\beta)=\widetilde\varphi_\beta(\sigma_t(u)u^*).$$
Thus,
\begin{equation}
\label{xr81}
\forall t\in\RR,\;\;\;F(t)=e^{i\theta t}=F(t+i\beta).
\end{equation}
Hence $F$ is the holomorphic function $z\mapsto e^{i\theta z}$ and, evaluating equation (\ref{xr81}) at $t=0$, one gets
$$e^{-\theta\beta}=1.$$
Thus, $\theta=0$, so $u$ is fixed by the flow $(\widetilde\sigma_t)$. Hence, by equation (\ref{xr64}), the unitary $u$ belongs to the centralizer $M_{\widetilde\varphi_\beta}$ of $\widetilde\varphi_\beta$. Moreover, by equation (\ref{xr64}), any element of $M_{\widetilde\varphi_\beta}$ is fixed by the flow $(\widetilde\sigma_t)$, and so commutes with $u$, by definition of $u$. Hence, $u$ belongs to the center of $M_{\widetilde\varphi_\beta}$. Thus, by Lemma \ref{xr76}, one deduces that $u\in\CC$, so, as an automorphism of $M$,
\begin{equation}
 \label{xr82}
\widetilde\sigma_{t_0}=1.
\end{equation}
By equation (\ref{riemann3}), for any sufficiently large $n$, there exist finite places $\pp$ and $\qq$ of $k$ such that $\N\pp=q^n$ and $\N\qq=q^{n+1}$. We then have
$\sigma_{t_0}(\mu_\qq\mu_\pp^*)=q^{(n+1)it_0-nit_0}\mu_\qq\mu_\pp^*=q^{it_0}\mu_\qq\mu_\pp^*$. On the other hand, equation (\ref{xr82}) gives $\sigma_{t_0}(\mu_\qq\mu_\pp^*)=\mu_\qq\mu_\pp^*.$ Thus, we get $1=q^{it_0}$, so $t_0\in 2\pi/(\log q)\ZZ$, which completes the proof.
\end{proof}

\newpage
\addcontentsline{toc}{section}{References}

\bibliographystyle{plain}
\bibliography{BC-for-function-fields}

\vskip 1cm

\noindent Benoît Jacob\\
\noindent jacob.benoit.1@gmail.com

\end{document}